\documentclass[letterpaper,11pt]{article}
\usepackage[margin=1in]{geometry}
\usepackage[utf8]{inputenc}
\usepackage{amsmath}
\usepackage{url}
\usepackage{amsthm}
\usepackage{amsfonts}
\usepackage{float}
\usepackage{graphicx}
\usepackage{bm}
\usepackage{subcaption}
\usepackage{amssymb}
\usepackage{bbm}
\usepackage{color}
\usepackage{algorithm}
\usepackage{algpseudocode}
\usepackage{booktabs}
\usepackage{multirow}

\newcommand{\eat}[1]{}

\renewcommand{\epsilon}{\varepsilon}

\definecolor{cof}{RGB}{219,144,71}
\definecolor{pur}{RGB}{186,146,162}
\definecolor{greeo}{RGB}{91,173,69}
\definecolor{greet}{RGB}{52,111,72}

\def \bP {\mathbb{P}}
\def \bE {\mathbb{E}}
\def \bR {\mathbb{R}}

\def \var {\mathsf{Var}}
\def\1{\mathbbm{1}}

\theoremstyle{plain}
\newtheorem{theorem}{Theorem}[section]
\newtheorem{lemma}[theorem]{Lemma}
\newtheorem{example}[theorem]{Example}
\newtheorem{remark}[theorem]{Remark}
\newtheorem{corollary}[theorem]{Corollary}
\newtheorem{definition}[theorem]{Definition}

\newtheorem{assumption}{Assumption}

\usepackage[
            CJKbookmarks=true,
            bookmarksnumbered=true,
            bookmarksopen=true,
            colorlinks=true,
            citecolor=red,
            linkcolor=blue,
            anchorcolor=red,
            urlcolor=blue
            ]{hyperref}

\usepackage{xspace}

\newcommand{\stepa}[1]{\overset{\rm (a)}{#1}}
\newcommand{\stepb}[1]{\overset{\rm (b)}{#1}}
\newcommand{\stepc}[1]{\overset{\rm (c)}{#1}}
\newcommand{\stepd}[1]{\overset{\rm (d)}{#1}}
\newcommand{\stepe}[1]{\overset{\rm (e)}{#1}}

\newcommand{\naturals}{\mathbb{N}}

\newcommand{\Poi}{\mathsf{Poi}}

\definecolor{myblue}{rgb}{.8, .8, 1}
\definecolor{mathblue}{rgb}{0.2472, 0.24, 0.6} 
\definecolor{mathred}{rgb}{0.6, 0.24, 0.442893}
\definecolor{mathyellow}{rgb}{0.6, 0.547014, 0.24}

\newcommand{\calA}{{\mathcal{A}}}

\newcommand{\calC}{{\mathcal{C}}}

\newcommand{\calK}{{\mathcal{K}}}
\newcommand{\calL}{{\mathcal{L}}}
\newcommand{\calM}{{\mathcal{M}}}
\newcommand{\calN}{{\mathcal{N}}}

\newcommand{\calP}{{\mathcal{P}}}

\newcommand{\calT}{{\mathcal{T}}}

\newcommand{\calW}{{\mathcal{W}}}
\newcommand{\calX}{{\mathcal{X}}}
\newcommand{\calY}{{\mathcal{Y}}}

\usepackage{cleveref}
\crefname{lemma}{Lemma}{Lemmas}
\Crefname{lemma}{Lemma}{Lemmas}
\crefname{thm}{Theorem}{Theorems}
\Crefname{thm}{Theorem}{Theorems}
\crefformat{equation}{(#2#1#3)}

\begin{document}

\title{On the Statistical Complexity of Sample Amplification}
\author{Brian Axelrod, Shivam Garg, Yanjun Han, Vatsal Sharan, and Gregory Valiant\thanks{B. Axelrod and G. Valiant are with the Department of Computer Science, Stanford University, emails: \url{{baxelrod,valiant}@cs.stanford.edu}. S. Garg is with the Microsoft Research, email: \url{shigarg@microsoft.com}. Y. Han is with the Courant Institute of Mathematical Sciences and the Center for Data Science, New York University, email: \url{yanjunhan@nyu.edu}. V. Shaman is with the Department of Computer Science, University of Southern California, email:  \url{vsharan@usc.edu}. }}

\maketitle

\begin{abstract}
The ``sample amplification'' problem formalizes the following question: Given $n$ i.i.d. samples drawn from an unknown distribution $P$, when is it possible to produce a larger set of $n+m$ samples which cannot be distinguished from $n+m$ i.i.d. samples drawn from $P$? In this work, we provide a firm statistical foundation for this problem by deriving generally applicable amplification procedures, lower bound techniques and connections to existing statistical notions. Our techniques apply to a large class of distributions including the exponential family, and establish a rigorous connection between sample amplification and distribution learning.
\end{abstract}

\section{Introduction}

Consider the following problem of manufacturing more data: an amplifier is given $n$ samples drawn i.i.d. from an unknown distribution $P$, and the goal is to generate a larger set of $n+m$ samples which are indistinguishable from $n+m$ i.i.d. samples from $P$. How large can $m$ be as a function of $n$ and the distribution class in question? Are there sound and systematic ways to generate a larger set of samples? Is this task strictly \emph{easier} than the learning task, in the sense that the number of samples required for generating $n+1$ samples is smaller than that required for learning $P$?

In our preliminary work \cite{axelrod2019sample}, we formalized this problem as the \emph{sample amplification} problem, considering total variation (TV) as the measure for indistinguishability.

\begin{definition}[Sample Amplification]\label{def:sample_amplification}
Let $\mathcal{P}$ be a class of probability distributions over a domain $\mathcal{X}$. We say that $\mathcal{P}$ admits an $(n,n+m,\varepsilon)$ \emph{sample amplification procedure} if there exists a (possibly randomized) map $T_{\mathcal{P},n,m,\varepsilon} : \mathcal{X}^n \rightarrow \mathcal{X}^{n+m}$ such that
\begin{align}\label{eq:target_sample_amplification}
\sup_{P\in \mathcal{P}}\|P^{\otimes n}\circ T_{\mathcal{P},n,m,\varepsilon}^{-1} - P^{\otimes (n+m)} \|_{\text{\rm TV}} \le \varepsilon. 
\end{align}
\end{definition} 

An equivalent formulation to view Definition \ref{def:sample_amplification} is as a game between two parties: an amplifier and a verifier. The amplifier gets $n$ samples drawn i.i.d. from the unknown distribution $P$ in the class $\mathcal{P}$, and her goal is to generate a larger dataset of $n+m$ samples which must be accepted with good probability by any verifier that also accepts a real dataset of $n+m$ i.i.d. samples from $P$ with good probability. Here, the verifier is computationally unbounded and knows the distribution $P$, but does not observe the amplifier's original set of $n$ samples. 

Along with being a natural statistical task, the sample amplification framework is also
relevant from a practical standpoint. Currently, there is an enormous trend in the machine learning community to train models on datasets that have been enlarged in various ways.  There are relatively transparent and classical approaches to achieve this, such as leveraging known invariances such as rotation or translation invariance to augment the dataset by including transformed versions of each original datapoint \cite{simard2003best,krizhevsky2012imagenet,cubuk2019autoaugment,shorten2019survey,cubuk2020randaugment}.
More recently, deep generative models have been used to both directly enlarge training data and construct larger datasets consisting of samples with properties that are rare in naturally occurring datasets \cite{antoniou2017data, calimeri2017biomedical, frid2018gan,madani2018chest,yi2018data,sandfort2019data,chatziagapi2019data,yi2019generative,han2019combining,chlap2021review, luzi2022boomerang, chen2022deep, lu2023machine, bai2022constitutional}.
More opaque approaches such as MixUp \cite{zhang2018mixup} and related techniques \cite{tokozume2018between,yun2019cutmix,verma2019manifold,hendrycks2019augmix} which add a significant fraction of new datapoints that are explicitly \emph{not} supported in the true data distribution are also very popular since they seem to improve the performance of the trained models. Given this current zoo of widely implemented approaches to enlarging datasets, there is a clear motivation for bringing a more principled statistical understanding to such approaches.  One natural starting point is the statistical setting we consider that asks the extent to which datasets can be enlarged in a perfect sense---where it is not possible to distinguish the enlarged dataset from a set of i.i.d. draws. Moreover, this work lays a foundation for the ambitious broader goal of understanding how various amplification techniques interact with downstream learning algorithms and statistical estimators, and developing amplification techniques that are optimal for certain classes of such algorithms and estimators.

In \cite{axelrod2019sample}, a subset of the authors introduced the  sample amplification problem, and studied two classes of distributions: the Gaussian location model and discrete distribution model. For these examples, they characterized the statistical complexity of sample amplification   and showed that it is strictly smaller than that of learning. In this paper, we work towards a general understanding of the statistical complexity of the sample amplification problem, and its relationship with learning. The main contributions of this paper are as follows: 

\begin{table}[]
\begin{tabular}{lll}
\toprule
\textbf{Distribution Class}                                    & \textbf{Amplification}              & \textbf{Procedure} \\ \midrule
Gaussian with unknown mean and fixed covariance       & \multirow{2}{*}{$(n, n+\Theta(n\varepsilon/\sqrt{d}))$}              & \multirow{2}{*}{Sufficiency/Learning}    \\
(UB: Example \ref{example:mean}, \ref{example:mean_computation}, \ref{example:mean_shuffle}; LB: Theorem \ref{thm:Gaussian}, \ref{thm:lower_product})\smallskip & & \\
Gaussian with unknown mean and covariance             & \multirow{2}{*}{$(n, n+\Theta(n\varepsilon/d))$}                     & \multirow{2}{*}{Sufficiency}             \\
(UB: Example \ref{example:gaussian_covariance}, \ref{example:mean_covariance}; LB: Example \ref{example:lower_covariance}) \smallskip & & \\
Gaussian with $s$-sparse mean and identity covariance & \multirow{2}{*}{$(n, n+\Theta(n\varepsilon/\sqrt{s \log d}))$} & \multirow{2}{*}{Learning}                \\
(UB: Example \ref{example:shuffle_sparse_gaussian}; LB: Example \ref{example:lower_sparse_gaussian}) \smallskip & & \\
Discrete distributions with support size at most $k$  & \multirow{2}{*}{$(n, n+\Theta(n\varepsilon/\sqrt{k}))$}              & \multirow{2}{*}{Learning}              \\
(UB: Example \ref{example:shuffle_discrete}; LB: \cite[Theorem 1]{axelrod2019sample}) \smallskip & & \\
Poissonized discrete distributions with support at most $k$ & \multirow{2}{*}{($n$, $n + \Theta(\sqrt{n}\varepsilon
+ n\varepsilon/\sqrt{k}))$} & \multirow{2}{*}{Learning} \\
(UB: Example \ref{example:discrete_poisson}; LB: Example \ref{example:discrete_poisson}) \smallskip & & \\
$d$-dim. product of Exponential distributions         & \multirow{2}{*}{$(n, n+\Theta(n\varepsilon/\sqrt{d}))$}              & \multirow{2}{*}{Sufficiency/Learning}    \\
(UB: Example \ref{example:exponential}, \ref{example:shuffle_exponential}; LB: Theorem \ref{thm:lower_product}) \smallskip & & \\
Uniform distribution on $d$-dim. rectangle            & \multirow{2}{*}{$(n, n+\Theta(n\varepsilon/\sqrt{d}))$}              & \multirow{2}{*}{Sufficiency/Learning}   \\
(UB: Example \ref{example:uniform}, \ref{example:shuffle_uniform}; LB: Theorem \ref{thm:lower_product}) \smallskip & & \\
$d$-dim. product of Poisson distributions             & \multirow{2}{*}{$(n, n+\Theta(n\varepsilon/\sqrt{d}))$}              & \multirow{2}{*}{Sufficiency+Learning}   \\
(UB: Example \ref{example:poisson}; LB: Theorem \ref{thm:lower_product})  & & \\
\bottomrule
\end{tabular}
\caption{Sample amplification achieved by the presented procedures. Results include matching upper bounds (UB) and lower bounds (LB), with appropriate pointers to specific examples or theorems for details.}\label{tab:examples} 
\end{table}

\begin{enumerate}
    \item \textbf{Amplification via sufficiency.} Our first contribution is a simple yet powerful procedure for sample amplification, i.e. apply the sample amplification map only to sufficient statistics. Specifically, the learner computes a sufficient statistic $T_n$ from $X^n$, maps $T_n$ properly to some $T_{n+m}$, and generates new samples $\widehat{X}^{n+m}$ from some conditional distribution conditioned on $T_{n+m}$. Surprisingly, this simple idea leads to a much cleaner procedure than \cite{axelrod2019sample} under Gaussian location models which is also exactly optimal (cf. Theorem~\ref{thm:Gaussian}). The range of applicability also extends to general exponential families, with rate-optimal sample amplification performances. Specifically, for general $d$-dimensional exponential families with a mild moment condition, the sufficiency-based procedure achieves an $(n, n+O(n\varepsilon/\sqrt{d}), \varepsilon)$ sample amplification for large enough $n$, which by our matching lower bounds in Section \ref{sec:lower_bound} is asymptotically minimax rate-optimal. 
    
    \item \textbf{Amplification via learning.} Our second contribution is another general sample amplification procedure that does not require the existence of a sufficient statistic, and also sheds light on the relationship between learning and sample amplification. This procedure essentially draws new samples from a rate-optimal estimate of the true distribution, and outputs a random permutation of the old and new samples. The procedure achieves an $(n, n+O(\varepsilon\sqrt{n/r_{\chi^2}(\mathcal{P},n)}),\varepsilon) $ sample amplification, where $r_{\chi^2}(\mathcal{P},n) $ denotes the minimax risk for learning $P\in \mathcal{P}$ under the expected $\chi^2$ divergence given $n$ samples. This shows that learning $P$ to $\chi^2$ divergence $O(n/\varepsilon^2)$ is sufficient for non-trivial sample amplification. 
    
    In addition, we show that for the special case of product distributions, it is important that the permutation step be applied coordinatewise to achieve the optimal sample amplification. Specifically, if $\mathcal{P} = \prod_{j=1}^d \mathcal{P}_j$, this procedure achieves a better sample amplification 
    \begin{align*}
    \left(n, n+O\left(\varepsilon\sqrt{\frac{n}{\sum_{j=1}^d r_{\chi^2}(\mathcal{P}_j, n)}}\right), \varepsilon\right). 
    \end{align*}
    We have summarized several examples in Table \ref{tab:examples} where the sufficiency and/or learning based sample amplification procedures are optimal. Note that there is no golden rule for choosing one idea over the other, and there exists an example where the above two ideas must be combined. 
    \item \textbf{Minimax lower bounds.} Complementing our sample amplification procedures, we provide a general recipe for proving lower bounds for sample amplification. This recipe is intrinsically different from the standard techniques of proving lower bounds for hypothesis testing, for the sample amplification problem differs significantly from an estimation problem. In particular, specializing our recipe to product models shows that, for essentially \emph{all} $d$-dimensional product models, an $(n, n+Cn\varepsilon/\sqrt{d},\varepsilon)$ sample amplification is impossible for some absolute constant $C<\infty$ \emph{independent of} the product model. 
    
    For non-product models, the above powerful result does not directly apply, but proper applications of the general recipe could still lead to tight lower bounds for sample amplification. Specifically, we obtain matching lower bounds for all examples listed in Table \ref{tab:examples}, including the non-product examples.  
\end{enumerate}

We now provide several numerical simulations to suggest the potential utility of sample amplification. Recall that a practical motivation of sample amplification is to produce an enlarged dataset that can be fed into a \emph{distribution-agnostic} algorithm in downstream applications. Here, we consider the following basic experiments in that vein: 

\begin{figure}
	\centering
	\begin{subfigure}[b]{\textwidth}
        \centering
		\includegraphics[width=0.55\linewidth]{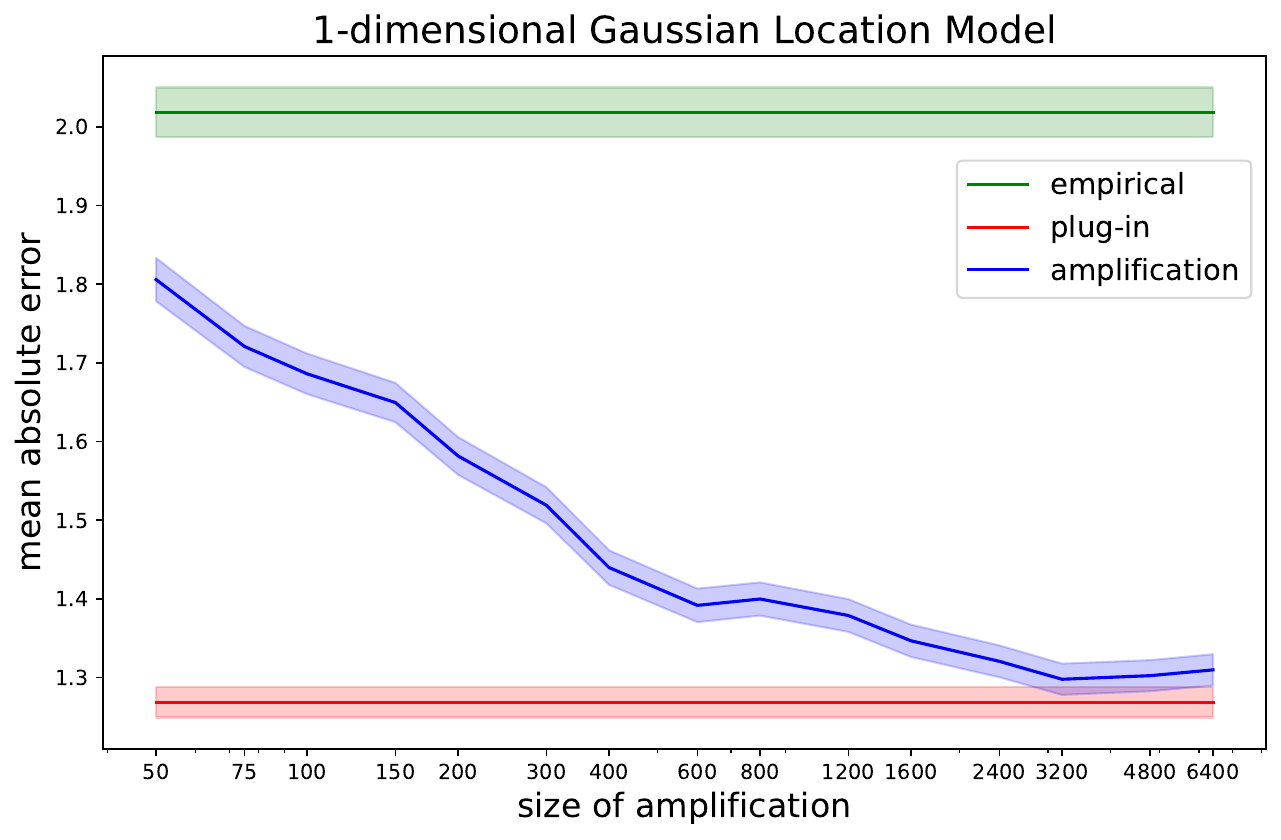}\caption{Mean absolute errors for estimation of the fourth moment $\bE[X^4]$ in a one-dimensional Gaussian location model $X_1,\cdots,X_n\sim \calN(\mu,1)$ with $n=100, \mu=1$.}\label{fig:a}
	\end{subfigure}
	\begin{subfigure}[b]{\textwidth}
         \centering
		\includegraphics[width=0.55\linewidth]{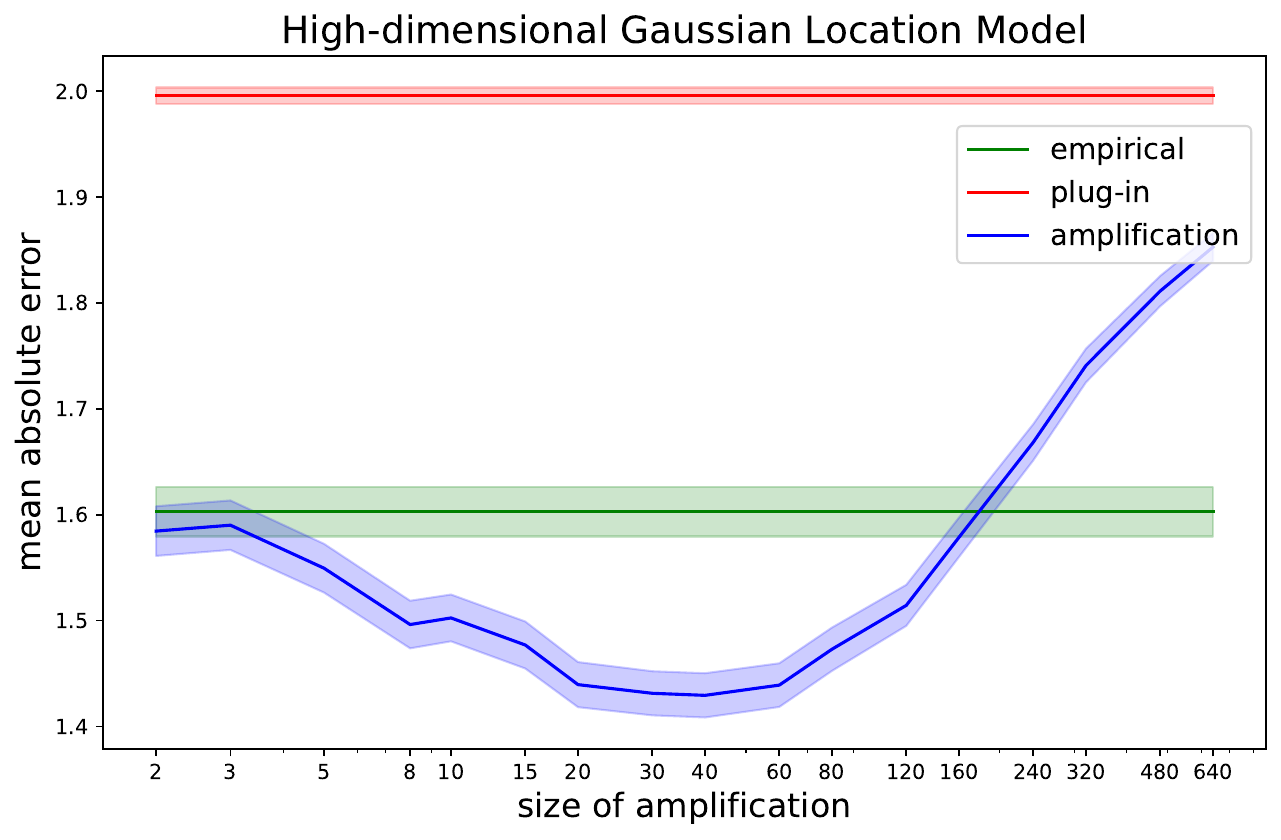}\caption{Mean absolute errors for estimation of the squared $L_2$ norm $\bE[\|X\|_2^2]$ in a $d$-dimensional Gaussian location model $X_1,\cdots,X_n\sim \calN(\mu,I_d)$ with $n=50, d=100, \mu={\bf 1}/\sqrt{d}$. }\label{fig:b}
	\end{subfigure}
	\begin{subfigure}[b]{\textwidth}
         \centering
		\includegraphics[width=0.55\linewidth]{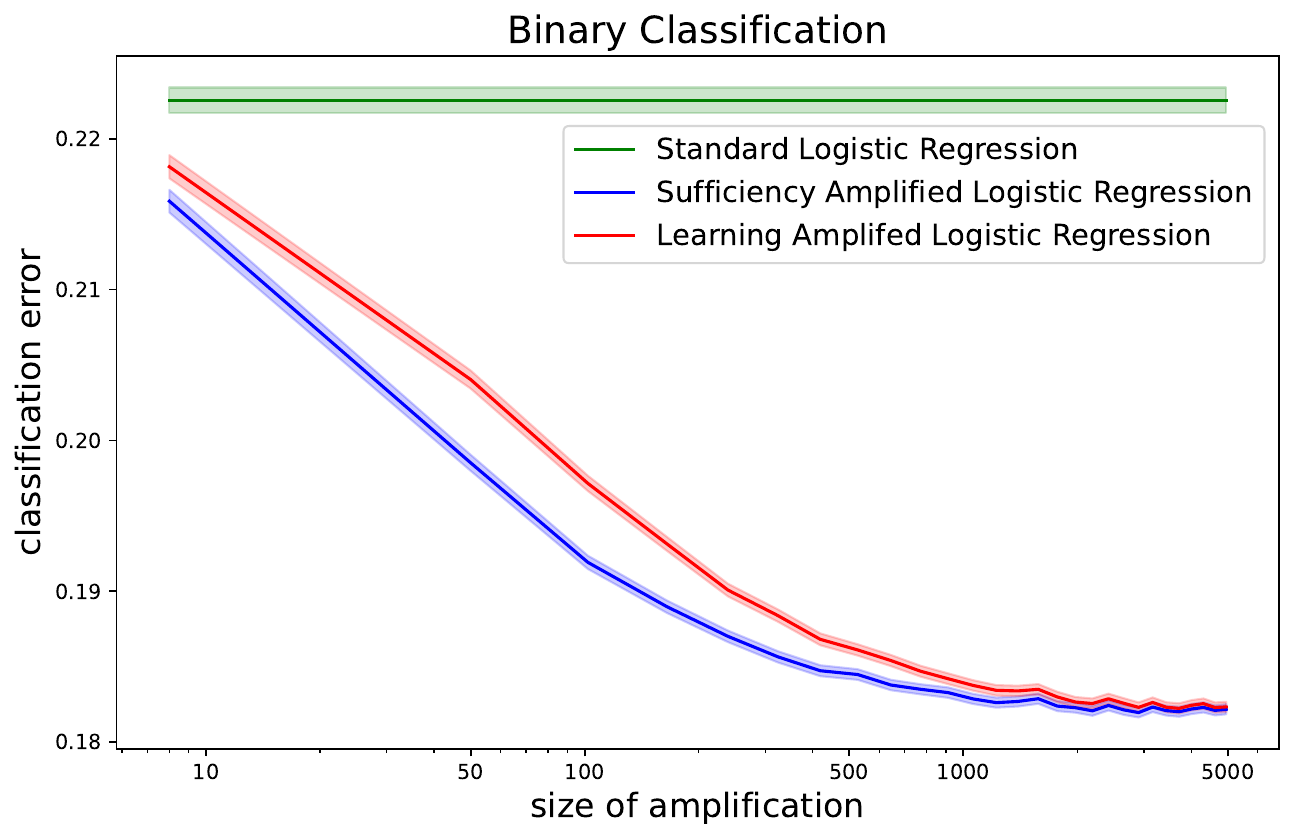}\caption{Classification errors for binary classification between two clusters of covariates $X_1,\cdots,X_{n/2}\sim \calN(e_1, I_d)$ and $X_{n/2+1},\cdots,X_n\sim \calN(-e_1,I_d)$, with $n=50$ and $d=10$.}\label{fig:c}
	\end{subfigure}
	\caption{Sample amplification experiments. The $x$-axis corresponds to the amount of amplification, $m$, and the shaded area depicts the 95\% confidence interval based on 5,000 simulations.}\label{fig:plots}
\end{figure}

\begin{itemize}
    \item Fourth moment estimation for one-dimensional Gaussian: here we observe $X_1,\cdots,X_n\sim \calN(\mu,1)$ with $n=100$ and $\mu=1$, and we consider three estimators. The empirical estimator operates in a distribution-agnostic fashion and is simply the empirical fourth moment $n^{-1}\sum_{i=1}^n X_i^4$. The plug-in estimator uses the knowledge of Gaussianity: it first estimates $\widehat{\mu} = \bar{X}$ and then uses $\bE_{X\sim \calN(\widehat{\mu},1)}[X^4] = \widehat{\mu}^4 + 6\widehat{\mu}^2 + 3$. The amplified estimator first amplifies the sample $X^n$ into $Y^{n+m}$ via sufficiency (cf. Example \ref{example:mean}), and then uses the empirical estimator $(n+m)^{-1}\sum_{j=1}^{n+m} Y_j^4$ based on the enlarged sample $Y^{n+m}$. The plots of the mean absolute errors (MAEs) are displayed in Figure \ref{fig:a}. We observe that although the empirical estimator based on the original sample $X^n$ has a large MAE, its performance is  improved based on the amplified sample $Y^{n+m}$. 
    \item Squared $L_2$ norm estimation for high-dimensional Gaussian: here we observe $X_1,\cdots,X_n\sim \calN(\mu,I_d)$ with $n=50, d=100$ and $\mu={\bf 1}/\sqrt{d}$, and we again consider three estimators for $\bE[\|X\|_2^2]$. As before, the empirical estimator is simply $n^{-1}\sum_{i=1}^n \|X_i\|_2^2$, and the plug-in estimator uses the knowledge $\bE_{X\sim \calN(\widehat{\mu},I_d)}[\|X\|_2^2] = \|\widehat{\mu}\|_2^2 + d$ and estimates $\widehat{\mu} = \bar{X}$. As for the amplified estimator, it first amplifies the sample $X^n$ into $Y^{n+m}$ via sufficiency (cf. Example \ref{example:mean}), and then uses the empirical estimator based on $Y^{n+m}$. The plots of the mean absolute errors are displayed in Figure \ref{fig:b}. Here the empirical estimator outperforms the plug-in estimator due to a smaller bias, while the sample amplification further reduces the MAE as long as the size of amplification $m$ is not too large. This could be explained by the bias-variance tradeoff, where the amplified estimator interpolates between the empirical estimator (with no bias) and the plug-in estimator (with the smallest asymptotic variance).
    \item Binary classification: here we observe two clusters of covariates $X_1,\cdots,X_{n/2}\sim \calN(e_1, I_d)$ (with label $1$) and $X_{n/2+1},\cdots,X_n\sim \calN(-e_1,I_d)$ (with label $-1$), with $n=50, d=10$ and $e_1$ being the first basis vector. The target is to train a classifier with a high classification accuracy on the test data with the same distribution. The standard classifier is via logistic regression, which does not use the knowledge of Gaussianity. To apply sample amplification, we first amplify the sample in each class via either sufficiency (cf. Example \ref{example:mean}) or learning (cf. Example \ref{example:mean_shuffle}), and then run logistic regression on the amplified datasets. Figure \ref{fig:c} displays the classification errors of all three approaches, and shows that both amplification procedures help reduce the classification error even for small values of $m$. 
\end{itemize}

The above experiments demonstrate the potential for  sample amplification to leverage knowledge of the data distribution to produce a larger dataset that is then fed into downstream distribution-agnostic algorithms. Some experiments (e.g. Figure \ref{fig:b}) also suggest a limit beyond which the amplification procedure alters the data distribution too much. We believe that rigorously examining amplification through the lens of the performance of downstream estimators and algorithms, including those illustrated in our numerical simulations, would be a fruitful direction for future work.

\section{Connections, limitations and future work}
As discussed above, it is commonplace in machine learning to increase the size of  datasets using various heuristics, often resulting in large gains in downstream learning performance. However, a clear statistical understanding of when this is possible and what techniques are useful for this is missing. A natural starting point to get a better understanding is the formulation we consider that asks the extent to which datasets can be amplified in a perfect sense---where any verifier who knows the true distribution is not able to distinguish the amplified dataset from a set of i.i.d. draws.

A limitation of the sample amplification formulation described above is that the additive amplification factor $m$ is rather small (e.g., $O(n\varepsilon/\sqrt{d})$ for $d$-dimensional exponential families). Moreover, we show matching lower bounds demonstrating that this factor cannot be improved even when $n$ is large enough to learn the distribution to non-trivial accuracy. However, it might be possible to achieve larger amplification factors with restricted verifiers, for instance, the class of verifiers corresponding to learning algorithms used for downstream tasks (see \cite{axelrod2019sample} for other possible classes of verifiers). Investigating the sample amplification problem with such restricted verifiers may be a practically fruitful future direction.

Despite this limitation, the sample amplification formulation does yield high-level insights that can inform the way datasets are amplified in practice. For instance, from the results in this paper, we know that sample amplification is possible for a broad class of distributions even when learning is not possible. 
Moreover, both our sufficiency or learning based approaches modify the original data points in general, conforming to the lower bound in \cite{axelrod2019sample} that optimal amplification may be impossible if the amplifier returns a superset of the input dataset. These observations show that the folklore way of enlarging datasets by learning the data distribution and adding more samples from the learned distribution can be far from optimal. 

\textbf{Connections with other statistical notions.}
An equivalent view of Definition \ref{def:sample_amplification} is through Le Cam's distance \cite{le1972limits}, a classical concept in statistics. The formal definition of Le Cam's distance $\Delta(\mathcal{M}, \mathcal{N})$ is summarized in Definition \ref{def:le_cam_distance}; roughly speaking, it measures the fundamental difference in power in the statistical models $\mathcal{M}$ and $\mathcal{N}$, without resorting to specific estimation procedures. The sample amplification problem is equivalent to the study of Le Cam's distance $\Delta(\mathcal{P}^{\otimes n}, \mathcal{P}^{\otimes (n+m)})$ between product models, where \eqref{eq:target_sample_amplification} is precisely equivalent to $\Delta(\mathcal{P}^{\otimes n}, \mathcal{P}^{\otimes (n+m)}) \le \varepsilon$. However, in the statistics literature, Le Cam's distance was mainly used to study the \emph{asymptotic} equivalence, where a typical target is to show that $\lim_{n\to\infty} \Delta(\calM_n, \calN_n)=0$ for certain sequences of statistical models. For example, showing that localized regular statistical models converge to Gaussian location models is the fundamental idea behind the H\'{a}jek--Le Cam asymptotic statistics; see \cite{le1972limits,LeCam1986asymptotic,le1990asymptotics} and \cite[Chapter 9]{Vandervaart2000}. In nonparametric statistics, there is also a rich line of research \cite{brown1996asymptotic,brown2002asymptotic,brown2004equivalence,ray2019asymptotic} establishing asymptotic (non-)equivalences, based on Le Cam's distance, between density models, regression models, and Gaussian white noise models. In the above lines of work, only asymptotic results were typically obtained with a fixed dimension and possibly slow convergence rate. In contrast, we aim to obtain a non-asymptotic characterization of $\Delta(\mathcal{P}^{\otimes n}, \mathcal{P}^{\otimes (n+m)})$ in $(n,m)$ and the dimension of the problem, a task which is largely underexplored in the literature. 

Another related angle is from reductions between statistical models. Over the past decade there has been a growing interest in constructing polynomial-time reductions between various statistical models (typically from the planted clique) to prove statistical-computational gaps, see, e.g. \cite{berthet2013optimal,ma2015computational,brennan2019optimal,brennan2020reducibility}. The sample amplification falls into the reduction framework, and aims to perform reductions from a product model $\mathcal{P}^{\otimes n}$ to another product model $\mathcal{P}^{\otimes (n+m)}$. While previous reduction techniques were mainly constructive and employed to prove computational lower bounds, in this paper we also develop general tools to prove limitations of all possible reductions purely from the statistical perspective. 

\textbf{Organization.}
The rest of this paper is organized as follows. Section \ref{sec:prelim} lists some notations and preliminaries for this paper, and in particular introduces the concept of Le Cam's distance. Section~\ref{sec.sufficiency} introduces a sufficiency-based procedure for sample amplification, with asymptotic properties for general exponential families and non-asymptotic performances in several specific examples. Section \ref{sec.shuffling} is devoted to a  learning-based procedure for sample amplification, with a general relationship between sample amplification and the $\chi^2$ estimation error, as well as its applications in several examples. Section \ref{sec:lower_bound} presents the general idea of establishing lower bounds for sample amplification, with a universal result specializing to product models. Section \ref{sec:discussion} discusses more examples in sample amplification and learning, and shows that these tasks are in general non-comparable. More concrete examples of both the upper and lower bounds, auxiliary lemmas, and proofs are relegated to the appendices.

\section{Preliminaries}\label{sec:prelim}

We use the following notations throughout this paper. For a random variable $X$, let $\calL(X)$ be the law (i.e. probability distribution) of $X$. For a probability distribution $P$ on a probability space $\Omega$ and a measurable map $T: \Omega \to \Omega'$, let $P\circ T^{-1}$ denotes the pushforward probability measure, i.e. $\calL(T(X))$ with $\calL(X)=P$. For a probability measure $P$, let $P^{\otimes n}$ be the $n$-fold product measure. For a positive integer $n$, let $[n]\triangleq \{1,\cdots,n\}$, and $x^n \triangleq (x_1,\cdots,x_n)$. We adopt the following asymptotic notations: for two non-negative sequences $(a_n)$ and $(b_n)$, we use $a_n = O(b_n)$ to denote that $\limsup_{n\to\infty} a_n/b_n < \infty$, and $a_n = \Omega(b_n)$ to denote $b_n = O(a_n)$, and $a_n = \Theta(b_n)$ to denote both $a_n = O(b_n)$ and $b_n = O(a_n)$. We also use the notations $O_c, \Omega_c, \Theta_c$ to denote the respective meanings with hidden constants depending on $c$. For probability measures $P,Q$ defined on the same probability space, the total variation (TV) distance, Hellinger distance, Kullback--Leibler (KL) divergence, and the chi-squared divergence are defined as follows: 
\begin{align*}
\|P-Q\|_{\text{TV}} &= \frac{1}{2}\int |\mathrm{d}P - \mathrm{d}Q|, \qquad \qquad H(P,Q) = \left( \frac{1}{2}\int (\sqrt{\mathrm{d}P} - \sqrt{\mathrm{d}Q})^2 \right)^{\frac{1}{2}},  \\
D_{\text{KL}}(P\|Q) &= \int \mathrm{d}P\log\frac{\mathrm{d}P}{\mathrm{d}Q}, \qquad \qquad 
\chi^2(P\|Q) = \int \frac{(\mathrm{d}P-\mathrm{d}Q)^2}{\mathrm{d}Q}.
\end{align*}
We will frequently use the following inequalities between the above quantities \cite[Chapter 2]{tsybakov2009introduction}: 
\begin{align}
H^2(P,Q) \le \|P-Q\|_{\text{TV}} \le H(P,Q)\sqrt{2-H^2(P,Q)}, \label{eq:TV_Hellinger} \\
\|P-Q\|_{\text{TV}} \le \sqrt{\frac{1}{2}D_{\text{KL}}(P\|Q)} \le \sqrt{\frac{1}{2}\log(1+\chi^2(P\|Q))}. \label{eq:TV_KL_chi} 
\end{align}

Next we define several quantities related to Definition \ref{def:sample_amplification}. For a given distribution class $\mathcal{P}$ and sample sizes $n$ and $m$, the \emph{minimax error of sample amplification} is defined as
\begin{align}\label{eq:SA_distance}
\varepsilon^\star(\calP,n,m) \triangleq \inf_T \sup_{P\in \calP} \| P^{\otimes (n+m)} - P^{\otimes n}\circ T^{-1} \|_{\text{TV}}, 
\end{align}
where the infimum is over all (possibly randomized) measurable mappings $T: \calX^n\to \calX^{n+m}$. For a given error level $\varepsilon$, the \emph{maximum size of sample amplification} is the largest $m$ such that there exists an $(n,n+m,\varepsilon)$ sample amplification, i.e. 
\begin{align}\label{eq:SA_size}
m^\star(\calP, n, \varepsilon) \triangleq \max\{m\in \naturals: \varepsilon^\star(\calP,n,m) \le \varepsilon \}. 
\end{align}
For the ease of presentation, we often choose $\varepsilon$ to be a small constant (say $0.1$) and abbreviate the above quantity as $m^\star(\calP,n)$; we remark that all our results work for a generic $\varepsilon\in (0,1)$. Finally, we also define the \emph{sample amplification complexity} as the smallest $n$ such that an amplification from $n$ to $n+1$ samples is possible: 
\begin{align}\label{eq:SA_complexity}
n^\star(\calP) \triangleq \min\{n\in \naturals: m^\star(\calP,n) \ge 1 \}. 
\end{align}
Note that all the above notions are instance-wise in the distribution class $\calP$.



The minimax error of sample amplification \eqref{eq:SA_distance} is precisely known as the \emph{Le Cam's distance} in the statistics literature. We adopt the standard framework of statistical decision theory \cite{Wald1950statistical}. A statistical model (or experiment) $\calM$ is a tuple $(\calX,(P_{\theta})_{\theta\in\Theta} )$ where an observation $X\sim P_\theta$ is drawn for some $\theta\in\Theta$. A \emph{decision rule} $\delta$ is a regular conditional probability kernel from $\calX$ to the family of probability distributions on a general action space $\calA$, and there is a (measurable) loss function $L: \Theta\times \calA\to \bR_+$. The \emph{risk function} of a given decision rule $\delta$ is defined as
\begin{align}\label{eq:risk}
R_\calM(\theta, \delta, L) \triangleq \bE_\theta[L(\theta,\delta(X))] = \int_\calX \int_\calA L(\theta, a) \delta(da\mid x) P_\theta(dx). 
\end{align}

Based on the definition of risk functions, we are ready to define a metric, known as Le Cam's distance, between statistical models. 
\begin{definition}[Le Cam's distance; see \cite{le1972limits,LeCam1986asymptotic,le1990asymptotics}]\label{def:le_cam_distance}
For two statistical models $\calM = (\calX,(P_\theta)_{\theta\in\Theta})$ and $\calN = (\calY,(Q_\theta)_{\theta\in\Theta})$, \emph{Le Cam's distance} $\Delta(\calM,\calN)$ is defined as
\begin{align*}
\Delta(\calM,\calN) &= \max\left\{\sup_L\sup_{\delta_\calN} \inf_{\delta_\calM} \sup_{\theta\in \Theta} |R_\calM(\theta,\delta_\calM,L) - R_\calN(\theta,\delta_\calN,L)|, \right. \\
&\qquad\qquad\qquad \left.\sup_L\sup_{\delta_\calM} \inf_{\delta_\calN} \sup_{\theta\in \Theta} |R_\calM(\theta,\delta_\calM,L) - R_\calN(\theta,\delta_\calN,L)|  \right\} \\
&= \max \left\{ \inf_{T_1}\sup_\theta \|P_\theta\circ T_1^{-1}  - Q_\theta \|_{\text{\rm TV}}, \  \inf_{T_2}\sup_\theta \|Q_\theta\circ T_2^{-1}  - P_\theta \|_{\text{\rm TV}} \right\},
\end{align*}
where the loss function is taken over all measurable functions $L: \Theta\times \calA\to [0,1]$. 
\end{definition}

In the language of model deficiency introduced in \cite{le1964sufficiency}, Le Cam's distance is the smallest $\varepsilon>0$ such that the model $\calM$ is $\varepsilon$-deficient to the model $\calN$, and $\calN$ is also $\varepsilon$-deficient to $\calM$. In the sample amplification problem, $(P_\theta)_{\theta\in \Theta} = \{P^{\otimes n}: P\in \calP\}, (Q_\theta)_{\theta\in \Theta} = \{P^{\otimes (n+m)}: P\in \calP\}$. Here, choosing $T_2(x^{n+m}) = x^n$ in Definition \ref{def:le_cam_distance} shows that $\calN$ is $0$-deficient to $\calM$, and the remaining quantity involving $T_1$ exactly reduces to the minimax error of sample amplification in \eqref{eq:SA_distance}. Therefore, studying the complexity of sample amplification is equivalent to the characterization of the quantity $\Delta(\calP^{\otimes n}, \calP^{\otimes (n+m)})$. 

\section{Sample amplification via sufficient statistics}\label{sec.sufficiency}
The first idea we present for sample amplification is the classical idea of reduction by sufficiency. Albeit simple, the sufficiency-based idea reduces the problem of generating multiple random vectors to a simpler problem of generating only a few vectors, achieves the optimal complexity of sample amplification in many examples, and is easy to implement. 

\subsection{The general idea}
We first recall the definition of sufficient statistics: in a statistical model $\calM=(\calX,(P_\theta)_{\theta\in\Theta})$ and $X\sim P_\theta$, a statistic $T=T(X)\in \calT$ is \emph{sufficient} if and only if both $\theta - X - T$ and $\theta-T-X$ are Markov chains. A classical result in statistical decision theory is \emph{reduction by sufficiency}, i.e. only the sufficient statistic needs to be maintained to perform statistical tasks as $P_{X|T,\theta}$ does not depend on the unknown parameter $\theta$. In terms of Le Cam's distance, let $\calM \circ T^{-1} = (\calT,(P_\theta \circ T^{-1})_{\theta\in \Theta})$ be the statistical experiment associated with $T$, then sufficiency of $T$ implies that $\Delta(\calM,\calM \circ T^{-1})=0$. Hereafter, we will call $\calM \circ T^{-1}$ the \emph{$T$-reduced model}, or simply \emph{reduced model} in short. 

\begin{algorithm}[t]
	\caption{\textsc{Sample amplification via sufficiency}}\label{alg:sufficiency}
	\begin{algorithmic}[1]
		\State \textbf{Input:} samples $X_1,\cdots,X_n$, a given transformation $f$ between sufficient statistics
		\State Compute the sufficient statistic $T_n = T_n(X_1,\cdots,X_n)$. 
		\State Apply $f$ to the sufficient statistic and compute $\widehat{T}_{n+m} = f(T_n)$. 
		\State Generate $(\widehat{X}_1,\cdots,\widehat{X}_{n+m}) \sim P_{X^{n+m} \mid T_{n+m}} (\cdot \mid  \widehat{T}_{n+m})$. 
		\State \textbf{Output:} amplified samples $(\widehat{X}_1,\cdots,\widehat{X}_{n+m})$. 
	\end{algorithmic}
\end{algorithm}

Reduction by sufficiency could be applied to sample amplification in a simple way, with a general algorithm displayed in Algorithm~\ref{alg:sufficiency}. Suppose that both models $\calP^{\otimes n}$ and $\calP^{\otimes (n+m)}$ admit sufficient statistics $T_n = T_n(X^n)$ and $T_{n+m} = T_{n+m}(X^{n+m})$, respectively. Algorithm~\ref{alg:sufficiency} claims that it suffices to perform sample amplification on the reduced models $\calP^{\otimes n} \circ T_n^{-1}$ and $\calP^{\otimes (n+m)} \circ T_{n+m}^{-1}$, i.e. construct a randomization map $f$ from $T_n$ to $T_{n+m}$. Concretely, the algorithm decomposes into three steps: 
\begin{enumerate}
	\item Step I: map $X^n$ to $T_n$. This step is straightforward: we simply compute $T_n = T_n(X_1,\cdots,X_n)$. 
	\item Step II: apply a randomization map in the reduced model. Upon choosing the map $f$, we simply compute $\widehat{T}_{n+m} = f(T_n)$ with the target that the TV distance $\|\calL(\widehat{T}_{n+m}) - \calL(T_{n+m}) \|_{\text{TV}}$ is uniformly small. The concrete choice of $f$ depends on specific models. 
	\item Step III: map $T_{n+m}$ to $X^{n+m}$. By sufficiency of $T_{n+m}$, the conditional distribution $P_{X^{n+m} \mid T_{n+m}}$ does not depend on the unknown model. Therefore, after replacing the true statistic $T_{n+m}$ by $\widehat{T}_{n+m}$, it is feasible to generate $\widehat{X}^{n+m} \sim P_{X^{n+m} \mid T_{n+m}} (\cdot \mid  \widehat{T}_{n+m})$. To simulate this random vector, it suffices to choose any distribution $P_0\in \calP$ and generate $\widehat{X}^{n+m} \sim  (P_0^{\otimes (n+m)} \mid T_{n+m}(\widehat{X}^{n+m}) = \widehat{T}_{n+m} )$. This step may suffer from computational issues which will be discussed in Section \ref{subsec:computation}.
\end{enumerate}

The validity of this idea simply follows from
\begin{align*}
\Delta(\calP^{\otimes n}, \calP^{\otimes (n+m)}) = \Delta(\calP^{\otimes n}\circ T_n^{-1}, \calP^{\otimes (n+m)}\circ T_{n+m}^{-1}), 
\end{align*}
or equivalently, under each $P\in \calP$, 
\begin{align*}
\| \calL(\widehat{X}^{n+m}) - \calL(X^{n+m} ) \|_{\text{TV}} &= \| \calL(\widehat{T}_{n+m})\times P_{X^{n+m} \mid T_{n+m} }- \calL(T_{n+m} )\times P_{X^{n+m} \mid T_{n+m} } \|_{\text{TV}} \\
&\stepa{=} \| \calL(\widehat{T}_{n+m})- \calL(T_{n+m} ) \|_{\text{TV}} = \| \calL(f(T_n))- \calL(T_{n+m} ) \|_{\text{TV}},
\end{align*}
where (a) is due to the identity $\|P_XP_{Y|X} - Q_XP_{Y|X}\|_{\text{TV}} = \|P_X-Q_X\|_{\text{TV}}$. In other words, it suffices to work on reduced models and find the map $f$ between sufficient statistics. 

This idea of reduction by sufficiency simplifies the design of sample amplification procedures. Unlike in original models where $X^n$ and $X^{n+m}$ typically take values in spaces of different dimensions, in reduced models the sufficient statistics $T_n$ and $T_{n+m}$ are usually drawn from the same space. A simple example is as follows. 
\begin{example}[Gaussian location model with known covariance]\label{example:mean}
Consider the observations $X_1,\cdots,X_n$ from the Gaussian location model $P_\theta = \calN(\theta,\Sigma)$ with an unknown mean $\theta\in \bR^d$ and a known covariance $\Sigma\in \bR^{d\times d}$. To amplify to $n+m$ samples, note that the sample mean vector is a sufficient statistic here, with 
\begin{align*}
T_n(X_1,\cdots,X_n) = \frac{1}{n}\sum_{i=1}^n X_i \sim \calN(\theta, \Sigma/n). 
\end{align*}
Now consider the identity map between sufficient statistics $\widehat{T}_{n+m} = T_n$ used with algorithm \ref{alg:sufficiency}. The amplified samples $(\widehat{X}_1,\cdots,\widehat{X}_{n+m})$ are drawn from  $\calN(0,\Sigma)$ conditioned on the event that $T_{n+m}(\widehat{X}^{n+m}) = \widehat{T}_{n+m} = T_n(X^n)$. For every mean vector $\theta\in \bR^d$ we can upper bound the amplification error of this approach: 
\begin{align*}
\| \calL(\widehat{T}_{n+m}) - \calL(T_{n+m}) \|_{\text{\rm TV}} &= \| \calL(T_n) - \calL(T_{n+m}) \|_{\text{\rm TV}}  \\
&=  \| \calN(\theta, \Sigma/n) - \calN(\theta, \Sigma/(n+m)) \|_{\text{\rm TV}} \\
&\stepa{\le} \sqrt{\frac{1}{2} D_{\text{\rm KL}}(\calN(\theta, \Sigma/n) \| \calN(\theta, \Sigma/(n+m))) } \\
&= \sqrt{\frac{d}{4}\left(\frac{m}{n} - \log\left(1 + \frac{m}{n} \right)\right)} = O\left(\frac{m\sqrt{d}}{n}\right),
\end{align*}
where (a) is due to \eqref{eq:TV_KL_chi}, and the last step holds whenever $m = O(n)$. Therefore, we could amplify $\Omega(n/\sqrt{d})$ additional samples based on $n$ observations, and the complexity of sample amplification in \eqref{eq:SA_complexity} is $n^\star = O(\sqrt{d})$. In contrast, learning this distribution within a small TV distance requires $n=\Omega(d)$ samples, which is strictly harder than sample amplification. This example recovers the upper bound of \cite{axelrod2019sample} with a much simpler analysis, and in later sections we will show that this approach is exactly minimax optimal. 
\end{example}

We make two remarks for the above example. First, the amplified samples $\widehat{X}^{n+m}$ are no longer independent, either marginally or conditioned on $X^n$. Therefore, the above approach is fundamentally different from first estimating the distribution and then generating independent samples from the estimated distribution. Second, the amplified samples do not contain the original samples as a subset. In contrast, a tempting approach for sample amplification is to add $m$ fake samples to the original $n$ observations. However, \cite{axelrod2019sample} showed that any sample amplification approach containing the original samples cannot succeed if $n = o(d/\log d)$ in the above model, and our approach conforms to this result. More examples will be presented in \Cref{subsec:example_sufficiency}. 

\subsection{Computational issues}\label{subsec:computation}
A natural computational question in Algorithm \ref{alg:sufficiency} is how to sample $\widehat{X}^{n+m} \sim P_{X^{n+m} \mid T_{n+m}} (\cdot \mid  \widehat{T}_{n+m})$ in a computationally efficient way. With an additional mild assumption that the sufficient statistic $T$ is also complete (which is easy to find in exponential families), the conditional distribution $P_{X\mid T}$ could be efficiently sampled if we could find a statistic $S = S(X)$ with the following two properties: 
\begin{enumerate}
	\item $S$ is \emph{ancillary}, i.e. $\calL(S)$ is independent of the model parameter $\theta$; 
	\item There is a (measurable) bijection $g$ between $(T,S)$ and $X$, i.e. $X=g(T,S)$ almost surely. 
\end{enumerate}
In fact, if such an $S$ exists, then under any $\theta\in \Theta$, 
\begin{align*}
\calL( X \mid T=t ) \stepa{=} \calL( g(T,S) \mid T=t ) \stepb{=} \calL(g(t,S)),
\end{align*}
where (a) is due to the assumed bijection $g$ between $(T,S)$ and $X$, and (b) is due to a classical result of Basu \cite{basu1955statistics,basu1958statistics} that $S$ and $T$ are independent. Therefore, by the ancillarity of $S$, we could sample $X\sim P_{\theta_0}$ with any $\theta_0\in\Theta$ and compute the statistic $S$ from $X$, then $g(t,S)$ follows the desired conditional distribution $P_{X|T=t}$. An example of this procedure is illustrated below. 

\begin{example}[Computation in Gaussian location model]\label{example:mean_computation}
Consider the setting of Example \ref{example:mean} where $P_\theta = \calN(\theta,\Sigma)^{\otimes (n+m)}, T_{n+m} = (n+m)^{-1}\sum_{i=1}^{n+m} X_i$, and the target is to sample from the distribution $P_{X^{n+m} | T_{n+m}}$. In this model, $T_{n+m}$ is complete and sufficient, and we choose $S=S(X^{n+m})=(S_1,\cdots,S_{n+m-1})$ with $S_i = X_{i+1} - X_1$ for all $i$. Clearly $S$ is ancillary, and $X^{n+m}$ could be recovered from $(T_{n+m}, S)$ via
\begin{align*}
X_1 = T_{n+m} - \frac{\sum_{i=1}^{n+m-1} S_i}{n+m}, \qquad X_{i+1} = X_1 + S_i, \quad i\in [n+m-1].
\end{align*}
Therefore, the choice of $S$ satisfies both conditions. Consequently, we can draw $Z^{n+m} \sim \calN(0,\Sigma)^{\otimes (n+m)}$, compute $S = S(Z^{n+m})$ (where $S_i=Z_{i+1}-Z_1$), and recover $X^{n+m}$ from $(T_{n+m}, S)$. 
\end{example}

The proper choice of $S$ depends on specific models and may require some effort to find; we refer to \Cref{subsec:example_sufficiency} for more examples. We remark that in general there is no golden rule to find $S$. One tempting approach is to find a \emph{maximal} ancillary statistic $S$ such that any other ancillary statistic $S'$ must be a function of $S$. This idea is motivated by the existence of the minimal sufficient statistic under mild conditions and a known computationally efficient approach to compute it. However, for ancillary statistics there is typically no such a maximal one in the above sense, and there may exist uncountably many ``maximal'' ancillary statistics which are not equivalent to each other. From the measure theoretic perspective, this is due to the fact that the family of all ancillary sets is not closed under intersection and thus not a $\sigma$-algebra. In addition, even if a proper notion of ``maximal'' or ``essentially maximal'' could be defined, there is no guarantee that such an ancillary statistic satisfies the bijection condition, and it is hard to determine whether a given ancillary statistic is maximal or not. We refer to \cite{basu1959family,lehmann1992ancillarity} for detailed discussions on ancillarity from mathematical statistics. 

There is also another sampling procedure of $P_{X^{n}|T_n}$ in the conditional inference literature \cite{williams1982some}. Specifically, for each $i\in [n]$, this approach sequentially generates the observation $X_i$ from the one-dimensional distribution $P_{X_i \mid X^{i-1}, T_n}$, which is a simple task as long as we know its CDF. Although this works in simple models such as the Gaussian location model above, in more complicated models exact computation of the CDF is typically not feasible. 

\subsection{General theory for exponential families}\label{subsec:exponential_family}
In this section, we show that a general $(n, n +\Omega(n\varepsilon/\sqrt{d}),\varepsilon)$ sample amplification phenomenon holds for a rich class of exponential families, and is achieved by the sufficiency-based procedure in Algorithm \ref{alg:sufficiency}. Specifically, we consider the following natural exponential family.

\begin{definition}[Exponential family]\label{def:exponential}
 The \emph{exponential family} $(\calX,(P_\theta)_{\theta\in\Theta})$ of probability measures is determined by
\begin{align*}
dP_\theta(x) = \exp(\theta^\top T(x) - A(\theta)) d\mu(x), 
\end{align*}
where $\theta\in \Theta$ is the natural parameter with $\Theta = \{\theta\in \bR^d: A(\theta)<\infty \}$, $T(x)$ is the sufficient statistic, $A(\theta)$ is the log-partition function, and $\mu$ is the base measure. 
\end{definition}

The exponential family includes many widely used probability distributions such as Gaussian, Gamma, Poisson, Exponential, Beta, etc. In the exponential family, the statistic $T(x)$ is sufficient and complete, and several well-known identities include $\bE_{\theta}[T(X)] = \nabla A(\theta)$, and $\mathsf{Cov}_{\theta}[T(X)] = \nabla^2 A(\theta)$. We refer to \cite{diaconis1979conjugate} for a mathematical theory of the exponential family. 

To establish a general theory of sample amplification for exponential families, we shall make the following assumptions on the exponential family. 
\begin{assumption}[Continuity]\label{assump.continuity}
The parameter set $\Theta$ has a non-empty interior. Under each $\theta\in\Theta$, the probability distribution $\calL(T(X))$ is absolutely continuous with respect to the Lebesgue measure. 
\end{assumption}
\begin{assumption}[Moment condition $\mathsf{M}_k$]\label{assump.moment}
For a given integer $k>0$, it holds that
\begin{align*}
\sup_{\theta\in\Theta} \bE_\theta\left[ \left\| (\nabla^2 A(\theta))^{-1/2}(T(X) - \nabla A(\theta) ) \right\|_2^k \right] < \infty. 
\end{align*}
We call it the moment condition $\mathsf{M}_k$. 
\end{assumption}

Assumption \ref{assump.continuity} requires an exponential family of continuous distributions. The motivation is that for continuous exponential family, the sufficient statistics $T_n(X)$ and $T_{n+m}(X)$ for different sample sizes take continuous values in the same space, and it is easier to construct a general transformation. We will propose a different sample amplification approach for discrete statistical models in Section \ref{sec.shuffling}. Assumption \ref{assump.moment} is a moment condition on the normalized statistic $(\nabla^2 A(\theta))^{-1/2}(T(X) - \nabla A(\theta) )$, whose moments always exist as the moment generating function of $T(X)$ exists around the origin. The moment condition $\mathsf{M}_k$ claims that the above normalized statistic has a uniformly bounded $k$-th moment for all $\theta\in \Theta$, which holds in several examples (such as Gaussian, exponential, Pareto) or by considering a slightly smaller $\Theta_0\subseteq \Theta$ bounded away from the boundary. The following lemma presents a sufficient criterion for the moment condition $\mathsf{M}_k$. 
\begin{lemma}\label{lemma:self-concordance}
If the log-partition function $A(\theta)$ satisfies
\begin{align*}
\sup_{\theta\in\Theta} \sup_{u\in \bR^d \backslash \{0\} } \frac{|\nabla^3 A(\theta)[u;u;u]|}{(\nabla^2 A(\theta)[u;u])^{3/2}} \le M < \infty, 
\end{align*}
then the exponential family satisfies the moment condition $\mathsf{M}_k$ for all $k\in \naturals$. Here for a $k$-tensor $T$ and vectors $u_1,\cdots,u_k$, $T[u_1;\cdots;u_k]$ denotes the value of $\langle T, u_1\otimes \cdots \otimes u_k \rangle$. 
\end{lemma}

The condition in Lemma \ref{lemma:self-concordance} is called the self-concordant condition, which is a key concept in the interior point method for convex optimization \cite{nesterov2003introductory}. For example, all quadratic functions and logarithmic functions are self-concordant (which correspond to Gaussian, exponential, and Pareto distributions), and the self-concordance is always fulfilled when $\Theta$ is compact. 

Given any exponential family $\calP$ satisfying Assumptions \ref{assump.continuity} and \ref{assump.moment}, we will show that a simple sample amplification procedure gives a size $\Omega(n/\sqrt{d})$ of sample amplification. Let $X_1,\cdots,X_n$ be i.i.d. samples drawn from $P_\theta$ taking a general form in Definition \ref{def:exponential}, then it is clear that the sample average
\begin{align*}
T_n(X^n) \triangleq \frac{1}{n}\sum_{i=1}^n T(X_i)
\end{align*}
is a sufficient statistic by the factorization theorem. We will apply the general Algorithm \ref{alg:sufficiency} with an identity map between sufficient statistics, i.e. $\widehat{T}_{n+m} = T_n$. The next theorem shows the performance of this approach. 

\begin{theorem}\label{thm:upper_exp_general}
If the exponential family $\calP$ satisfies Assumptions \ref{assump.continuity} and \ref{assump.moment} with $k=3$, then for $\theta\in \Theta$, it holds that
\begin{align*}
\varepsilon^\star(\calP,n,m) \le \|\calL(T_n) - \calL(T_{n+m}) \|_{\text{\rm TV}} \le \frac{C}{\sqrt{n}} + \frac{m\sqrt{d}}{n}, 
\end{align*}
where $C<\infty$ is an absolute constant depending only on $d$ and the moment upper bound in Assumption \ref{assump.moment}. In particular, for sufficiently large $n$, a sample amplification of size $\Omega(n/\sqrt{d})$ is achievable. 
\end{theorem}


Theorem \ref{thm:upper_exp_general} shows that the above simple procedure could achieve a sample amplification from $n$ to $n+\Omega(n/\sqrt{d})$ samples in general continuous exponential families, provided that $n$ is large enough. The main idea behind the proof of Theorem \ref{thm:upper_exp_general} is also simple. We show that the distribution of $T_n$ is approximately $G_n \sim \calN(\nabla A(\theta), \nabla^2 A(\theta)/n)$ by CLT, apply the same CLT for $T_{n+m}$, and then compute the TV distance between two Gaussians as in Example \ref{example:mean}. Theorem \ref{thm:upper_exp_general} is then a direct consequence of the triangle inequality: 
\begin{align*}
&\|\calL(T_n) - \calL(T_{n+m}) \|_{\text{\rm TV}} \\
&\le \|\calL(T_n) - \calL(G_n) \|_{\text{\rm TV}} + \|\calL(G_n) - \calL(G_{n+m}) \|_{\text{\rm TV}} + \|\calL(T_{n+m}) - \calL(G_{n+m}) \|_{\text{\rm TV}}. 
\end{align*}
Note that Assumption \ref{assump.continuity} ensures a vanishing TV distance for the Gaussian approximation, and Assumption \ref{assump.moment} enables us to apply Berry--Esseen type arguments and obtain an $O(1/\sqrt{n})$ convergence rate for the Gaussian approximation. 

The main drawback of Theorem \ref{thm:upper_exp_general} is that there is a hidden constant $C$ depending on the dimension $d$, thus it does not mean that an $(n,n+1,\varepsilon)$ sample amplification is possible as long as $n = \Omega(\sqrt{d}/\varepsilon)$. To tackle this issue, we need to improve the first term in Theorem \ref{thm:upper_exp_general} and find the best possible dependence of the constant $C$ on $d$. We remark that this is a challenging task in probability theory: although the convergence rates of both TV \cite{prohorov1952local,sirazhdinov1962convergence,bally2014distances,bally2016asymptotic} and KL \cite{barron1986entropy,bobkov2014berry} in the CLT result were studied, almost all of them solely focused on the convergence rate on $n$, leaving the tight dependence on $d$ still open. Moreover, direct computation of the quantity $\|\calL(T_n) - \calL(G_n)\|_{\text{\rm TV}}$ shows that even if the random vector $T_n$ has independent components, this quantity is typically at least $\Omega(\sqrt{d/n})$. Therefore, $C=\Omega(\sqrt{d})$ under this proof technique, and a vanishing first term in Theorem \ref{thm:upper_exp_general} requires that $n=\Omega(d)$, which is already larger than the anticipated sample amplification complexity $n=O(\sqrt{d})$. 

To overcome the above difficulties, we make the following changes to both the assumption and analysis. First, to avoid the unknown dependence on $d$, we additionally assume a \emph{product} exponential family, i.e. $P_\theta(dx) = \prod_{i=1}^d p_{\theta_i}(dx_i)$, where each $p_{\theta_i}(x_i)$ is a one-dimensional exponential family. Exploiting the product structure enables to find a constant $C$ depending linearly on $d$. Second, we improve the $O(1/\sqrt{n})$ dependence on $n$ by applying a higher-order CLT result to $T_n$ and $T_{n+m}$, known as the \emph{Edgeworth expansion} \cite{bhattacharya2010normal}. For any $k\ge 2$ and $n\in \naturals$, the signed measure of the Edgeworth expansion on $\bR^d$ is
\begin{align}\label{eq:edgeworth}
\Gamma_{n,k}(dx) = \gamma(x)\left( 1 + \sum_{\ell=1}^{\lfloor k/3 \rfloor} \frac{\calK_\ell(x)}{n^{\ell/2}} \right)dx, 
\end{align}
where $\gamma(x)$ is the density of a standard normal random variable on $\bR^d$, and $\calK_m(x)$ is a polynomial of degree $3m$ which depends only on the first $3m$ moments of the distribution. We note that unlike CLT, the general Edgeworth expansion is a signed measure with possibly negative densities; however, it is close to Gaussian with an $O(n^{-1/2})$ approximation error. Such a higher-order expansion enables us to have better Berry-Esseen type bounds, but upper bounding $\|\Gamma_{n,k} - \Gamma_{n+m,k}\|_{\text{\rm TV}}$ becomes more complicated and requires to handle the Gaussian part and the correction part separately; see Appendix \ref{append:edgeworth} for details. In particular, we could improve the error dependence on $n$ from $O(1/\sqrt{n})$ to $O(1/n^2)$. 

Formally, the next theorem shows a better sample amplification performance for product exponential families. 
\begin{theorem}\label{thm:upper_exp_product}
  Let $(\calX,\calP = (P_\theta)_{\theta\in \Theta} )$ be a product exponential family, where each one-dimensional component satisfies Assumptions \ref{assump.continuity} and \ref{assump.moment} with $k=10$. Then for $\theta\in \Theta$, it holds that
	\begin{align*}
	\varepsilon^\star(\calP,n,m) \le \|\calL(T_n) - \calL(T_{n+m}) \|_{\text{\rm TV}} \le C\left(\frac{d}{n^2} + \frac{m\sqrt{d}}{n} \right), 
	\end{align*}
	where $C<\infty$ is an absolute constant independent of $(n,d)$. In particular, as long as $n=\Omega(\sqrt{d}/\varepsilon)$, an $(n,n+m,\varepsilon)$ sample amplification of size $m=\Omega(n\varepsilon/\sqrt{d})$ is achievable. 
\end{theorem}
Theorem \ref{thm:upper_exp_product} shows that for product exponential family, we not only achieve the amplification size $m=\Omega(n\varepsilon/\sqrt{d})$, but also attain a sample complexity $n=O(\sqrt{d}/\varepsilon)$ for sample amplification. This additional result on sample complexity is important in the sense that, even if distribution learning is impossible, it is possible to perform sample amplification. Although the independence or even the exponential family assumption could be strong practically, in \Cref{subsec:example_sufficiency} we show that both phenomena $m=\Omega(n\varepsilon/\sqrt{d})$ and $n=O(\sqrt{d}/\varepsilon)$ hold in many natural models. 
\section{Sample amplification via learning}\label{sec.shuffling}
The sufficiency-based approach of sample amplification is not always desirable. First, models outside the exponential family typically do not admit non-trivial sufficient statistics, and therefore the reduction by sufficiency may not be very helpful. Second, the identity map applied to the sufficient statistics only works for continuous models, and incurs a too large TV distance when the underlying model is discrete. Third, previous approaches are not directly related to learning the model, so a general relationship between learning and sample amplification is largely missing. In this section, we propose another sample amplification approach, and show that how a good learner helps to obtain a good sample amplifier. 

\subsection{The general result}
For a class of distribution $\calP$ and $n$ i.i.d. observations drawn from an unknown $P\in \calP$, we define the following notion of the $\chi^2$-estimation error. 

\begin{definition}[$\chi^2$-estimation error]\label{def:learning_error}
For a class of distributions $\calP$ and sample size $n$, the \emph{$\chi^2$-estimation error} $r_{\chi^2}(\calP,n)$ is defined to be the minimax estimation error under the expected $\chi^2$-divergence: 
\begin{align*}
r_{\chi^2}(\calP,n) \triangleq \inf_{\widehat{P}_n}\sup_{P\in \calP} \bE_P[\chi^2(\widehat{P}_n, P)], 
\end{align*}
where the infimum is taken over all possible distribution estimators $\widehat{P}_n$ based on $n$ samples. 
\end{definition}

Roughly speaking, the $\chi^2$-estimation error in the above definition characterizes the complexity of the distribution class $\calP$ in terms of distribution learning under the $\chi^2$-divergence. The main result of this section is to show that, the error of sample amplification in \eqref{eq:SA_distance} could be upper bounded by using the $\chi^2$-estimation error.

\begin{theorem}\label{thm:shuffling_general}
For general $\calP$ and $n,m\ge 0$, it holds that
\begin{align*}
\varepsilon^\star(\calP,n,m) \le \sqrt{\frac{m^2}{n}\cdot r_{\chi^2}(\calP,n/2) }. 
\end{align*}	
\end{theorem}

The following corollary is immediate from Theorem \ref{thm:shuffling_general}. 

\begin{corollary}\label{cor:shuffling_general}
	An $(n,n+m,\varepsilon)$ sample amplification is possible if $m = O(\varepsilon\sqrt{n/r_{\chi^2}(\calP,n/2)})$. Moreover, the sample complexity of amplification in \eqref{eq:SA_complexity} satisfies
	\begin{align*}
	n^\star(\calP) = O\left( \min\left\{n\in \naturals: r_{\chi^2}(\calP,n/2) \le n \right\} \right). 
	\end{align*}
\end{corollary}

\begin{remark}
Although the error of sample amplification in \eqref{eq:SA_distance} is measured under the TV distance, the same result holds for the squared root of the KL divergence (which by \eqref{eq:TV_KL_chi} is no smaller than the TV distance). 
\end{remark}

The above result provides a quantitative guarantee that the sample amplification is easier than learning (under the $\chi^2$-divergence). Specifically, the sample complexity of learning is the smallest $n\in \naturals$ such that $r_{\chi^2}(\calP,n) = O(1)$, while Corollary \ref{cor:shuffling_general} shows that the complexity for amplification is at most the smallest $n\in\naturals$ such that $r_{\chi^2}(\calP,n/2)=O(n)$. As $r_{\chi^2}(\calP,n)$ is non-increasing in $n$, this means that the learning complexity is in general larger. 

When the distribution class $\calP$ has a product structure $\calP = \prod_{j=1}^d \calP_j$, the next theorem shows a better relationship between the amplification error and the learning error. 

\begin{theorem}\label{thm:shuffling_product}
For $\calP = \prod_{j=1}^d \calP_j$ and $n,m\ge 0$, it holds that
\begin{align*}
\varepsilon^\star(\calP,n,m) \le \sqrt{\frac{m^2}{n} \sum_{j=1}^d r_{\chi^2}(\calP_j,n/2) }. 
\end{align*}	
\end{theorem}
	
\begin{corollary}\label{cor:shuffling_product}
For product models, an $(n,n+m,\varepsilon)$ sample amplification is achievable if 
\begin{align*}
m = O\left(\varepsilon\sqrt{\frac{n}{\sum_{j=1}^d r_{\chi^2}(\calP_j,n/2)}}\right). 
\end{align*}
Moreover, the sample complexity of amplification in \eqref{eq:SA_complexity} satisfies
\begin{align*}
n^\star(\calP) = O\left( \min\left\{n\in \naturals: \sum_{j=1}^d r_{\chi^2}(\calP_j,n/2) \le n \right\} \right). 
\end{align*}
\end{corollary}	

We observe that the result of Theorem \ref{thm:shuffling_product} typically improves over Theorem \ref{thm:shuffling_general} for product models. In fact, since
\begin{align*}
\chi^2\left(\prod_{j=1}^d P_j , \prod_{j=1}^d Q_j \right) = \prod_{j=1}^d (1+\chi^2(P_j,Q_j)) - 1\ge \sum_{j=1}^d \chi^2(P_j,Q_j), 
\end{align*}
the inequality $\sum_{j=1}^d r_{\chi^2}(\calP_j, n/2) \le r_{\chi^2}(\calP,n/2)$ typically holds. Moreover, there are scenarios where we have $\sum_{j=1}^d r_{\chi^2}(\calP_j, n/2) \ll r_{\chi^2}(\calP,n/2)$, thus Theorem \ref{thm:shuffling_product} provides a substantial improvement over Theorem \ref{thm:shuffling_general}. For example, when $\calP = (\calN(\theta,I_d))_{\theta\in \bR^d}$, it could be verified that $r_{\chi^2}(\calP_j, n/2) = O(1/n)$ for each $j\in [d]$, while $r_{\chi^2}(\calP,n/2) = \exp(\Omega(d/n))-1$. Hence, in the important regime $\sqrt{d}\ll n\ll d$ where learning is impossible but the sample amplification is possible, Theorem \ref{thm:shuffling_product} is strictly stronger than Theorem \ref{thm:shuffling_general}. 

\begin{remark}
In the above Gaussian location model, there is an alternative way to conclude that Theorem \ref{thm:shuffling_product} is strictly stronger than Theorem \ref{thm:shuffling_general}. We will see that the shuffling approach achieving the bound in Theorem \ref{thm:shuffling_general} keeps all the observed samples, whereas \cite{axelrod2019sample} shows that all such approaches must incur a sample complexity $n=\Omega(d/\log d)$ for the Gaussian model. In contrast, Theorem \ref{thm:shuffling_product} and Corollary \ref{cor:shuffling_product} give a sample complexity $n=O(\sqrt{d})$ of amplification in the Gaussian location model. 
\end{remark}

\subsection{The shuffling approach}
This section presents the sample amplification approaches to achieve Theorems \ref{thm:shuffling_general} and \ref{thm:shuffling_product}. The idea is simple: we find a good distribution learner $\widehat{P}_n$ which attains the rate-optimal $\chi^2$-estimation error, draw additional $m$ samples $Y_1,\cdots,Y_m$ from $\widehat{P}_n$, and shuffle them with the original samples $X_1,\cdots,X_n$ uniformly at random. This approach suffices to achieve the sample amplification error in Theorem \ref{thm:shuffling_general}, while for Theorem \ref{thm:shuffling_product} an additional trick is applied: instead of shuffling the whole vectors, we independently shuffle each coordinate instead. For technical reasons, in both approaches we apply the sample splitting: the first $n/2$ samples are used for the estimation of $\widehat{P}_n$, while the second $n/2$ samples are used for shuffling. The algorithms are summarized in Algorithms \ref{alg:shuffling_general} and \ref{alg:shuffling_product}. 

\begin{algorithm}[t]
	\caption{\textsc{Sample amplification via shuffling: general model}}\label{alg:shuffling_general}
	\begin{algorithmic}[1]
		\State \textbf{Input:} samples $X_1,\cdots,X_n$, a given class of distributions $\calP$.
		\State Based on samples $X_1,\cdots,X_{n/2}$, find an estimator $\widehat{P}_n$ such that
		\begin{align*}
		 \sup_{P\in \calP} \bE_P[\chi^2(\widehat{P}_n, P)] \le C\cdot r_{\chi^2}(\calP,n/2).
		\end{align*}
		\State Draw $m$ additional samples $Y_1,\cdots,Y_m$ from $\widehat{P}_n$. 
		\State Uniformly at random, shuffle the pool of $X_{n/2+1},\cdots,X_{n},Y_1,\cdots,Y_m$ to obtain $(Z_1,\cdots,Z_{n/2+m})$.
		\State \textbf{Output:} amplified samples $(X_1,\cdots,X_{n/2},Z_1,\cdots,Z_{n/2+m})$. 
	\end{algorithmic}
\end{algorithm}

\begin{algorithm}[t]
	\caption{\textsc{Sample amplification via shuffling: product model}}\label{alg:shuffling_product}
	\begin{algorithmic}[1]
		\State \textbf{Input:} samples $X_1,\cdots,X_n$, a given class of product distributions $\calP=\prod_{j=1}^d \calP_j$
		\For{$j=1,2,\cdots,d$}
		\State Based on samples $X_{1,j},\cdots,X_{n/2,j}$, find an estimator $\widehat{P}_{n,j}$ such that
		\begin{align*}
		\sup_{P_j\in \calP_j} \bE_{P_j}[\chi^2(\widehat{P}_{n,j}, P_j)] \le C\cdot r_{\chi^2}(\calP_j,n/2).
		\end{align*}
		\State Draw $m$ additional samples $Y_{1,j},\cdots,Y_{m,j}$ from $\widehat{P}_{n,j}$.
		\State Uniformly at random, shuffle $X_{n/2+1,j},\cdots,X_{n,j},Y_{1,j},\cdots,Y_{m,j}$ to obtain $(Z_{1,j},\cdots,Z_{n/2+m,j})$. 
		\EndFor
		\State For each $i\in [n/2+m]$, form the vector $Z_i = (Z_{i,1},\cdots,Z_{i,d})$.
		\State \textbf{Output:} amplified samples $(X_1,\cdots,X_{n/2},Z_1,\cdots,Z_{n/2+m})$. 
	\end{algorithmic}
\end{algorithm}

The following lemma is the key to analyze the performance of the shuffling approach. 
\begin{lemma}\label{lemma:shuffle_chi}
	Let $X_1,\cdots,X_n$ be i.i.d. drawn from $P$, and $Y_1,\cdots,Y_m$ be i.i.d. drawn from $Q$ independent of $(X_1,\cdots,X_n)$. Let $(Z_1,\cdots,Z_{n+m})$ be a uniformly random permutation of  $(X_1,\cdots,X_n,Y_1,\cdots,Y_m)$ , and $P_{\text{\rm mix}}$ be the distribution of the random mixture $(Z_1,\cdots,Z_{n+m})$. Then
	\begin{align*}
	\chi^2\left(P_{\text{\rm mix}}, P^{\otimes (n+m)} \right) \le \left(1 + \frac{m}{n+m}\chi^2(Q,P) \right)^m - 1. 
	\end{align*}
\end{lemma}

Based on Lemma \ref{lemma:shuffle_chi}, the advantage of random shuffling is clear: if we simply append $Y_1,\cdots,Y_m$ to the end of the original sequence $X_1,\cdots,X_n$, then the $\chi^2$-divergence is exactly $(1+\chi^2(Q,P))^m-1$. By comparing with the upper bound in Lemma \ref{lemma:shuffle_chi}, we observe that a smaller coefficient $m/(n+m)$ is applied to the individual $\chi^2$-divergence after a random shuffle. The proofs of Theorems \ref{thm:shuffling_general} and \ref{thm:shuffling_product} are then clear, where we simply take $Q=\widehat{P}_n$ and apply the above lemma. Note that the statement of Lemma \ref{lemma:shuffle_chi} requires that $Y_1,\cdots,Y_m$ be independent of $X_1,\cdots,X_n$, which is exactly the reason why we apply sample splitting in Algorithms \ref{alg:shuffling_general} and \ref{alg:shuffling_product}. The proof of Lemma \ref{lemma:shuffle_chi} is presented in Appendix \ref{appendix:lemma}, and the complete proofs of Theorems \ref{thm:shuffling_general} and \ref{thm:shuffling_product} are relegated to Appendix \ref{appendix:thm}. We also include concrete examples of the shuffling approach in \Cref{subsec:example_shuffling}. 
\section{Minimax lower bounds}\label{sec:lower_bound}
In this section we establish minimax lower bounds for sample amplification in different statistical models. Section \ref{subsec:lower_general} presents a general and tight approach for establishing the lower bound, which leads to an exact sample amplification result for the Gaussian location model. Based on this result, we show that for $d$-dimensional continuous exponential families, the sample amplification size cannot exceed $\omega(n\varepsilon/\sqrt{d})$ for sufficiently large sample size $n$. Section \ref{subsec:lower_product} provides a specialized criterion for product models, where we show that $n=\Omega(\sqrt{d}/\varepsilon)$ and $m=O(n\varepsilon/\sqrt{d})$ are always valid lower bounds, with hidden constants independent of all parameters. \Cref{subsec:lower_example} lists several concrete examples where our general idea could be properly applied to provide tight and non-asymptotic results. 

\subsection{General idea}\label{subsec:lower_general}
The main tool to establish the lower bound is the first equality in the Definition \ref{def:le_cam_distance} of Le Cam's distance. Specifically, for a class of distributions $\calP=(P_\theta)_{\theta\in\Theta}$, let $\mu$ be a given prior distribution on $\Theta$, and $L: \Theta\times \calA\to [0,1]$ be a given non-negative loss function upper bounded by one. Given $n$ i.i.d. samples from an unknown distribution in $\calP$, define the following \emph{Bayes risk} and \emph{minimax risk}: 
\begin{align*}
r_{\text{\rm B}}(\calP,n,L,\mu) &= \inf_{\widehat{\theta}} \int_{\Theta} \bE_{\theta}[L(\theta, \widehat{\theta}(X^n))] \mu(d\theta), \\
r(\calP,n,L) &= \inf_{\widehat{\theta}} \sup_{\theta\in\Theta} \bE_{\theta}[L(\theta, \widehat{\theta}(X^n))], 
\end{align*}
where the infimum is over all possible estimators $\widehat{\theta}(\cdot)$ taking value in $\calA$. The following result is a direct consequence of Definition \ref{def:le_cam_distance}. 

\begin{lemma}\label{lemma:general_lower_bound}
For any integer $n,m>0$, any class of distributions $\calP = (P_\theta)_{\theta\in \Theta}$, any prior $\mu$ on $\Theta$, and any loss function $L: \Theta\times \calA\to [0,1]$, the minimax error of sample amplification $\varepsilon^\star(\calP,n,m)$ in \eqref{eq:SA_distance} satisfies that
\begin{align*}
\varepsilon^\star(\calP,n,m) &\ge r_{\text{\rm B}}(\calP,n,L,\mu) - r_{\text{\rm B}}(\calP,n+m,L,\mu), \\
\varepsilon^\star(\calP,n,m) &\ge r(\calP,n,L) - r(\calP,n+m,L). 
\end{align*}
\end{lemma}

Based on Lemma \ref{lemma:general_lower_bound}, it suffices to find an appropriate prior distribution $\mu$ and a loss function $L$, and then compute (or lower bound) the difference between the Bayes or minimax risks with different sample sizes. We note that the lower bound technique in \cite{axelrod2019sample}, albeit seemingly different, is a special case of Lemma \ref{lemma:general_lower_bound}. Specifically, the authors of \cite{axelrod2019sample} designed a set-valued mapping $A_n: \theta\to \calP(\calX^n)$ for each $n\in \mathbb{N}$ such that $\bP_\theta(X^{n+m} \in A_{n+m}(\theta))\ge 0.99$ for all $\theta\in\Theta$, while there is a prior distribution $\mu$ on $\Theta$ such that
\begin{align}\label{eq:coverage}
\bE_{X^n}\left[ \sup_{x^n\in \calX^n} \bP_{\theta \mid X^n}(x^n \in A_n(\theta)) \right] \le 0.5. 
\end{align}
If the above conditions hold, then an $(n,n+m)$ sample amplification is impossible. Note that the probability term in \eqref{eq:coverage} is the maximum coverage probability of the sets $A_n(\theta)$ where $\theta$ follows the posterior distribution $\bP_{\theta \mid X^n}$, which is a well-defined geometric object when both $A_n(\theta)$ and the posterior are known. To see that the above approach falls into our framework, consider the loss function $L: \Theta \times \cup_{n\ge 1}\calX^n\to [0,1]$ with $L(\theta,X^n) = \1(X^n\notin A_n(\theta))$. Then the first condition ensures that $r_{\text{B}}(\calP,n+m,L,\nu)\le 0.01$ for each prior $\nu$, and the second condition \eqref{eq:coverage} exactly states that $r_{\text{B}}(\calP,n,L,\mu)\ge 0.5$ for the chosen prior $\mu$. 

A first application of Lemma \ref{lemma:general_lower_bound} is an \emph{exact} lower bound in Gaussian location models. 
\begin{theorem}\label{thm:Gaussian}
	For the Gaussian location model $\calP = \{\calN(\theta,\Sigma)\}_{\theta\in\bR^d}$ with a fixed covariance $\Sigma\in \bR^{d\times d}$, the minimax error of sample amplification in \eqref{eq:SA_distance} is exactly
	\begin{align*}
	\varepsilon^\star(\calP,n,m) = \left\|\calN\left(0,\frac{I_d}{n}\right) - \calN\left(0,\frac{I_d}{n+m}\right) \right\|_{\text{\rm TV}}. 
	\end{align*}
	In particular, the sufficiency-based approach in Example \ref{example:mean} is exactly minimax optimal. 
\end{theorem}

Theorem \ref{thm:Gaussian} shows that an exact error characterization for the Gaussian location model is possible through the general lower bound approach in Lemma \ref{lemma:general_lower_bound}. This result is also asymptotically useful to a rich family of models: note that by CLT, the sufficient statistic in a continuous exponential family follows a Gaussian distribution asymptotically, with a vanishing TV distance. This idea was used in Section \ref{subsec:exponential_family} to establish the $O(n\varepsilon/\sqrt{d})$ upper bound, and the same observation could lead to an $\Omega(n\varepsilon/\sqrt{d})$ lower bound as well, under slightly different assumptions. Specifically, we drop Assumption \ref{assump.moment} while introducing an additional assumption.

\begin{assumption}[Linear independence]\label{assump.independence}
	The components of sufficient statistic $T(x)$ are linearly independent, i.e. $a^\top T(x)=0$ for $\mu$-almost all $x\in \calX$ implies $a=0$. 
\end{assumption}

Assumption \ref{assump.independence} ensures that the true dimension of the exponential family is indeed $d$. Whenever Assumption \ref{assump.independence} does not hold, we could transform it into a minimal exponential family with a lower dimension fulfilling this assumption. Note that when Assumptions \ref{assump.continuity} and \ref{assump.independence} hold, the mean mapping $\theta\mapsto \nabla A(\theta)$ is a diffeomorphism between $\Theta$ and $\nabla A(\theta)$; see, e.g. \cite[Theorem 1.22]{miescke2008statistical}. Therefore, $\nabla A(\cdot)$ is an open map, and the set $\{\nabla A(\theta): \theta\in \Theta \}$ contains a $d$-dimensional ball. This fact enables us to obtain a $d$-dimensional Gaussian location model after we apply the CLT.

The following theorem characterizes an asymptotic lower bound for every exponential family satisfying Assumptions \ref{assump.continuity} and \ref{assump.independence}. 

\begin{theorem}\label{thm:lower_exponential}
	Given a $d$-dimensional exponential family $\calP$ satisfying Assumptions \ref{assump.continuity} and \ref{assump.independence}, for every $n,m\in \naturals$, the minimax error of sample amplification satisfies
	\begin{align*}
	\varepsilon^\star(\calP,n,m)\ge c\cdot \left(\frac{m\sqrt{d}}{n}\wedge 1\right) - C\cdot \left(\frac{\log n}{n}\right)^{\frac{1}{3}},
	\end{align*}
	where $c>0$ is an absolute constant independent of $(n,m,d,\calP)$, and constant $C>0$ depends only on the exponential family (and thus on $d$). 
\end{theorem}

Theorem \ref{thm:lower_exponential} shows that there exists some $n_0>0$ depending only on the exponential family, such that sample amplification from $n$ to $n+\omega(n\varepsilon/\sqrt{d})$ samples is impossible for all $n>n_0$. However, similar to the nature of the upper bound in Theorem \ref{thm:upper_exp_general}, this asymptotic result does not imply that the sample amplification is impossible if $n=o(\sqrt{d}/\varepsilon)$. Nevertheless, in the following sections we show that the sample complexity lower bound $n=\Omega(\sqrt{d}/\varepsilon)$ indeed holds in product families, as well as several other concrete examples.

\subsection{Product models}\label{subsec:lower_product}
Although Lemma \ref{lemma:general_lower_bound} presents a lower bound argument in general, the computation of exact Bayes or minimax risks could be very challenging, and the usual rate-optimal analysis (i.e. bounding the risks within a multiplicative constant) will not lead to meaningful results. In addition, choosing the right prior and loss is a difficult task which may change from instance to instance. Therefore, it is helpful to propose specialized versions of Lemma \ref{lemma:general_lower_bound} which are easier to work with. Surprisingly, such a simple version exists for product models, which is summarized in the following theorem. 

\begin{theorem}\label{thm:lower_bound}
	Let $\varepsilon \in (0,1)$ and $P_{\theta} = \prod_{j=1}^d p_{\theta_j}$ be a  product model with $(\theta_1,\cdots,\theta_d)\in \prod_{j=1}^d \Theta_j$. Suppose for each $j\in [d]$, there exist two points $\theta_{j,+}, \theta_{j,-}\in \Theta_j$ such that
	\begin{align}
	\| p_{\theta_{j,+}}^{\otimes n} - p_{\theta_{j,-}}^{\otimes n} \|_{\text{\rm TV}} &\le \alpha_j - \frac{\epsilon}{\sqrt{d}}, \label{eq:small_TV} \\
	\| p_{\theta_{j,+}}^{\otimes (n+m)} - p_{\theta_{j,-}}^{\otimes (n+m)} \|_{\text{\rm TV}} &\ge \alpha_j + \frac{\epsilon}{\sqrt{d}}, \label{eq:large_TV} 
	\end{align}
	with $\alpha_j\in (\underline{\alpha}, \overline{\alpha})$, where $\underline{\alpha}, \overline{\alpha}\in (0,1)$ are absolute constants. Then there exists an absolute constant $c=c(\underline{\alpha},\overline{\alpha})>0$ such that
	\begin{align*}
	\varepsilon^\star(\calP,n,m) \ge c\varepsilon. 
	\end{align*}
\end{theorem}

Theorem \ref{thm:lower_bound} leaves the choices of the prior and loss function in Lemma \ref{lemma:general_lower_bound} implicit, and provides a simple criterion for product models. The usual routine of applying Theorem \ref{thm:lower_bound} is as follows: fix any constant $\alpha$ and a target error $\varepsilon$, find for each $j\in [d]$ two points $\theta_{j,+}, \theta_{j,-}\in \Theta_j$ such that the condition \eqref{eq:small_TV} holds for a given sample size $n$. Then the condition \eqref{eq:large_TV} becomes an inequality solely for $m$, from which we could solve the smallest $m_j\in \naturals$ such that \eqref{eq:large_TV} holds along the $j$-th coordinate. Finally, the sample amplification from $n$ to $n+m$ samples is impossible by the above theorem, where $m = \max_{j\in [d]} m_j$. Although Theorem \ref{thm:lower_product} is only for product models, similar ideas could also be applied to non-product models; we refer to \Cref{subsec:lower_example} for concrete examples. 

Theorem \ref{thm:lower_bound} also provides some intuition on why the sample complexity lower bound for amplification is typically smaller than that of learning. Specifically, for learning under the TV distance, a small TV distance $\|\prod_{j=1}^d p_{\theta_{j,+}}^{\otimes n} - \prod_{j=1}^d p_{\theta_{j,-}}^{\otimes n}  \|_{\text{TV}}$ between product distributions is required. This requirement typically leads to a much smaller individual TV distance $\| p_{\theta_{j,+}}^{\otimes n} - p_{\theta_{j,-}}^{\otimes n} \|_{\text{\rm TV}}$, e.g. $O(1/\sqrt{d})$ for many regular models. In contrast, the conditions \eqref{eq:small_TV} and \eqref{eq:large_TV} only require a constant individual TV distance, which leads to a smaller sample complexity $n$ in the sample amplification lower bound. To understand why a larger individual TV distance works for sample amplification, in the proof of Theorem \ref{thm:lower_bound} we consider the uniform prior on $2^d$ points $\prod_{j=1}^d \{\theta_{j,+}, \theta_{j,-} \}$. Under this prior, the test accuracy for each dimension is precisely $(1+\text{TV}_j)/2$, which is slightly smaller than $(1+\alpha)/2$ with $n$ samples, and slightly larger than $(1+\alpha)/2$ with $n+m$ samples (assuming $\alpha_j\equiv \alpha$). Therefore, if a unit loss is incurred when the fraction of correct tests does not exceed $(1+\alpha)/2$, the current scaling in \eqref{eq:small_TV}, \eqref{eq:large_TV} shows that there is an $\Omega(\varepsilon)$ difference in the expected loss under different sample sizes. In other words, such an aggregate voting test helps to have a larger individual TV distance. The details of the proof are deferred to Appendix \ref{appendix:thm}. 

Theorem \ref{thm:lower_bound} has a far-reaching consequence: with almost no assumption on the product model $\calP$, for any $c>0$ it always holds that $\varepsilon^\star(\calP,n,\lceil c\varepsilon n/\sqrt{d}\rceil) \ge c'\varepsilon$ for some absolute constant $c'>0$ \emph{independent} of the product model $\calP$. The only assumption (besides the product structure) we make on $\calP$ is as follows (here $n\in \naturals$ is a given sample size): 

\begin{assumption}\label{assump.hellinger}
Let $\calP$ possess the product structure as in Theorem \ref{thm:lower_bound}. For each $j\in [d]$, there exists two points $\theta_{j,+}, \theta_{j,-}\in \Theta_j$ such that $1/(10n) \le H^2(p_{\theta_{j,+}}, p_{\theta_{j,-}} ) \le 1/(5n)$. 
\end{assumption}

Assumption \ref{assump.hellinger} is a mild assumption that  requires that for each coordinate, the map $\theta_j \mapsto p_{\theta_j}$ is continuous under the Hellinger distance. This assumption is satisfied for almost all practical models, either discrete or continuous, and is invariant with model reparametrizations or bijective transformation of observations. We note that the coefficients $1/10$ and $1/5$ are not essential, and could be replaced by any smaller constants. The next theorem states that if Assumption \ref{assump.hellinger} holds, we always have a lower bound $n=\Omega(\sqrt{d})$ for the sample complexity and an upper bound $m=O(n/\sqrt{d})$ for the size of sample amplification.

\begin{theorem}\label{thm:lower_product}
Let $\calP$ be a product model satisfying Assumption \ref{assump.hellinger}. Then for any $c>0$, there is some $c'>0$ depending only on $c$ (thus independent of $n,d,\varepsilon,\calP$) such that
\begin{align*}
\varepsilon^\star\left(\calP,n, \left\lceil \frac{c\varepsilon n}{\sqrt{d}} \right\rceil \right) \ge c'\varepsilon. 
\end{align*}
\end{theorem}

Theorem \ref{thm:lower_product} is a general lower bound for sample amplification in product models, with intriguing properties that it is instance-wise in the model $\calP$, while the constants $c$ and $c'$ are \emph{independent} of $\calP$. As a result, the sample complexity is uniformly $\Omega(\sqrt{d}/\varepsilon)$, and the maximum size of sample amplification is uniformly $O(n\varepsilon/\sqrt{d})$ for all product models. In comparison, the matching upper bound in Theorem \ref{thm:upper_exp_product} for product models has a hidden constant depending on the statistical model. We note that it is indeed natural to have sample amplification results independent of the underlying statistical model. For example, it is clear by definition that sample amplifications are invariant with bijective transformation of observations. However, Assumption \ref{assump.moment} depends on such transformations, so it possibly contains some redundancy. In contrast, Assumption \ref{assump.hellinger} remains invariant, which is therefore more natural. 

The proof idea of Theorem \ref{thm:lower_product} is best illustrated for the case $d=1$. Using the two points $\theta_+, \theta_-$ in Assumption \ref{assump.hellinger}, one could show that the TV distance between $n$ copies of $p_{\theta_+}$ and $p_{\theta_-}$ is bounded from above by a small constant. Similarly, for a large $C>0$, the TV distance between $Cn$ copies of them is lower bounded by a large constant. Consequently, if $m=(C-1)n$, Theorem \ref{thm:lower_bound} applied with $d=1$ gives an $\Omega(1)$ lower bound on $\varepsilon^\star(\calP,n,m)$. What happens if $m=c\varepsilon n$ with a small $c$? The idea is to consider the TV distances between $n, n+c\varepsilon n, n+2c\varepsilon n, \cdots, Cn$ copies of $p_{\theta_+}$ and $p_{\theta_-}$, which is an increasing sequence. Now by the pigeonhole principle, there must be two adjacent TV distances differing by at least $\Omega(c\varepsilon/C)=\Omega(\varepsilon)$, and Theorem \ref{thm:lower_bound} could be applied to this pair of sample sizes. This idea is easily generalized to any dimensions, with the full proof in Appendix \ref{appendix:thm}. 

We note that the lower bounds in Theorem \ref{thm:lower_exponential} (as well as \ref{thm:Gaussian}) and Theorem \ref{thm:lower_product} are on two different ends of the spectrum. In Theorems \ref{thm:Gaussian} and \ref{thm:lower_exponential}, an asymptotic setting (i.e. $d$ fixed and $n\to\infty$) is essentially considered, and a Gaussian limit is crucially used as long as there is local asymptotic normality. In comparison, Theorem \ref{thm:lower_product} deals with a high-dimensional scenario ($n,d$ can grow together) but restricts to a product model. However, looking at product submodels and/or exploiting its proof techniques could still lead to tight lower bounds for several non-product models, as shown in the examples in \Cref{subsec:lower_example}.

\section{Discussions on sample amplification versus learning}\label{sec:discussion}

In all the examples we have seen in the previous sections, there is always a squared root relationship between the statistical complexities of sample amplification and learning. Specifically, when the dimensionality of the problem is $d$, the complexity of learning the distribution (under a small TV distance) is typically $n=\Theta(d)$, whereas that of the sample amplification is typically $n=\Theta(\sqrt{d})$. In this section, we give examples where this relationship could break down in either direction, thus show that there is no universal scaling between the sample complexities of amplification and learning.  

\subsection{An example where the complexity of sample amplification is $o(\sqrt{d})$}
We first provide an example where the distribution learning is hard, but an $(n,n+1,0.1)$ sample amplification is easy. Consider the following class $\calP_{d,t}$ of discrete distributions: 
\begin{align*}
    \calP_{d,t} = \left\{(p_0,\cdots,p_d): p_i\ge 0, \sum_{i=0}^d p_i = 1, p_0 = t \right\}, 
\end{align*}
where it is the same as the class of all discrete distributions over $d+1$ points, except that the learner has the perfect knowledge of $p_0 = t$ for some known $t\in [1/(2\sqrt{d}), 1/2]$. It is a classical result (see, e.g. \cite{han2015minimax}) that the sample complexity of learning the distribution over $\calP_{d,t}$ with a small TV distance is still $n=\Theta(d)$, regardless of $t$. However, the next theorem shows that the complexity of sample amplification is much smaller. 

\begin{theorem}
\label{thm:distribution_top_element}
For the class $\calP_{d,t}$ with $t\in [1/(2\sqrt{d}), 1/2]$, an $(n,n+1,0.1)$ sample amplification is possible if and only if
\begin{align*}
    n = \Omega\left(\frac{1}{t}\right). 
\end{align*}
\end{theorem}

Note that for the choice of $t = \Theta(d^{-\alpha})$ with $\alpha\in [0,1/2]$, the complexity of sample amplification could possibly be $n=\Theta(d^\alpha)$ for every $\alpha\in [0,1/2]$, showing that it could be $o(\sqrt{d})$ with an arbitrary polynomial scale in $d$. Moreover, if $t = o(1/\sqrt{d})$, the complexity of sample amplification reduces to $n=\Theta(\sqrt{d})$, the case without the knowledge of $t$. The main reason why sample amplification is easier here is that the additional fake sample could be chosen as the first symbol, which has a large probability. In contrast, learning the distribution requires the estimation of all other probability masses, so the existence of a probable symbol does not help much in learning. 

\subsection{An example where the complexity of sample amplification is $\omega(\sqrt{d})$}
Next we provide an example where the complexity of sample amplification is the same as that of learning. Consider a low-rank covariance estimation model: $X_1,\cdots,X_n\sim \calN(0, \Sigma)$, where $\Sigma\in \bR^{p\times p}$ could be written as $\Sigma = UU^\top$ with $U\in \bR^{p\times d}$ and $U^\top U = I_d$. In other words, the covariance matrix $\Sigma$ is isotropic on some $d$-dimensional subspace. Here $n\ge d$ samples suffice to estimate $\Sigma$ and thus the whole distribution perfectly, for the $d$-dimensional subspace could be recovered using $d$ i.i.d. samples with probability one. Therefore, the complexity of learning the distribution is $n=d$. The following theorem states that this is also the complexity of sample amplification. 

\begin{theorem}\label{thm:lowrank_covariance}
For the above low-rank covariance estimation model with $p\ge d+1$, an $(n,n+1,0.1)$ sample amplification is possible if and only if $n\ge d$. 
\end{theorem}

Theorem \ref{thm:lowrank_covariance} shows that as opposed to learning, sample amplification fails to exploit the low-rank structure in the covariance estimation problem. As a result, the complexity of sample amplification coincides with that of learning in this example. Note that sample amplification is always no harder than learning: the learner could always estimate the distribution, generate one observation from the distribution and append it to the original samples. Therefore, Theorem \ref{thm:lowrank_covariance} provides an example where the relationship between sample amplification and learning is the worst possible. 

\subsection{An example where the TV distance is not the right metric}
Finally we provide an example showing that the TV distance is not the right metric for the learning-based approach in Section \ref{sec.shuffling}, and thereby partially illustrate the necessity of using the $\chi^2$ divergence. This example also goes beyond parametric families for sample amplification. Let $\calP$ be the class of all $L$-Lipschitz densities supported on $[0,1]$, i.e. the density $f$ satisfies $|f(x)-f(y)|\le L|x-y|$ for all $x,y\in [0,1]$. For $c\in (0,1)$, also let $\calP_c\subseteq \calP$ be the subclass of densities lower bounded by $c$ everywhere, i.e. $f(x)\ge c$ for all $x\in [0,1]$. It is a classical result (see, e.g. \cite{tsybakov2009introduction}) that the minimax density estimation error under TV distance is $\Theta(n^{-1/3})$ for both $\calP$ and $\calP_c$. The next theorem shows that the sample complexities for amplification are actually different.

\begin{theorem}\label{thm:lipschitz}
Let $L\ge 8$ and $c\in (0,1)$ be fixed. It holds that 
\begin{align*}
m^\star(\calP_c, n) \asymp n^{5/6}, \quad \text{ while } \quad m^\star(\calP, n) \lesssim n^{3/4}. 
\end{align*}
\end{theorem}

Theorem \ref{thm:lipschitz} shows that, although assuming a density lower bound does not alter the TV estimation error, it boosts the size of amplified samples from $O(n^{3/4})$ to $\Theta(n^{5/6})$. In fact, the $\chi^2$-estimation error is also reduced from $\calP$ to $\calP_c$: in the proof of Theorem \ref{thm:lipschitz} we show that $r_{\chi^2}(\calP_c, n)\lesssim n^{-2/3}$, but $m^\star(\calP, n) \lesssim n^{3/4}$ together with Theorem \ref{thm:shuffling_general} imply that $r_{\chi^2}(\calP, n)\gtrsim n^{-1/2}$. Therefore, this is an example suggesting that measuring the estimation error under the $\chi^2$ divergence might be a better indicator for the complexity of sample amplification than the TV distance.

\paragraph{Acknowledgements.} Thank you to anonymous reviewers for helpful feedback on earlier drafts of this paper. 

\paragraph{Funding.} Shivam Garg conducted this research while affiliated with Stanford University and was supported by a Stanford Interdisciplinary Graduate Fellowship. Yanjun Han was supported by a Simons-Berkeley research fellowship and the Norbert Wiener postdoctoral fellowship in statistics at MIT IDSS. Vatsal Sharan was supported by NSF CAREER Award CCF-2239265 and an Amazon Research Award. Gregory Valiant was supported by NSF Awards AF-2341890, CCF-1704417, CCF-1813049, UT Austin's Foundation of ML NSF AI Institute, and a Simons Foundation Investigator Award. 

\appendix

\section{Concrete examples of sample amplification}\label{append:examples}
In this section we include concrete examples of sample amplification omitted in the main text due to space limitations, including the non-asymptotic upper bounds in \Cref{sec.sufficiency} and \ref{sec.shuffling}, and the lower bounds for non-product models in \Cref{sec:lower_bound}. 

\subsection{Concrete examples of amplification via sufficiency}\label{subsec:example_sufficiency}
In this section, in contrast to the mostly asymptotic results in Section \ref{subsec:exponential_family}, we investigate several non-asymptotic examples of sample amplification in concrete models. We show that for many natural models, including exponential family with dependent coordinates and non-exponential family which are not covered in the general theory, the sufficiency-based sample amplification approach could still amplify $\Omega(n/\sqrt{d})$ additional samples. We also illustrate the computational idea in Section \ref{subsec:computation} via more involved examples. 

Our first example concerns the Gaussian model with a known mean but an unknown covariance. It is a folklore result that estimating the unknown covariance in a vanishing Frobenius norm requires $n=\Omega(d^2)$ samples \cite[Corollary 1.2]{devroye2020minimax}, which is also the sample complexity for learning the distribution within a small TV distance. The following example shows that $n=O(d)$ samples suffice for sample amplification. 
\begin{example}[Gaussian covariance model with known mean]\label{example:gaussian_covariance}
Consider the i.i.d. observations $X_1,\cdots,X_n$ drawn from $\calN(0,\Sigma)$ with zero mean and an unknown covariance $\Sigma\in \bR^{d\times d}$. Here a minimal sufficient statistic is the sample covariance matrix
\begin{align*}
\widehat{\Sigma}_n = \frac{1}{n}\sum_{i=1}^n X_iX_i^\top. 
\end{align*}
Lemma \ref{lemma:gaussian_covariance} shows that $\|\calL(\widehat{\Sigma}_n) - \calL(\widehat{\Sigma}_{n+m})\|_{\text{\rm TV}} \le \varepsilon$ as long as $m = O(\varepsilon n/d)$, therefore drawing samples from $P_{X^{n+m}\mid \widehat{\Sigma}_{n+m} }$ achieves sample amplification of size $m = \Omega(n/d)$. This coincides with the general $O(n/\sqrt{D})$ result where $D\asymp d^2$ is the parameter dimension. 

In order to sample from this conditional distribution, consider the following statistic
\begin{align*}
S_{n+m} = [(n+m)\widehat{\Sigma}_{n+m}]^{-1/2}[X_1, X_2, \cdots, X_{n+m}] \in \bR^{d\times (n+m)}, 
\end{align*}
which always exists even if $\widehat{\Sigma}_{n+m}$ is not invertible. Clearly there is a bijection between $X^{n+m}$ and $(\widehat{\Sigma}_{m+n}, S_{n+m})$, and Lemma \ref{lemma:gaussian_covariance} shows that $S_{n+m}$ is an ancillary statistic. In particular, $S_{n+m}$ always follows the uniform distribution on the following set: 
\begin{align*}
A = \{U\in \bR^{d\times (n+m)}: UU^\top = I_d \}. 
\end{align*}
Consequently  an $(n, n+\Omega(n\varepsilon/d),\varepsilon)$ sample amplification is efficiently achievable in the Gaussian covariance model using the following algorithm:
\begin{enumerate}
    \item Given samples $[X_1, X_2, \cdots, X_{n}]$ compute $\widehat{\Sigma}_n = \frac{1}{n}\sum_{i=1}^n X_iX_i^\top$. 
    \item Sample $[Z_1, Z_2, \cdots, Z_{n+m}]$ where $Z_i\sim N(0,I)$ for all $i$. Compute $\widehat{\Sigma}_{n+m} = \frac{1}{n}\sum_{i=1}^n Z_iZ_i^\top$.  Based on these compute $S_{n+m} = [(n+m)\widehat{\Sigma}_{n+m}]^{-1/2}[Z_1, Z_2, \cdots, Z_{n+m}]$.
    \item Output $(n+m)$ samples $X^{n+m}=[(n+m)\widehat{\Sigma}_{n}] ^{1/2}S_{n+m}$.
\end{enumerate}


\end{example}
\begin{lemma}\label{lemma:gaussian_covariance}
Under the notations of Example \ref{example:gaussian_covariance}, for $n\ge 4\max\{m,d\}$ it holds that
\begin{align*}
\| \calL(\widehat{\Sigma}_n) - \calL(\widehat{\Sigma}_{n+m}) \|_{\text{\rm TV}} \le \frac{2md}{n}. 
\end{align*}
In addition, $S_{n+m}$ is uniformly distributed on $A = \{U\in \bR^{d\times (n+m)}: UU^\top = I_d \}$. 
\end{lemma}

The next example shows that the sample amplification result does not change much if both the mean and covariance are unknown, though the sampling procedure becomes slightly more involved.
\begin{example}[Gaussian model with unknown mean and covariance]\label{example:mean_covariance}
Next we consider the most general Gaussian model where $X_1,\cdots,X_n\sim \calN(\theta,\Sigma)$ with unknown mean vector $\theta\in\bR^d$ and unknown covariance matrix $\Sigma\in \bR^{d\times d}$. In this case, a minimal sufficient statistic is the pair $(\overline{X}_n, \widehat{\Sigma}_n)$, with
\begin{align*}
\overline{X}_n = \frac{1}{n}\sum_{i=1}^n X_i, \qquad \widehat{\Sigma}_n = \frac{1}{n-1}\sum_{i=1}^n (X_i - \overline{X}_n)(X_i - \overline{X}_n)^\top.  
\end{align*}
Lemma \ref{lemma:gaussian_mean_covariance} shows that $\|\calL(\overline{X}_n, \widehat{\Sigma}_n) - \calL(\overline{X}_{n+m}, \widehat{\Sigma}_{n+m})\|_{\text{\rm TV}} \le \varepsilon$ as long as $m = O(n\varepsilon/d)$, therefore drawing amplified samples from $P_{X^{n+m}\mid (\overline{X}_{n+m}, \widehat{\Sigma}_{n+m})}$ achieves sample amplification of size $m = \Omega(n\varepsilon/d)$. For the computation, consider the following statistic
\begin{align*}
S_{n+m} = [(n+m-1)\widehat{\Sigma}_{n+m}]^{-1/2}[X_1-\overline{X}_{n+m}, \cdots, X_{n+m} - \overline{X}_{n+m}] \in \bR^{d\times (n+m)}, 
\end{align*}
and it is clear that the whole samples $X^{n+m}$ could be recovered from $(\overline{X}_{n+m}, \widehat{\Sigma}_{n+m}, S_{n+m})$. Again, Lemma \ref{lemma:gaussian_mean_covariance} shows that $S_{n+m}$ is ancillary and uniformly distributed on the following set (assuming $n+m-1\ge d$): 
\begin{align*}
A = \{U\in \bR^{d\times (n+m)}: UU^\top = I_d, U{\bf 1} = {\bf 0} \}.
\end{align*}
\end{example}

\begin{lemma}\label{lemma:gaussian_mean_covariance}
	Under the notations of Example \ref{example:mean_covariance}, for $n\ge 4\max\{m,d\}$ it holds that
	\begin{align*}
	\| \calL(\overline{X}_n, \widehat{\Sigma}_n) - \calL(\overline{X}_{n+m}, \widehat{\Sigma}_{n+m}) \|_{\text{\rm TV}} \le \frac{3md}{n-1}. 
	\end{align*}
	In addition, $S_{n+m}$ is uniformly distributed on $A = \{U\in \bR^{d\times (n+m)}: UU^\top = I_d, U{\bf 1} = {\bf 0} \}$. 
\end{lemma}

The following example concerns the product exponential distributions, or general Gamma distributions with a fixed shape parameter. 
\begin{example}[Product exponential distribution]\label{example:exponential}
In this example, we consider the product exponential model where $X_1,\cdots,X_n \sim \prod_{i=1}^d \text{\rm Exp}(\lambda_i)$ with unknown rate vector $\lambda=(\lambda_1,\cdots,\lambda_d)$. Again, in this model the sample mean $\overline{X}_n$ is sufficient, and follows a product Gamma distribution $\prod_{i=1}^d\text{\rm Gamma}(n,n\lambda_i)$. Consequently, 
\begin{align*}
D_{\text{\rm KL}}( \calL(\overline{X}_n) \| \calL(\overline{X}_{n+m}) ) = d\cdot \left( m -  (n+m)\log\left(1+\frac{m}{n}\right) + \log\frac{\Gamma(n+m)}{\Gamma(n)} - m\psi(n) \right), 
\end{align*}
where $\Gamma(x)=\int_0^\infty t^{x-1}e^{-t}dt$ and $\psi(x)=\frac{d}{dx}[\log \Gamma(x)]$ are the gamma and digamma functions, respectively. Using the definition of $f(n,m,d)$ in \eqref{eq:Wishart_KL}, we note that the above KL divergence is precisely $d\cdot f(2n,2m,1)$, thus by the proof of Lemma \ref{lemma:gaussian_covariance} is further at most $O(dm^2/n^2)$. Consequently, sample amplification is possible as long as $n = \Omega(\sqrt{d}/\varepsilon)$ and $m=O(n\varepsilon/\sqrt{d})$. Alternatively, the same result could also follow from Theorem \ref{thm:upper_exp_product}, for both assumptions \ref{assump.continuity} and \ref{assump.moment} hold for the exponential distribution. 

To draw amplified samples $X_1,\cdots,X_{n+m}$ conditioned on $\overline{X}_{n+m}$, note that the following statistic $S\in \bR^{d\times(n+m)}$ with $i$-th row being
\begin{align*}
S_i = \left(\frac{X_{1,i}}{\overline{X}_{n+m,i}}, \frac{X_{2,i}}{\overline{X}_{n+m,i}}, \cdots, \frac{X_{n+m,i}}{\overline{X}_{n+m,i}} \right), 
\end{align*}
is ancillary. In fact, each $S_i$ follows the Dirichlet distribution $\text{\rm Dir}(1,1,\cdots,1)$. Then computational efficiency follows from the obvious fact that $X^{n+m}$ is determined by $(\overline{X}_{n+m}, S)$. 
\end{example}

Our final example is a non-exponential family which is not even differentiable in quadratic mean, but the $n=O(\sqrt{d})$ sample complexity still holds.
\begin{example}[Uniform distribution over a rectangle]\label{example:uniform}
In this example, let $X_1,\cdots,X_n$ be i.i.d. samples from the uniform distribution on an unknown rectangle $\prod_{j=1}^d [a_j, b_j]$ in $\bR^d$. Note that this is not an exponential family. However, the sufficiency-based sample amplification could still be applied in this case. Specifically, here the sufficient statistics are $X_{\min}^n$ and $X_{\max}^n$, where for $j\in [d]$, 
\begin{align*}
X_{\min,j}^n = \min_{i\in [n]} X_{i,j}, \qquad X_{\max,j}^n = \max_{i\in [n]} X_{i,j}. 
\end{align*}
It is not hard to find that, the joint density of $(X_{\min,j}^n, X_{\max,j}^n)$ is
\begin{align*}
f_n(a,b) = \frac{n(n-1)(b-a)^{n-2}}{(b_j - a_j)^{n}}\cdot \1(a_j\le a\le b\le b_j). 
\end{align*}
Then after some algebra, the {\rm KL} divergence between the sufficient statistics is 
\begin{align*}
D_{\text{\rm KL}}(\calL(X_{\min}^n, X_{\max}^n) \|  \calL(X_{\min}^{n+m}, X_{\max}^{n+m})) = d\left(\frac{m}{n} - \log\left(1+\frac{m}{n}\right) + \frac{m}{n-1} - \log\left(1+\frac{m}{n-1}\right)\right), 
\end{align*}
which is $O(dm^2/n^2)$ when $m = O(n)$. Therefore, sample amplification for uniform distributions is still possible whenever $m = O(n\varepsilon/\sqrt{d})$ and $n = \Omega(\sqrt{d}/\varepsilon)$, same as exponential families. 

To draw amplified samples $X_1,\cdots,X_{m+n}$ conditioned on $(X_{\min}^{n+m}, X_{\max}^{m+n})$, note that the following statistic $S\in \bR^{d\times (n+m)}$ with $i$-th row being 
\begin{align*}
S_i = \left( \frac{X_{1,i} - X_{\min}^{n+m} }{ X_{\max}^{m+n} - X_{\min}^{m+n} }, \cdots, \frac{X_{m+n,i} - X_{\min}^{n+m} }{ X_{\max}^{m+n} - X_{\min}^{m+n} } \right),
\end{align*}
is ancillary. This is due to the invariance property of the uniform distribution: if $X\sim \mathsf{U}(0,1)$, then $aX+b\sim \mathsf{U}(b,a+b)$. Since $(X_1,\cdots,X_{n+m})$ is determined by $(X_{\min}^{n+m}, X_{\max}^{n+m}, S)$, the computational efficiency follows. 
\end{example}

Although all of the above examples work exclusively for continuous models and apply an identity map between sufficient statistics, we remark that more sophisticated maps based on learning could be useful and work for discrete models. The idea of sample amplification via learning is presented in Section \ref{sec.shuffling}, and an example which combines both the sufficiency and learning ideas could be found in Example \ref{example:poisson}. 

\subsection{Concrete examples of amplification via learning}\label{subsec:example_shuffling}
In this section we show how learning-based approaches achieve optimal performances of sample amplification in several examples, including both continuous and discrete models. In some scenarios the following strengthened lemma may be useful to deal with too large $\chi^2$-divergence. 

\begin{lemma}\label{lemma:modified_chi2}
The same results of Theorem \ref{thm:shuffling_general} and \ref{thm:shuffling_product} hold for the following modification of the $\chi^2$-estimation error: 
\begin{align*}
r_{\bar{\chi}^2}(\calP,n) = \inf_{\widehat{P}_n}\sup_{P\in \calP} \bE_P\left[\chi^2(\widehat{P}_n,P) \wedge n \right]. 
\end{align*}
\end{lemma}

The idea behind Lemma \ref{lemma:modified_chi2} is that for some models, it might happen that $\chi^2(\widehat{P}_n,P)=\infty$ with a small probability. However, for the TV distance we always have $\| \widehat{P}_n - P \|_{\text{TV}}\le 1$, so a large $\chi^2$-divergence could still lead to a meaningful TV distance. The proof of Lemma \ref{lemma:modified_chi2} could be found in Appendix \ref{appendix:lemma}. 

The first example is again the Gaussian location model in Example \ref{example:mean}, where we show that the shuffling-based approach also achieves the complexity $n=O(\sqrt{d}/\varepsilon)$ and the size $m=\Omega(n\varepsilon/\sqrt{d})$. 

\begin{example}[Gaussian location model with known covariance, continued]\label{example:mean_shuffle}
Consider the setting of Example \ref{example:mean}, where the family of distributions is $\calP = \{\calN(\theta,I_d)\}_{\theta\in \bR^d}$. Here $\calP$ has a product structure, with $\calP_j = \{ \calN(\theta_j,1)\}_{\theta_j\in \bR}$ for each $j\in [d]$. To find an upper bound on the $\chi^2$-estimation error, consider the distribution estimator $\widehat{P}_{n,j} = \calN(\widehat{\theta}_j,1)$, where $\widehat{\theta}_j = n^{-1}\sum_{i=1}^n X_{i,j}$. Consequently, 
\begin{align*}
\chi^2( \widehat{P}_{n,j} , P_j ) = \chi^2( \calN(\widehat{\theta}_j, 1) , \calN(\theta_j,1) ) = \exp( (\widehat{\theta}_j - \theta_j)^2 ) - 1,
\end{align*}
and using $\widehat{\theta}_j - \theta_j \sim \calN(0,1/n)$, we have
\begin{align*}
\bE[ \chi^2(\widehat{P}_{n,j}, P_j) ] = \sqrt{\frac{n}{n-2}} - 1 = O\left(\frac{1}{n}\right)
\end{align*}
whenever $n\ge 3$. Consequently, $r_{\chi^2}(\calP_j,n)= O(1/n)$ for all $j\in [d]$, and Theorem \ref{thm:shuffling_product} implies that an $(n,n+m,\varepsilon)$ sample amplification is possible if $n=\Omega(\sqrt{d}/\varepsilon)$ and $m=O(n\varepsilon/\sqrt{d})$. 
\end{example}

The next example is the discrete distribution model considered in \cite{axelrod2019sample}. 
\begin{example}[Discrete distribution model]\label{example:shuffle_discrete}
Let $\calP$ be the class of all discrete distributions supported on $k$ elements. In this case, a natural learner is the empirical distribution $\widehat{P}_n = (\widehat{p}_1,\cdots,\widehat{p}_k)$, with
\begin{align*}
\widehat{p}_j = \frac{1}{n}\sum_{i=1}^n \1(X_i = j), \qquad j\in [k]. 
\end{align*}
Consequently,
\begin{align*}
\bE[ \chi^2( \widehat{P}_n, P) ] = \bE\left[\sum_{j=1}^k \frac{(\widehat{p}_j - p_j)^2}{p_j} \right] = \sum_{j=1}^k \frac{1-p_j}{n} = \frac{k-1}{n},
\end{align*}
meaning that $r_{\chi^2}(\calP,n)\le (k-1)/n$. Hence, Theorem \ref{thm:shuffling_general} implies that sample amplification is possible whenever $n=\Omega(\sqrt{k}/\varepsilon)$ and $m=O(n\varepsilon/\sqrt{k})$. In this case, Algorithm \ref{alg:shuffling_general} essentially subsamples from the original data and add them back with random shuffling, which is the same algorithm in \cite{axelrod2019sample}. However, the analysis here is much simplified. 
\end{example}

The following example revisits the uniform distribution in Example \ref{example:uniform}, and again shows that the same amplification performance could be achieved by random shuffling. 

\begin{example}[Uniform distribution, continued]\label{example:shuffle_uniform}
Consider the setting in Example \ref{example:uniform}, where the distribution class $\calP$ is the family of all uniform distributions on a rectangle $\prod_{j=1}^d [a_j,b_j]$. For each $j$, a natural distribution estimator $\widehat{P}_{n,j}$ is simply the uniform distribution on $[X_{\min,j}^n, X_{\max,j}^n]$, where these quantities are defined in Example \ref{example:uniform}. For this learner, it holds that
\begin{align*}
\chi^2\left( \widehat{P}_{n,j}, P_j \right) = \frac{(a_j-b_j)^2}{(X_{\max,j}^n - X_{\min,j}^n)^2} - 1. 
\end{align*}
Note that Example \ref{example:uniform} shows that the joint density of $(X_{\min,j}^n, X_{\max,j}^n)$ is given by
\begin{align*}
f_n(a,b) = \frac{n(n-1)(b-a)^{n-2}}{(b_j - a_j)^{n}}\cdot \1(a_j\le a\le b\le b_j), 
\end{align*}
the expected $\chi^2$-divergence could then be computed as
\begin{align*}
\bE\left[\chi^2\left( \widehat{P}_{n,j}, P_j \right) \right] &= -1 + \iint_{a_j\le a\le b\le b_j} \frac{n(n-1)(b-a)^{n-4}}{(b_j-a_j)^{n-2}} dadb = \frac{4n-6}{(n-2)(n-3)}, 
\end{align*}
which is $O(n^{-1})$ if $n\ge 4$. Hence, we have $r_{\chi^2}(\calP_j,n) = O(1/n)$ for each $j\in [d]$, and Theorem \ref{thm:shuffling_product} shows an $(n,n+m,\varepsilon)$ sample amplification if $n=\Omega(\sqrt{d}/\varepsilon)$ and $m=O(n\varepsilon/\sqrt{d})$. 
\end{example}

The next example is the exponential distribution model studied in Example \ref{example:exponential}, where the modified $\chi^2$-learning error $r_{\bar{\chi}^2}(\calP,n)$ in Lemma \ref{lemma:modified_chi2} will be useful. 

\begin{example}[Exponential distribution, continued]\label{example:shuffle_exponential}
Consider the setting in Example \ref{example:exponential}, where the distribution class $\calP$ is a product of exponential distributions with unknown rate parameters. In this case, a natural distribution learner is to estimate each rate parameter as $\widehat{\lambda}_{n,j} = n/\sum_{i=1}^n X_{i,j}$, and use $\text{\rm Exp}(\widehat{\lambda}_{n,j})$ to estimate the truth $\text{\rm Exp}(\lambda_j)$. Note that
\begin{align*}
\chi^2\left(\text{\rm Exp}(\widehat{\lambda}_{n,j}), \text{\rm Exp}(\lambda_j) \right) = \frac{(\widehat{\lambda}_{n,j} - \lambda_j)^2}{\lambda_j(2\widehat{\lambda}_{n,j} - \lambda_j)}
\end{align*}
whenever $2\widehat{\lambda}_{n,j}>\lambda_j$. However, if $2\widehat{\lambda}_{n,j} \le \lambda_j$ the $\chi^2$-divergence will be unbounded, which occurs with a small but positive probability.

To address this issue, we note that when $\lambda_j=1$, the sub-exponential concentration claims that $|\sum_{i=1}^n X_{i,j} - n|\le n/3$ with probability at least $1-\exp(-\Omega(n))$. By a simple scaling, the above event implies that $\widehat{\lambda}_{n,j}/\lambda_j \in [3/4,3/2]$. Hence, 
\begin{align*}
\bE\left[\chi^2\left(\text{\rm Exp}(\widehat{\lambda}_{n,j}), \text{\rm Exp}(\lambda_j) \right)\wedge n \right] \le 2\cdot \bE\left[ \left(\frac{\widehat{\lambda}_{n,j} }{\lambda_j} -1\right)^2 \right] + n\cdot \exp(-\Omega(n)) = O\left(\frac{1}{n}\right), 
\end{align*}
which means that $r_{\bar{\chi}^2}(\calP_j,n) = O(1/n)$. Therefore, by Lemma \ref{lemma:modified_chi2}, the sample amplification is possible whenever $n=\Omega(\sqrt{d}/\varepsilon)$ and $m=O(n\varepsilon/\sqrt{d})$, the same as Example \ref{example:exponential}. 
\end{example}

The following example considers an interesting non-product model, i.e. the Gaussian distribution with a sparse mean vector. 

\begin{example}[Sparse Gaussian model]\label{example:shuffle_sparse_gaussian}
Consider the Gaussian location model $\calP = \{\calN(\theta,I_d) \}_{\theta\in\Theta}$, with an additional constraint that the mean vector $\theta$ is $s$-sparse, i.e. $\Theta = \{\theta\in \bR^d: \|\theta\|_0\le s\}$. For the learning problem, it is well-known (cf. \cite[Theorem 1]{donoho1994ideal}) that the soft-thresholding estimator $\widehat{\theta}_n$ with
\begin{align*}
\widehat{\theta}_{n,j} = \text{\rm sign}\left(\frac{1}{n}\sum_{i=1}^n X_{i,j} \right)\cdot \left( \left|\frac{1}{n}\sum_{i=1}^n X_{i,j}\right| - \sqrt{\frac{C\log n}{n}}\right)_+
\end{align*}
and any constant $C> 2$ achieves that $\sup_{\theta\in\Theta} \bE[\|\widehat{\theta}_n - \theta \|_2^2] = O(s\log d/n)$. Therefore, the sample complexity of learning sparse Gaussian distributions is $n=O(s\log d)$. 

For the complexity of sample amplification, for each $j\in [d]$ we apply $\calN(\widehat{\theta}_{n,j}, 1)$ as the distribution estimator. The $\chi^2$-estimation performance of this estimator is summarized in Lemma \ref{lmm:sparse_gaussian_shuffle}. Therefore, by a simple adaptation of Lemma \ref{lemma:modified_chi2}, an $(n,n+m)$ sample amplification is possible as long as $n=\Omega(\sqrt{s\log d}/\varepsilon)$ and $m = O(n\varepsilon/\sqrt{s\log d})$. 
\end{example} 

\begin{lemma}\label{lmm:sparse_gaussian_shuffle}
Under the setting of Example \ref{example:shuffle_sparse_gaussian}, it holds that 
\begin{align*}
\sup_{\theta\in\Theta} \sum_{j=1}^d \bE\left[ \chi^2\left( \calN(\widehat{\theta}_{n,j}, 1), \calN(\theta_j,1) \right) \wedge n \right] = O\left( \frac{s\log d}{n}\right). 
\end{align*}
\end{lemma}

The final example is a special example where the sole application of either sufficient statistics or learning will fail to achieve the optimal sample amplification. The solution is to use a combination of both ideas. 
\begin{example}[Poisson distribution]\label{example:poisson}
Consider the product Poisson model $\prod_{j=1}^d \mathsf{Poi}(\lambda_j)$ with $\lambda\in \bR_+^d$. We first show that a na\"{i}ve application of either sufficiency-based or shuffling-based idea will not lead to a desired sample amplification. In fact, the sufficient statistic here is $T_n = \sum_{i=1}^n X_i$ which follows a product Poisson distribution $\prod_{j=1}^d \mathsf{Poi}(n\lambda_j)$; as $T_n$ takes discrete values, applying any linear map between $T_n$ to $T_{n+m}$ will not result in a small TV distance between sufficient statistics. 

The argument for the shuffling-based approach is subtler. A natural distribution estimator for $\mathsf{Poi}(\lambda_j)$ is $\mathsf{Poi}(\widehat{\lambda}_{n,j})$, with $\widehat{\lambda}_{n,j}=n^{-1}\sum_{i=1}^n X_{i,j}$. Lemma \ref{lmm:shuffle_poisson} shows that this distribution estimator could suffer from an expected $\chi^2$-estimation error far greater than $\Omega(1/n)$, so we cannot conclude an $(n, n+ \Omega(n\varepsilon/\sqrt{d}), \varepsilon)$ sample amplification from Theorem \ref{thm:shuffling_product}. 

Now we show that a combination of the sufficient statistic and learning leads to the rate-optimal sample amplification in this model. Specifically, we split samples and compute the empirical rate parameter $\widehat{\lambda}_{n/2}$ based on the first $n/2$ samples. Next, conditioned on the first half of samples, the sufficient statistic for the remaining half is $T_{n/2} = \sum_{i=n/2+1}^n X_i$. Define 
\begin{align*}
\widehat{T}_{n/2+m} = T_{n/2} + Z, 
\end{align*}
where $Z\sim \prod_{j=1}^d \mathsf{Poi}(m\widehat{\lambda}_{n/2,j})$ is independent of $T_{n/2}$ conditioning on the first half samples. Finally, we generate the $(n/2+m)$ amplified samples from the conditional distribution and append them to $(X_1,\cdots,X_{n/2})$. By the second statement of Lemma \ref{lmm:shuffle_poisson}, this achieves an $(n,n+m,\varepsilon)$ sample amplification whenever $n=\Omega(\sqrt{d}/\varepsilon)$ and $m=O(n\varepsilon/\sqrt{d})$.
\end{example}

\begin{lemma}\label{lmm:shuffle_poisson}
Under the settings of Example \ref{example:poisson}, there exists some $\lambda_j>0$ such that
\begin{align*}
    \bE\left[ \chi^2\left( \mathsf{Poi}(\widehat{\lambda}_{n,j}) , \mathsf{Poi}(\lambda_j) \right) \wedge n \right] = \Omega\left(\frac{1}{\log n}\right) \gg \frac{1}{n}. 
\end{align*}
In addition, the proposed sufficiency+learning approach satisfies
\begin{align*}
    \| \calL(\widehat{X}^{n+m}) - \calL(X^{n+m}) \|_{\text{\rm TV}} \le \frac{m\sqrt{2d}}{n}. 
\end{align*}
\end{lemma}

\subsection{Non-asymptotic examples of lower bounds}\label{subsec:lower_example}
In this section, we apply the lower bound ideas to several concrete examples in Section \ref{sec.sufficiency} and \ref{sec.shuffling}, and show that the previously established upper bounds for sample amplification are indeed tight. Since Theorem \ref{thm:lower_product} already handles all product models (including the Gaussian location model in Example \ref{example:mean}, \ref{example:mean_computation}, and \ref{example:mean_shuffle}, exponential model in Example \ref{example:exponential} and \ref{example:shuffle_exponential}, uniform model in Example \ref{example:uniform} and \ref{example:shuffle_uniform}, and Poisson model in Example \ref{example:poisson}), we are only left with the remaining non-product models. In the sequel, we will make use of the general Lemma \ref{lemma:general_lower_bound} to prove non-asymptotic lower bounds in these examples.

The lower bound of the discrete distribution model was obtained in \cite{axelrod2019sample}, where the learning approach in Example \ref{example:shuffle_discrete} is rate optimal. Our first example concerns the ``Poissonized'' version of the discrete distribution model, and we show that the results become slightly different. 

\begin{example}[``Poissonized'' discrete distribution model]\label{example:discrete_poisson}
We consider the following ``Poissonized'' discrete distribution model, where we have $n$ i.i.d. samples drawn from $\prod_{j=1}^k \mathsf{Poi}(p_j)$, with $(p_1,\cdots,p_k)$ being an unknown probability vector. Although the Poissonization does not affect the optimal rate of estimation in many problems, Lemma \ref{lmm:lower_poissonization} shows the following distinction when it comes to sample amplification: the optimal amplification size is $m=\Theta(n\varepsilon/\sqrt{k}+\sqrt{n}\varepsilon)$ under the Poissonized model, while it is $m=\Theta(n\varepsilon/\sqrt{k})$ under the non-Poissonized model \cite{axelrod2019sample}. 

The complete proof of Lemma \ref{lmm:lower_poissonization} is relegated to the appendix, but we briefly comment on why Theorem \ref{thm:lower_product} is not directly applicable when $k\gg n$. The reason here is that to apply Theorem \ref{thm:lower_product}, we will construct a parametric submodel which is a product model: 
\begin{align*}
P_{\theta} = \prod_{j=1}^{k_0} \left[\mathsf{Poi}\left(\frac{1}{k} + \theta_j \right) \times \mathsf{Poi}\left(\frac{1}{k} - \theta_j \right)\right], \qquad \theta \in \Theta \triangleq \left[-\frac{1}{k}, \frac{1}{k}\right]^{k_0},
\end{align*}
where $k=2k_0$. However, for $\theta, \theta'\in [0,1/k]$, the range of the squared Hellinger distance
\begin{align*}
H^2\left(\mathsf{Poi}\left(\frac{1}{k} + \theta \right) \times \mathsf{Poi}\left(\frac{1}{k} - \theta \right), \mathsf{Poi}\left(\frac{1}{k} + \theta' \right) \times \mathsf{Poi}\left(\frac{1}{k} - \theta' \right)\right)
\end{align*}
is only $[0, \Theta(1/k)]$, so Assumption \ref{assump.hellinger} does not hold when $k\gg n$. This is precisely the subtlety in the Poissonized model. 
\end{example}

\begin{lemma}\label{lmm:lower_poissonization}
Under the Poissonized discrete distribution model, an $(n,n+m,\varepsilon)$ sample amplification is possible if and only if $n=\Omega(1)$ and $m=O(n\varepsilon/\sqrt{k}+\sqrt{n}\varepsilon)$. 
\end{lemma}

Our next example is the sparse Gaussian model in Example \ref{example:shuffle_sparse_gaussian}, where we establish the tightness of $n=\Omega(\sqrt{s\log d}/\varepsilon)$ and $m=O(n\varepsilon/\sqrt{s\log d})$ directly using Lemma \ref{lemma:general_lower_bound}. 

\begin{example}[Sparse Gaussian location model]\label{example:lower_sparse_gaussian}
Consider the setting of the sparse Gaussian location model in Example \ref{example:shuffle_sparse_gaussian}, and we aim to prove a matching lower bound for sample amplification. In the sequel, we first handle the case $s=1$ to reflect the main idea, and postpone the case of general $s$ to Lemma \ref{lemma:gaussian_prob}.  

For $s=1$, we apply Lemma \ref{lemma:general_lower_bound} to a proper choice of the prior $\mu$ and loss $L$. Fixing a parameter $t>0$ to be chosen later, let $\mu$ be the uniform distribution on the finite set of vectors $\{te_1, \cdots, te_d \}$, where $e_1,\cdots,e_d$ are the canonical vectors in $\bR^d$. Moreover, for an estimator $\widehat{\theta}\in \bR^d$ of the unknown mean vector, the loss function is chosen to be $L(\theta,\widehat{\theta}) = \1(\widehat{\theta} \neq \theta)$. For the above prior and loss, it is clear that the maximum likelihood estimator $\widehat{\theta} = te_{\widehat{j}}$ with
\begin{align*}
\widehat{j} = \arg\max_{j\in [d]} \sum_{i=1}^n X_{i,j}
\end{align*}
is the Bayes estimator, and the Bayes risk admits the following expression: 
\begin{align*}
r_{\text{\rm B}}(\calP,n,\mu,L) = 1 - \bE\left[\frac{\exp(\sqrt{n}t(\sqrt{n}t+Z_1))}{\exp(\sqrt{n}t(\sqrt{n}t+Z_1)) + \sum_{j=2}^d \exp(\sqrt{n}tZ_j) } \right] \triangleq 1 - p_d(\sqrt{n}\cdot t), 
\end{align*}
where $Z_1,\cdots,Z_d\sim \calN(0,1)$ are i.i.d. standard normal random variables. Similarly, the Bayes risk under $n+m$ samples is $r_{\text{\rm B}}(\calP,n+m,\mu,L) = 1 - p_d(\sqrt{n+m}\cdot t)$, and it remains to investigate the property of the function $p_d(\cdot)$. Lemma \ref{lemma:gaussian_prob} summarizes a useful property for the lower bound, i.e. the function $p_d(z)$ enjoys a phase transition around $z\sim \sqrt{2\log d}$. 

Based on Lemma \ref{lemma:gaussian_prob}, the pigeonhole principle implies that for every given $c>0$, there exists some $z\in [\sqrt{2\log d}-C, \sqrt{2\log d}+C]$ such that $p_d(z+c\varepsilon) - p_d(z) = \Omega_c(\varepsilon)$. Therefore, if $m=\lceil cn\varepsilon/\sqrt{\log d} \rceil$, choosing $t = z/\sqrt{n}$ for the above $z$ yields that $p_d(\sqrt{n+m}\cdot t) - p_d(\sqrt{n}\cdot t)=\Omega_c(\varepsilon)$, and consequently 
\begin{align*}
\varepsilon^\star(\calP,n,m) \ge r_{\text{\rm B}}(\calP,n,\mu,L) - r_{\text{\rm B}}(\calP,n+m,\mu,L) = p_d(\sqrt{n+m}\cdot t) - p_d(\sqrt{n}\cdot t)=\Omega_c(\varepsilon).
\end{align*}
Therefore, we must have $n=\Omega(\sqrt{\log d}/\varepsilon)$ and $m=O(n\varepsilon/\sqrt{\log d})$ for sample amplification. 
\end{example}

The following lemma summarizes the phase-transition property of $p_d$ used in the above example. 
\begin{lemma}\label{lemma:gaussian_prob}
For the function $p_d(z)$ in Example \ref{example:lower_sparse_gaussian}, there exists an absolute constant $C$ independent of $d$ such that 
\begin{align*}
p_d(\sqrt{2\log d} - C) &\le 0.1, \\
p_d(\sqrt{2\log d} + C) &\ge 0.9. 
\end{align*}
Moreover, for general $s<d/2$, an $(n,n+m)$ sample amplification is possible under the sparse Gaussian model only if $n=\Omega(\sqrt{s\log(d/s)}/\varepsilon)$ and $m=O(n\varepsilon/\sqrt{s\log(d/s)})$. 
\end{lemma}

Our final example proves the tightness of the sufficiency-based approaches in Examples \ref{example:gaussian_covariance} and \ref{example:mean_covariance} for Gaussian models with unknown covariance. The proof relies on the computation of the minimax risks in Lemma \ref{lemma:general_lower_bound}, as well as the statistical theory of invariance. 

\begin{example}[Gaussian model with unknown covariance]\label{example:lower_covariance}
We establish the matching lower bounds of sample amplification in Example \ref{example:gaussian_covariance} with a known mean, which imply the lower bounds in Example \ref{example:mean_covariance}. We make use of the minimax risk formulation of Lemma \ref{lemma:general_lower_bound}, and consider the following loss: 
\begin{align*}
L(\Sigma, \widehat{\Sigma}) = \1\left(\ell(\Sigma, \widehat{\Sigma}) \ge g(n+m-1-d,d) + C\cdot \frac{d}{n} \right), 
\end{align*}
where function $g$ is given by \eqref{eq:g_function} in the Lemma \ref{lemma:stein_loss} below, $C>0$ is a large absolute constant to be determined later, and $\ell(\Sigma, \widehat{\Sigma})$ is the Stein's loss
\begin{align*}
\ell(\Sigma, \widehat{\Sigma}) = \text{\rm tr}(\Sigma^{-1}\widehat{\Sigma}) - \log\det(\Sigma^{-1}\widehat{\Sigma}) - d. 
\end{align*}

To search for the minimax estimator under the above loss, similar arguments in \cite{james1961estimation} based on the theory of invariance show that it suffices to consider estimators taking the form
\begin{align}\label{eq:invariant_estimator}
\widehat{\Sigma}_n = L_nD_nL_n^\top, 
\end{align}
where $D_n\in \bR^{d\times d}$ is a diagonal matrix independent of the observations, and $L_n\in \bR^{d\times d}$ is the lower triangular matrix satisfying that $L_nL_n^\top = \sum_{i=1}^n X_iX_i^\top$. Moreover, the risk of the above estimator is invariant with the unknown $\Sigma$, so in the sequel we assume that $\Sigma = I_d$. Note that when $D_n = I_d/n$, the estimator $\widehat{\Sigma}_n$ is the sample covariance; however, other choices of the diagonal matrix $D_n$ could give a uniform improvement over the sample covariance, see \cite{james1961estimation}. 

The proof idea is to show that with $n+m$ samples, there exists an estimator $\widehat{\Sigma}_{n+m}$ with a specific choice of $D_{n+m}$ in \eqref{eq:invariant_estimator} such that (cf. Lemma \ref{lmm:invariant_estimator})
\begin{align*}
|\bE[\ell(\Sigma,\widehat{\Sigma}_{n+m})] - g(n+m+1-d,d)| \le \frac{5d}{n+m}, \qquad \sqrt{\var(\ell(\Sigma,\widehat{\Sigma}_{n+m}))} \le \frac{4d}{n+m}. 
\end{align*}
Consequently, by Chebyshev's inequality, for $C>0$ large enough $r(\calP,n+m,L) \le 0.1$. 

To lower bound the minimax risk $r(\calP,n,L)$ with $n$ samples, we need to enumerate over all possible estimators taking the form of \eqref{eq:invariant_estimator}. It turns out that for any choice of $D_n$, the first two moments of $\ell(\Sigma,\widehat{\Sigma}_n)$ admit explicit expressions, and it always holds that (cf. Lemma \ref{lmm:invariant_estimator})
\begin{align*}
\bE[\ell(\Sigma,\widehat{\Sigma}_{n})] \ge g(n+1-d,d) + C_1\left(\sqrt{\var(\ell(\Sigma,\widehat{\Sigma}_{n}))} - \frac{4d}{n}\right) - \frac{6d}{n}
\end{align*}
for any given constant $C_1>0$, provided that $n, d\ge C_2$ with $C_2$ depending only on $C_1$. Choosing $C_1>0$ large enough, Chebyshev's inequality leads to $\ell(\Sigma,\widehat{\Sigma}_n) \ge g(n+1-d,d) - 5C_1d/n$ with probability at least $0.9$ for every $\widehat{\Sigma}_n$ taking the form of \eqref{eq:invariant_estimator}. By the last statement of Lemma \ref{lmm:invariant_estimator}, for $m=C_3n/d$ with a large enough $C_3>0$, the above event implies that $r(\calP,n,L)\ge 0.9$.

Combining the above scenarios and applying the pigeonhole principle, an application of Lemma \ref{lemma:general_lower_bound} claims that $m=O(n\varepsilon/d)$ is necessary for an $(n,n+m,\varepsilon)$ sample amplification. 
\end{example}

The following lemma summarizes the necessary technical results for Example \ref{example:lower_covariance}. 
\begin{lemma}\label{lmm:invariant_estimator}
Let $n\ge d$. For $D_n = \text{\rm diag}(\lambda_1,\cdots,\lambda_d)$ with $\lambda_j = 1/(n+d+1-2j)$ for all $j\in [d]$, the corresponding estimator $\widehat{\Sigma}_{n}$ in \eqref{eq:invariant_estimator} satisfies
\begin{align*}
|\bE[\ell(\Sigma,\widehat{\Sigma}_{n})] - g(n+1-d,d)| \le \frac{5d}{n}, \qquad \sqrt{\var(\ell(\Sigma,\widehat{\Sigma}_{n}))} \le \frac{4d}{n}, 
\end{align*}
where 
\begin{align}\label{eq:g_function}
g(u,v) &\triangleq \frac{(u+2v)\log(u+2v) + u\log u}{2} - (u+v)\log (u+v). 
\end{align}
Meanwhile, for any choice of $D_n$ and any absolute constant $C_1>0$, the estimator $\widehat{\Sigma}_{n}$ in \eqref{eq:invariant_estimator} satisfies
\begin{align*}
\bE[\ell(\Sigma,\widehat{\Sigma}_{n})] \ge g(n+1-d,d) + C_1\left(\sqrt{\var(\ell(\Sigma,\widehat{\Sigma}_{n}))} - \frac{4d}{n}\right) - \frac{6d}{n},
\end{align*}
as long as $n, d\ge C_2$ for some large enough constant $C_2$ depending only on $C_1$. 

Finally, the function $g$ defined in \eqref{eq:g_function} satisfies the following inequality: for $n\ge 2d$ and $0\le m\le n$, it holds that
\begin{align*}
g(n+1-d,d) - g(n+m+1-d,d) \ge \frac{md^2}{13n^2}. 
\end{align*}
\end{lemma}
\section{Proof of main theorems}\label{appendix:thm}
\subsection{Proof of Theorem \ref{thm:upper_exp_general}}
We recall the following result in \cite[Theorem 2.6]{bally2016asymptotic}: given Assumptions \ref{assump.continuity} and \ref{assump.moment} with $k=3$, there exists a constant $C>0$ depending only on $d$ and the moment upper bound such that 
\begin{align*}
\sup_{\theta\in\Theta} \| \calL( \sqrt{n}[\nabla^2 A(\theta)]^{-1/2}(T_n - \nabla A(\theta)) ) - \calN(0,I_d) \|_{\text{TV}} \le \frac{C}{\sqrt{n}}. 
\end{align*}
By an affine transformation, it is then clear that
\begin{align*}
\sup_{\theta\in\Theta} \| \calL( T_n ) - \calN(\nabla A(\theta), \nabla^2 A(\theta)/n) \|_{\text{TV}} \le \frac{C}{\sqrt{n}}. 
\end{align*}
Moreover, the computation in Example \ref{example:mean} shows that 
\begin{align*}
\sup_{\theta\in\Theta} \| \calN(\nabla A(\theta), \nabla^2 A(\theta)/n) - \calN(\nabla A(\theta), \nabla^2 A(\theta)/(n+m))  \|_{\text{TV}} \le \frac{m\sqrt{d}}{n}. 
\end{align*}
Now the desired result follows from the above inequalities and a triangle inequality. 

\subsection{Proof of Theorem \ref{thm:upper_exp_product}}\label{append:edgeworth}
As the density in the Edgeworth expansion \eqref{eq:edgeworth} may be negative at some point, throughout the proof we extend the formal definition of the TV distance to any signed measures $P$ and $Q$, just as half of the $L_1$ distance. Under this new definition, it is clear that the triangle inequality $\|P-R\|_{\text{TV}} \le \|P-Q\|_{\text{TV}} + \|Q-R\|_{\text{TV}}$ still holds. The following tensorization property also holds: for general signed measures $P_1,\cdots,P_d$ and $Q_1,\cdots,Q_d$ with $\max_{i\in [d]}\max\{|P_i|(\Omega), |Q_i|(\Omega) \}\le r$, we have
\begin{align}\label{eq:tensorization}
&\| P_1\times \cdots\times P_d - Q_1 \times \cdots\times Q_d \|_{\text{TV}} \nonumber\\
&\le \sum_{i=1}^{d} \| P_1\times \cdots \times P_{i-1} \times Q_{i}\times \cdots\times Q_d - P_1\times \cdots \times P_{i} \times Q_{i+1}\times \cdots\times Q_d \|_{\text{TV}} \nonumber \\
&\le \sum_{i=1}^d \|P_i - Q_i\|_{\text{TV}}\cdot \prod_{j<i} |P_j|(\Omega) \cdot \prod_{k>i} |Q_k|(\Omega) \nonumber \\
&\le r^d\sum_{i=1}^d \|P_i - Q_i\|_{\text{TV}}. 
\end{align}

Fix any $\theta\in \Theta$, let $P_{n,i}$ (resp. $P_{n+m,i}$) be the probability distribution of the $i$-th coordinate of $T_n$ (resp. $T_{n+m}$), and $Q_{n,i}$ (resp. $Q_{n+m,i}$) be the signed measure of the corresponding Edgeworth expansion taking the form \eqref{eq:edgeworth} with $k=9$. We note that the polynomials $\calK_\ell(x)$ in \eqref{eq:edgeworth} are the same for $Q_{n,i}$ and $Q_{n+m,i}$, and their coefficients are uniformly bounded over $\theta\in\Theta$ thanks to Assumption \ref{assump.moment}. Then based on Assumptions \ref{assump.continuity} and \ref{assump.moment} with $k=10$, the result of \cite[Theorem 2.7]{bally2016asymptotic} claims that
\begin{align*}
\|P_{n,i} - Q_{n,i} \|_{\text{TV}} \le \frac{C}{n^2}, 
\end{align*}
with $C>0$ independent of $(n,d,\theta)$. Moreover, the signed measure of the Edgeworth expansion in \eqref{eq:edgeworth} could be negative only if $|x| = \Omega(\sqrt{n})$, and therefore the total variation of each $Q_{n,i}$ satisfies
\begin{align*}
|Q_{n,i}|(\bR)  = |\Gamma_{n,9}|(\bR)\le \Gamma_{n,9}([-c\sqrt{n}, c\sqrt{n}]) + |\Gamma_{n,9}|(\bR\backslash [-c\sqrt{n}, c\sqrt{n}]) \le 1 + \exp(-\Omega(n)), 
\end{align*}
where the last inequality follows from integrating the Gaussian tails. Finally, let $Q_{n,i}^+$ be the positive part of $Q_{n,i}$ in the Jordan decomposition, the tensorization property \eqref{eq:tensorization} leads to 
\begin{align}\label{eq:edgeworth_approx}
\left\|\calL(T_n) - \prod_{i=1}^d Q_{n,i}^+ \right\|_{\text{TV}} \le \left\|\calL(T_n) - \prod_{i=1}^d Q_{n,i} \right\|_{\text{TV}} \le \frac{Cd}{n^2}\cdot \left(1 + \exp(-\Omega(n)) \right)^d,
\end{align}
and similar result holds for $\calL(T_{n+m})$. 

Next it remains to upper bound the TV distance between $\prod_{i=1}^d Q_{n,i}^+$ and $\prod_{i=1}^d Q_{n+m,i}^+$. To this end, we also generalize the formal definition of the Hellinger distance between general measures which are not necessarily probabilities. Then the following inequality between generalized TV and Hellinger distance holds: for measures $P$ and $Q$ on $\Omega$,  
\begin{align}\label{eq:inequality_TV_hellinger}
\|P-Q\|_{\text{TV}} &= \frac{1}{2}\int |\mathrm{d}P-\mathrm{d}Q| = \frac{1}{2}\int |\sqrt{\mathrm{d}P}-\sqrt{\mathrm{d}Q}|(\sqrt{\mathrm{d}P}+\sqrt{\mathrm{d}Q}) \nonumber \\
&\le \frac{1}{2}\sqrt{\int (\sqrt{\mathrm{d}P}-\sqrt{\mathrm{d}Q})^2\cdot \int (\sqrt{\mathrm{d}P}+\sqrt{\mathrm{d}Q})^2} \nonumber \\
&\le H(P,Q)\cdot \sqrt{P(\Omega) + Q(\Omega)}. 
\end{align}
Also, the following tensorization property holds for the Hellinger distance between general measures: for (not necessarily probability) measures $P_i, Q_i$ on $\Omega_i$, 
\begin{align}\label{eq:tensorization_hellinger}
H^2\left(\prod_{i=1}^d P_i, \prod_{i=1}^d Q_i \right) = \frac{\prod_{i=1}^d P_i(\Omega_i) + \prod_{i=1}^d Q_i(\Omega_i)}{2} - \prod_{i=1}^d\left(\frac{P_i(\Omega_i) + Q_i(\Omega_i)}{2} - H^2(P_i, Q_i) \right). 
\end{align}
Consequently, it suffices to prove an upper bound on the Hellinger distance $H(P_i, Q_i)$, and the TV distance on the product measure is a direct consequence of \eqref{eq:inequality_TV_hellinger} and \eqref{eq:tensorization_hellinger}. 

To upper bound the individual Hellinger distance, note that after a proper affine transformation, the densities of $Q_{n,i}$ and $Q_{n+m,i}$ are as follows: 
\begin{align}
Q_{n,i}(dx) &= \gamma_n(x)\left(1 + \sum_{\ell =1}^3 \frac{\calK_{\ell,i} (\sqrt{n}\cdot x)}{n^{\ell/2}} \right)dx, \label{eq:Q_n}\\
Q_{n+m,i}(dx) &= \gamma_{n+m}(x)\left(1 + \sum_{\ell =1}^3 \frac{\calK_{\ell,i} (\sqrt{n+m}\cdot x)}{(n+m)^{\ell/2}} \right)dx, \label{eq:Q_n+m} 
\end{align}
where $\gamma_n$ is the density of $\calN(0,1/n)$, and $\calK_{\ell,i}$ is a polynomial of degree $3\ell$ with uniformly bounded coefficients. By \eqref{eq:Q_n} and \eqref{eq:Q_n+m}, we observe that there exists an absolute constant $c>0$ independent of $(n,m)$ such that $Q_{n,i}(x)/\gamma_n(x), Q_{n+m,i}(x) / \gamma_{n+m}(x)\in [1/2,3/2]$ whenever $|x|\le c$. Consequently, the squared Hellinger distance could then be expressed as
\begin{align*}
H^2(Q_{n,i}^+, Q_{n+m,i}^+) &= \frac{1}{2}\int_{|x|\le c} \left(\sqrt{Q_{n,i}(x)} - \sqrt{Q_{n+m,i}(x)}\right)^2 \mathrm{d}x + \frac{1}{2}\int_{|x|>c} \left(\sqrt{Q_{n,i}^+(x)} - \sqrt{Q_{n+m,i}^+(x)}\right)^2 \mathrm{d}x \\
&\le \int_{|x|\le c} (\sqrt{\gamma_n(x)} - \sqrt{\gamma_{n+m}(x)})^2\left(1 + \sum_{\ell =1}^3 \frac{\calK_{\ell,i} (\sqrt{n+m}\cdot x)}{(n+m)^{\ell/2}} \right) \mathrm{d}x \\
&\qquad + \int_{|x|\le c} \gamma_n(x)\left(\sqrt{1 + \sum_{\ell =1}^3 \frac{\calK_{\ell,i} (\sqrt{n}\cdot x)}{n^{\ell/2}}} - \sqrt{1 + \sum_{\ell =1}^3 \frac{\calK_{\ell,i} (\sqrt{n+m}\cdot x)}{(n+m)^{\ell/2}}} \right)^2 \mathrm{d}x \\
&\qquad + \int_{|x|>c} \left(Q_{n,i}^+(x) + Q_{n+m,i}^+(x) \right) \mathrm{d}x\\
&\equiv A_1 + A_2 + A_3. 
\end{align*}

Next we upper bound the terms $A_1, A_2$, and $A_3$ separately. 
\begin{enumerate}
	\item Upper bounding $A_1$: note that by the definition of $c$, the multiplication factor in $A_1$ is at most $3/2$. Therefore, 
	\begin{align}\label{eq:A1}
	A_1 &\le \frac{3}{2}\int_{|x|\le c} (\sqrt{\gamma_n(x)} - \sqrt{\gamma_{n+m}(x)})^2 \mathrm{d}x \nonumber \\
	&\le \frac{3}{2}\int_{\bR} (\sqrt{\gamma_n(x)} - \sqrt{\gamma_{n+m}(x)})^2 \mathrm{d}x \nonumber \\
	&\stepa{=} 3\left(1 - \left(\frac{n(n+m)}{(n+m/2)^2} \right)^{1/4}\right) \nonumber\\
	&\stepb{\le} \frac{3m^2}{4n^2}, 
	\end{align}
	where (a) is due to the direct computation of the squared Hellinger distance between $\calN(0,1/n)$ and $\calN(0,1/(n+m))$, and (b) makes use of the inequality $1-x^{1/4} \le 1-x$ for $x\in [0,1]$. 
	\item Upper bounding $A_2$: using $|\sqrt{a}-\sqrt{b}|\le |a-b|$ for $a,b\ge 1/2$, we have
	\begin{align*}
	A_2 &\le \int_{|x|\le c} \gamma_n(x)\left(\sum_{\ell=1}^3 \frac{\calK_{\ell,i} (\sqrt{n}\cdot x)}{n^{\ell/2}} - \sum_{\ell=1}^3\frac{\calK_{\ell,i} (\sqrt{n+m}\cdot x)}{(n+m)^{\ell/2}}  \right)^2 \mathrm{d}x \\
	&\le 3\sum_{\ell=1}^3 \int_{\bR} \gamma_n(x)\left(\frac{\calK_{\ell,i} (\sqrt{n}\cdot x)}{n^{\ell/2}} - \frac{\calK_{\ell,i} (\sqrt{n+m}\cdot x)}{(n+m)^{\ell/2}}  \right)^2 \mathrm{d}x \\
	&= 3\sum_{\ell=1}^3 \int_{\bR} \gamma(x)\left(\frac{\calK_{\ell,i} (x)}{n^{\ell/2}} - \frac{\calK_{\ell,i} (\sqrt{1+m/n}\cdot x)}{(n+m)^{\ell/2}}  \right)^2 \mathrm{d}x,
	\end{align*}
	where the last step is a change of measure, and $\gamma$ is the density of $\calN(0,1)$. Writing $\calK_{\ell,i}(x) = \sum_{j=0}^{3\ell} a_{i,j}x^j$, then 
	\begin{align*}
	\left| \frac{\calK_{\ell,i} (x)}{n^{\ell/2}} - \frac{\calK_{\ell,i} (\sqrt{1+m/n}\cdot x)}{(n+m)^{\ell/2}} \right| = \left| \sum_{j=0}^{3\ell} \frac{a_{i,j}x^j}{n^{\ell/2}} \left(1 - \left(1+\frac{m}{n}\right)^{\frac{j-\ell}{2}} \right) \right| \lesssim \frac{m(1+x^{3\ell})}{n^{\ell/2+1}}
	\end{align*}
	whenever $m= O(n)$. Combining the above two inequalities yields
	\begin{align}\label{eq:A2}
	A_2 \lesssim \frac{m^2}{n^3}. 
	\end{align}
	\item Upper bounding $A_3$: note that the tail of $\gamma_n(x)$ is at least $\exp(-\Omega(n))$ when $x\ge c$, integrating the tail leads to
	\begin{align}\label{eq:A3}
	A_3 = \exp(-\Omega(n)). 
	\end{align}
\end{enumerate}

In summary, a combination of \eqref{eq:A1}, \eqref{eq:A2}, and \eqref{eq:A3} leads to 
\begin{align*}
H^2(Q_{n,i}^+, Q_{n+m,i}^+) = O\left(\frac{m^2}{n^2} + \exp(-\Omega(n))\right). 
\end{align*}
Now using \eqref{eq:inequality_TV_hellinger}, \eqref{eq:tensorization_hellinger}, and $|Q_{n,i}|(\bR)\le 1+\exp(-\Omega(n))$, we conclude that
\begin{align}\label{eq:Hellinger_gaussian}
\left\|\prod_{i=1}^d Q_{n,i}^+ - \prod_{i=1}^d Q_{n+m,i}^+ \right\|_{\text{TV}} = O\left(\frac{m\sqrt{d}}{n} + (1+\exp(-\Omega(n)))^d \right). 
\end{align}
Hence, when $n=\Omega(\sqrt{d})$ and $m=O(n)$, the desired result follows from \eqref{eq:edgeworth_approx} and \eqref{eq:Hellinger_gaussian}. For the other scenarios, we simply use that the TV distance is upper bounded by one, and the result still holds. 

\subsection{Proof of Theorem \ref{thm:shuffling_general}}
Let $\widehat{P}_n$ be the distribution learned from the first $n/2$ samples which achieves the $\chi^2$-learning error $r_{\chi^2}(\calP,n/2)$, and $P_{\text{mix}}$ be the distribution of the shuffled samples $(Z_1,\cdots,Z_{n/2+m})$ in Algorithm \ref{alg:shuffling_general}. Note that both distributions depend on the first $n/2$ samples and are therefore random. Then the final TV distance of sample amplification is
\begin{align*}
\| P_{X^{n/2}}\times P_{\text{mix}}(X^{n/2}) - P^{\otimes (n+m)} \|_{\text{TV}} = \bE_{X^{n/2}} \| P_{\text{mix}}(X^{n/2}) - P^{\otimes (n/2+m)} \|_{\text{TV}}, 
\end{align*}
which is the expected TV distance between the mixture distribution and the product distribution. 

By Lemma \ref{lemma:shuffle_chi}, for any realization of $X^{n/2}$ it holds that
\begin{align*}
\chi^2\left(P_{\text{mix}}, P^{\otimes (n/2+m)} \right) \le \left(1 + \frac{m}{n/2+m}\chi^2(\widehat{P}_n, P) \right)^m - 1. 
\end{align*}
Since \eqref{eq:TV_KL_chi} shows that $D_{\text{KL}}(P\|Q)\le \log(1+\chi^2(P,Q))$, we have
\begin{align*}
D_{\text{KL}}(P_{\text{mix}} \| P^{\otimes (n/2+m)} ) \le m\log\left(1 + \frac{m}{n/2+m}\chi^2(\widehat{P}_n, P) \right) \le \frac{m^2}{n/2+m}\chi^2(\widehat{P}_n, P). 
\end{align*}
Consequently, by \eqref{eq:TV_KL_chi} again and the concavity of $x\mapsto \sqrt{x}$, we have
\begin{align*}
 \bE_{X^{n/2}} \| P_{\text{mix}}(X^{n/2}) - P^{\otimes (n/2+m)} \|_{\text{TV}} &\le \bE_{X^{n/2}} \sqrt{\frac{1}{2}D_{\text{KL}}(P_{\text{mix}}(X^{n/2}) \| P^{\otimes (n/2+m)} )} \\
 &\le \bE_{X^{n/2}} \sqrt{ \frac{m^2}{n+2m} \chi^2(\widehat{P}_n(X^{n/2}), P) } \\
 &\le \sqrt{\frac{m^2}{n+2m}  \bE_{X^{n/2}}[\chi^2(\widehat{P}_n(X^{n/2}), P)]} \\
 &\le \sqrt{\frac{m^2}{n+2m}\cdot r_{\chi^2}(\calP,n/2)}.
\end{align*}

\subsection{Proof of Theorem \ref{thm:shuffling_product}}
The main arguments are essentially the same as Theorem \ref{thm:shuffling_general}. Note that by the same argument, for each $j\in [d]$ we have
\begin{align*}
D_{\text{KL}}(P_{\text{mix},j} \| P_j^{\otimes (n/2+m)}) \le \frac{m^2}{n/2+m}\chi^2(\widehat{P}_{n,j},P_j). 
\end{align*}
Note that both $P_{\text{mix}}$ and $P$ have product structures, the chain rule of KL divergence implies that
\begin{align*}
D_{\text{KL}}(P_{\text{mix}} \| P^{\otimes (n/2+m)}) \le \frac{m^2}{n/2+m}\sum_{j=1}^d\chi^2(\widehat{P}_{n,j},P_j). 
\end{align*}
Now the rest of the proof follows from the same last few lines of that of Theorem \ref{thm:shuffling_general}. 

\subsection{Proof of Theorem \ref{thm:Gaussian}}
	The proof relies on Lemma \ref{lemma:general_lower_bound} and a classical statistical result known as Anderson's lemma. Without loss of generality we assume $\Sigma = I_d$. Choose $L(\theta,\widehat{\theta}) = \ell(\theta - \widehat{\theta})$, with a bowl-shaped (i.e. symmetric and quasi-convex) loss function $\ell(\cdot)\in [0,1]$. Then Anderson's lemma (see, e.g. \cite[Lemma 8.5]{Vandervaart2000}) implies that the minimax estimator under $L$ and $n$ samples is $\widehat{\theta}_n = n^{-1}\sum_{i=1}^n X_i$, thus
	\begin{align*}
	r(\calP,n,L) = \bE\left[\ell\left( \frac{Z}{\sqrt{n}} \right) \right], 
	\end{align*}
	with $Z\sim \calN(0,I_d)$. Choosing $\ell(x) = \1(\|x\|_2 \ge r)$ with the parameter $r>0$ determined by
	\begin{align*}
	\bP\left( \left| \frac{Z}{\sqrt{n}} \right| \ge r\right) - \bP\left( \left| \frac{Z}{\sqrt{n+m}} \right| \ge r\right) = \left\|\calN\left(0,\frac{I_d}{n}\right) - \calN\left(0,\frac{I_d}{n+m}\right) \right\|_{\text{\rm TV}}, 
	\end{align*}
	then clearly $\ell(\cdot)\in \{0,1\}$ is bowl-shaped. Hence, for this choice of $L$, Lemma \ref{lemma:general_lower_bound} gives that
	\begin{align*}
	\varepsilon^\star(\calP,n,m) \ge r(\calP,n,L) - r(\calP,n+m,L) = \left\|\calN\left(0,\frac{I_d}{n}\right) - \calN\left(0,\frac{I_d}{n+m}\right) \right\|_{\text{\rm TV}},
	\end{align*}
	and a matching upper bound is presented in Example \ref{example:mean}. 

\subsection{Proof of Theorem \ref{thm:lower_exponential}}
Based on the discussions above Theorem \ref{thm:lower_exponential}, we choose an arbitrary open ball $\Theta_0\subseteq \Theta$, and pick any $d$-dimensional ball $B_d(\mu_0;r)$ contained in $\nabla A(\Theta_0)$. Let $\theta_0$ be the center of $\Theta_0$, and after a proper affine transformation we assume that $\nabla A^2(\theta_0) = I_d$. Consider a truncated Gaussian prior $\nu$ with
\begin{align*}
\mu\sim \calN(\mu_0, c_nr_n^2) \mid \mu \in B_d(\mu_0; r_n), 
\end{align*}
where $c_n>0$ and $r_n\in (0,r)$ are parameters to be determined later. Using the diffeomorphism $\nabla A$, there is a prior $\nu_0$ on $\Theta_0$ which induces the above prior on $\nabla A(\Theta_0)$. Moreover, as $\overline{\Theta}_0$, the closure of $\Theta_0$, is a compact set, a weaker Assumption \ref{assump.moment} with the supremum restricted to $\Theta_0$ holds for $k=3$. 

Now we analyze the Bayes risks under the above prior $\nu_0$ and the loss 
\begin{align*}
L(\theta, \widehat{\mu}) = \ell(\nabla A(\theta) - \widehat{\mu}) \in [0,1], 
\end{align*}
where $\ell$ is a bowl-shaped function. By sufficiency, it remains to consider the class of estimators depending only on the sufficient statistic $T_n(X^n) = n^{-1}\sum_{i=1}^n T(X_i)$. Let $\widehat{\mu}(T_n)$ be such an estimator, then under each $\theta\in \Theta_0$, we have
\begin{align*}
\bE_\theta[L(\theta, \widehat{\mu}(T_n))] \ge \bE_\theta[L(\theta, \widehat{\mu}(Z_n))] - \|\calL(T_n) - \calL(Z_n) \|_{\text{TV}}, 
\end{align*}
where $Z_n \mid \theta \sim \calN(\nabla A(\theta), I_d/n)$. To upper bound this TV distance, the proof of Theorem \ref{thm:upper_exp_general} together with Assumption \ref{assump.moment} applied to $\Theta_0$ gives that
\begin{align*}
\|\calL(T_n) - \calN(\nabla A(\theta), \nabla^2 A(\theta)/n) \|_{\text{TV}} \le \frac{C_1}{\sqrt{n}}, 
\end{align*}
where $C_1>0$ only depends on the exponential family. Moreover, 
\begin{align*}
&\|\calN(\nabla A(\theta), I_d/n) - \calN(\nabla A(\theta), \nabla^2 A(\theta)/n) \|_{\text{TV}} \\
&= \| \calN(0,I_d) - \calN(0, \nabla^2 A(\theta)) \|_{\text{TV}} \\
&\stepa{\le} \frac{3}{2}\| \nabla^2 A(\theta) - I_d \|_{\text{F}} \stepb{\le} C_2r_n, 
\end{align*}
where (a) follows from \cite[Theorem 1.1]{devroye2018total}, and (b) makes use of the analytical property of $A(\theta)$ (see, e.g. \cite[Theorem 1.17]{miescke2008statistical}), the assumption that $ \|\nabla A(\theta) - \nabla A(\theta_0) \|_2\le r_n$, and $\nabla^2 A(\theta_0) = I_d$. Again, here the constant $C_2>0$ is independent of $n$. Combining the above inequalities yields that
\begin{align}\label{eq:Tn_to_Zn}
\bE_\theta[L(\theta, \widehat{\mu}(T_n))] \ge \bE_\theta[L(\theta, \widehat{\mu}(Z_n))] - \frac{C_1}{\sqrt{n}} - C_2r_n.  
\end{align}

Next we lower bound the Bayes risk $\bE_{\nu_0}\bE_\theta[L(\theta,\widehat{\mu}(Z_n))]$ when $T_n$ has been replaced by $Z_n$. Let $\nu'$ be the non-truncated Gaussian distribution $\calN(\mu_0, c_nr_n^2)$, and $G_n\sim \calN(0,I_d/n)$, then
\begin{align}\label{eq:mu_to_mu'}
\bE_{\nu_0}\bE_\theta[L(\theta,\widehat{\mu}(Z_n))] &= \bE_{\mu\sim \nu}[\ell(\mu - \widehat{\mu}(G_n + \mu) )] \nonumber\\
&\ge \bE_{\mu\sim \nu'}[\ell(\mu - \widehat{\mu}(G_n + \mu) )] - \nu'(\left\{ \mu\notin B_d(\mu_0;r_n) \right\}) \nonumber\\
&\stepc{\ge} \bE_{\mu\sim \nu'} [\ell(\mu - \widehat{\mu}(G_n + \mu) )] - e^{-C_3/c_n}\nonumber \\
&\stepd{\ge} \bE\left[\ell\left(\sqrt{\frac{nc_nr_n^2}{1+nc_nr_n^2}}\cdot G_n\right) \right] - e^{-C_3/c_n} \nonumber\\
&\stepe{\ge} \bE[\ell(G_n)] - e^{-C_3/c_n} - \frac{C_4}{1+nc_nr_n^2}. 
\end{align}
Here (c) follows from the Gaussian tail probability, (d) makes use of Anderson's lemma and the fact that under $\nu'$, the posterior distribution of $\mu$ given $Z_n$ is Gaussian with covariance $nc_nr_n^2I_d/(1+nc_nr_n^2)$. The final inequality (e) is due to the following upper bound on the TV distance: 
\begin{align*}
\|\calN(0,\Sigma) - \calN(0,c\Sigma) \|_{\text{TV}} \le \frac{3}{2}\|(c-1)I_d\|_{\text{F}} = \frac{3\sqrt{d}|c-1|}{2}. 
\end{align*}
Now combining \eqref{eq:Tn_to_Zn} and \eqref{eq:mu_to_mu'}, we obtain a lower bound on the Bayes risk: 
\begin{align*}
r_{\text{B}}(\calP,n,\nu_0,L) \ge \bE[\ell(G_n)] - \frac{C_1}{\sqrt{n}} - C_2r_n - e^{-C_3/c_n} - \frac{C_4}{1+nc_nr_n^2}. 
\end{align*}

Similarly, by reversing all the above inequalities, an upper bound of the Bayes risk with $n+m$ samples is also available: 
\begin{align*}
r_{\text{B}}(\calP,n+m,\nu_0,L) \le \bE[\ell(G_{n+m})] + \frac{C_1}{\sqrt{n}} + C_2r_n + e^{-C_3/c_n} + \frac{C_4}{1+nc_nr_n^2}. 
\end{align*}
By the proof of Theorem \ref{thm:Gaussian}, a proper choice of $\ell$ satisfies that
\begin{align*}
\bE[\ell(G_n)] - \bE[\ell(G_{n+m})] = \left\|\calN\left(0,\frac{I_d}{n}\right) - \calN\left(0,\frac{I_d}{n+m}\right) \right\|_{\text{TV}} = \Omega\left(\frac{m\sqrt{d}}{n}\wedge 1 \right), 
\end{align*}
where the last step is again due to \cite[Theorem 1.1]{devroye2018total}. Consequently, by choosing $c_n = \Theta(1/\log n)$ and $r_n = \Theta((\log n/n)^{1/3})$, the desired result follows from Lemma \ref{lemma:general_lower_bound}. 

\subsection{Proof of Theorem \ref{thm:lower_bound}}
We will apply Lemma \ref{lemma:general_lower_bound} to the uniform prior $\mu$ over $2^d$ points $\prod_{j=1}^d \{\theta_{j,+}, \theta_{j,-} \}$, and the loss function $L: \Theta\times \Theta \to [0,1]$ with 
\begin{align*}
L( (\theta_1,\cdots,\theta_d), (\widehat{\theta}_1,\cdots,\widehat{\theta}_d) ) = \1\left(\sum_{j=1}^d \1(\theta_j = \widehat{\theta}_j) \le \frac{d}{2} + \frac{1}{2}\sum_{j=1}^d \alpha_j \right). 
\end{align*}
We compute the Bayes risks $r_{\text{B}}(\calP,n,\mu,L)$ and $r_{\text{B}}(\calP,n+m,\mu,L)$ in this scenario. 

Under the uniform prior $\mu$, it is straightforward to see that the posterior distribution of $\theta$ given $X^n$ is a product distribution $\prod_{j=1}^d p_{\theta_j\mid X_j^n}$, with $p_{\theta_j \mid X_j^n}$ supported on two elements $\{\theta_{j,+}, \theta_{j,-} \}$, and
\begin{align*}
p_{\theta_j \mid X_j^n} (\theta_{j,+}) &= \frac{ \prod_{i=1}^n p_{\theta_{j,+}}(X_{i,j}) }{\prod_{i=1}^n p_{\theta_{j,+}}(X_{i,j}) + \prod_{i=1}^n p_{\theta_{j,-}}(X_{i,j})}, \\
p_{\theta_j \mid X_j^n} (\theta_{j,-}) &= \frac{ \prod_{i=1}^n p_{\theta_{j,-}}(X_{i,j}) }{\prod_{i=1}^n p_{\theta_{j,+}}(X_{i,j}) + \prod_{i=1}^n p_{\theta_{j,-}}(X_{i,j})}. 
\end{align*}
Then given $X^n$, the Bayes estimator $\widehat{\theta}(X^n)\in \Theta$ which minimizes the expected loss $L$ is a minimizer $\xi\in \Theta$ of the above expression
\begin{align*}
 \bP_{\theta\mid X^n}\left(\sum_{j=1}^d \1(\theta_j = \xi_j) \le  \frac{d}{2} + \frac{1}{2}\sum_{j=1}^d \alpha_j  \right) ,
\end{align*}
which is easily seen to be
\begin{align*}
\widehat{\theta}_j(X^n) = \widehat{\theta}_j(X_j^n) = \theta_{j,-} + (\theta_{j,+} - \theta_{j,-})\cdot \1\left(  \prod_{i=1}^n p_{\theta_{j,+}}(X_{i,j}) \ge  \prod_{i=1}^n p_{\theta_{j,-}}(X_{i,j}) \right). 
\end{align*}

For the above Bayes estimator, the random variables $\1(\theta_j = \widehat{\theta}_j(X^n))$ are mutually independent, with the mean value
\begin{align*}
p_{n,j} &= \bP(\theta_j = \widehat{\theta}_j(X^n)) \\
&= \frac{1}{2}p_{\theta_{j,+}}^{\otimes n} \left(  \prod_{i=1}^n p_{\theta_{j,+}}(X_{i,j}) < \prod_{i=1}^n p_{\theta_{j,-}}(X_{i,j}) \right) +  \frac{1}{2}p_{\theta_{j,-}}^{\otimes n} \left(  \prod_{i=1}^n p_{\theta_{j,+}}(X_{i,j}) \ge \prod_{i=1}^n p_{\theta_{j,-}}(X_{i,j}) \right) \\
&= \frac{1}{2}\left(1 + \|p_{\theta_{j,+}}^{\otimes n}  - p_{\theta_{j,-}}^{\otimes n} \|_{\text{TV}} \right). 
\end{align*}
Consequently, we have $p_{n,j}\le (1+\alpha_j)/2 - \varepsilon/(2\sqrt{d})$ for each $j\in [d]$ by \eqref{eq:small_TV}, and thus
\begin{align}\label{eq:bayes_lower}
r_{\text{B}}(\calP,n,\mu,L) &= \bP\left(\sum_{j=1}^d \mathsf{Bern}(p_{n,j}) \le \frac{d}{2} + \frac{1}{2}\sum_{j=1}^d \alpha_j \right) \nonumber \\
&\ge \bP\left(\sum_{j=1}^d \mathsf{Bern}\left( \frac{1+\alpha_j}{2} - \frac{\varepsilon}{2\sqrt{d}} \right) \le \frac{d}{2} + \frac{1}{2}\sum_{j=1}^d \alpha_j \right). 
\end{align}

An entirely symmetric argument leads to
\begin{align}\label{eq:bayes_upper}
r_{\text{B}}(\calP,n+m,\mu,L) \le  \bP\left(\sum_{j=1}^d \mathsf{Bern}\left( \frac{1+\alpha_j}{2} + \frac{\varepsilon}{2\sqrt{d}} \right) \le \frac{d}{2} + \frac{1}{2}\sum_{j=1}^d \alpha_j \right). 
\end{align}
To further lower bound \eqref{eq:bayes_lower} and upper bound \eqref{eq:bayes_upper}, the following lemma shows that we may assume that the above Bernoulli distributions have the same parameter. 

\begin{lemma}[Theorems 4 and 5 of \cite{hoeffding1956distribution}]\label{lemma:bernoulli_convolution}
Let $X_i\sim \mathsf{Bern}(p_i)$ be independent Bernoulli random variables for $i=1,\cdots,n$, with $\bar{p} = n^{-1}\sum_{i=1}^n p_i$. Then for $0\le k\le n\bar{p}-1$, we have
\begin{align*}
\bP\left(\sum_{i=1}^n X_i \le k \right) \le \bP\left(\mathsf{B}(n,\bar{p}) \le k\right), 
\end{align*}
and for $k\ge n\bar{p}$, we have
\begin{align*}
\bP\left(\sum_{i=1}^n X_i \le k \right) \ge \bP\left(\mathsf{B}(n,\bar{p}) \le k\right). 
\end{align*}
In particular, for integers $(b,c)$ with $b\le n\bar{p}\le c$, we have
\begin{align*}
\bP\left(b\le \sum_{i=1}^n X_i \le c \right) \ge \bP\left(b\le \mathsf{B}(n,\bar{p}) \le c\right). 
\end{align*}
\end{lemma}

For $d \ge 4/\varepsilon^2$, based on the second statement of Lemma \ref{lemma:bernoulli_convolution}, the quantity in \eqref{eq:bayes_lower} satisfies
\begin{align}\label{eq:bayes_lower_binomial}
r_{\text{B}}(\calP,n,\mu,L) \ge \bP\left(\mathsf{B}\left(d, \frac{1+\alpha}{2}  - \frac{\varepsilon}{2\sqrt{d}}\right) \le \frac{1+\alpha}{2}\cdot d \right), 
\end{align}
with $\alpha\triangleq d^{-1}\sum_{j=1}^d \alpha_j$. Similarly, based on the first statement of Lemma \ref{lemma:bernoulli_convolution}, for \eqref{eq:bayes_upper} we have
\begin{align}\label{eq:bayes_upper_binomial}
r_{\text{B}}(\calP,n+m,\mu,L) \le \bP\left(\mathsf{B}\left(d, \frac{1+\alpha}{2} + \frac{\varepsilon}{2\sqrt{d}} \right) \le \frac{1+\alpha}{2}\cdot d \right). 
\end{align}
Since
\begin{align*}
\frac{\mathrm{d} \bP(\mathsf{B}(n,x)=k)}{\mathrm{d}x} = n\left(\bP(\mathsf{B}(n-1,x)=k-1) - \bP(\mathsf{B}(n-1,x)=k)\right), 
\end{align*}
we invoke Lemma \ref{lemma:general_lower_bound} and lower bound the Bayes risk difference as
\begin{align*}
\varepsilon^\star(\calP,n,m) &\ge r_{\text{B}}(\calP,n,\mu,L) - r_{\text{B}}(\calP,n+m,\mu,L) \\
&\ge \bP\left(\mathsf{B}\left(d, \frac{1+\alpha}{2}  - \frac{\varepsilon}{2\sqrt{d}}\right) \le \frac{1+\alpha}{2}\cdot d \right)  - \bP\left(\mathsf{B}\left(d, \frac{1+\alpha}{2} + \frac{\varepsilon}{2\sqrt{d}}\right) \le \frac{1+\alpha}{2}\cdot d \right) \\
&= -\sum_{0\le k\le (1+\alpha)d/2} \int_{\frac{1+\alpha}{2}  - \frac{\varepsilon}{2\sqrt{d}}}^{\frac{1+\alpha}{2} + \frac{\varepsilon}{2\sqrt{d}}} \frac{\mathrm{d} \bP(\mathsf{B}(d,x)=k)}{\mathrm{d}x} \mathrm{d}x \\
&= d \int_{\frac{1+\alpha}{2}  - \frac{\varepsilon}{2\sqrt{d}}}^{\frac{1+\alpha}{2} + \frac{\varepsilon}{2\sqrt{d}}} \bP\left(\mathsf{B}(d-1, x) = \left\lfloor \frac{1+\alpha}{2}\cdot d \right\rfloor \right)\mathrm{d}x \\
&\stepa{\ge} d \int_{\frac{1+\alpha}{2}  - \frac{\varepsilon}{2\sqrt{d}}}^{\frac{1+\alpha}{2} + \frac{\varepsilon}{2\sqrt{d}}} \frac{c(\underline{\alpha}, \overline{\alpha})}{\sqrt{d}}\mathrm{d}x = c(\underline{\alpha}, \overline{\alpha}) \varepsilon, 
\end{align*}
where (a) is due to 
\begin{align}\label{eq:stirling}
\min_{p_0\le p\le p_1, |k-np|\le C\sqrt{n}}\bP(\mathsf{B}(n,p) = k) = \Omega_{p_0, p_1,C}\left(\frac{1}{\sqrt{n}}\right)
\end{align}
for any $p_0, p_1\in (0,1)$ and $C>0$, by Stirling's approximation. 

For $d < 4/\varepsilon^2$, we first note the following identity:
\begin{align*}
&\frac{\mathrm{d}}{\mathrm{d}x}\bigg|_{x=0} \bP\left(\sum_{i=1}^n \mathsf{Bern}(p_i + x) = k\right) \\
&= \frac{\mathrm{d}}{\mathrm{d}x}\bigg|_{x=0} \sum_{w\in \{0,1\}^n: \sum_{i=1}^n w_i = k} \prod_{i=1}^n (p_i+x)^{w_i}(1-p_i-x)^{1-w_i}  \\
&= \sum_{i=1}^n \left(\sum_{w_{\backslash i}\in \{0,1\}^{n-1}: \sum_{j\neq i} w_j = k-1} \prod_{j\neq i} p_j^{w_j}(1-p_j)^{1-w_j} - \sum_{w_{\backslash i}\in \{0,1\}^{n-1}: \sum_{j\neq i} w_j = k} \prod_{j\neq i} p_j^{w_j}(1-p_j)^{1-w_j}\right) \\
&=\sum_{i=1}^n \left[\bP\left(\sum_{j\neq i} \mathsf{Bern}(p_j) = k-1\right) - \bP\left(\sum_{j\neq i} \mathsf{Bern}(p_j) = k\right)\right]. 
\end{align*}
Based on \eqref{eq:bayes_lower} and \eqref{eq:bayes_upper}, we then have
\begin{align*}
&\varepsilon^\star(\calP,n,m) \ge r_{\text{B}}(\calP,n,\mu,L) - r_{\text{B}}(\calP,n+m,\mu,L) \\
&\ge  \bP\left(\sum_{j=1}^d \mathsf{Bern}\left( \frac{1+\alpha_j}{2} - \frac{\varepsilon}{2\sqrt{d}} \right) \le \frac{(1+\alpha)d}{2}\right) -  \bP\left(\sum_{j=1}^d \mathsf{Bern}\left( \frac{1+\alpha_j}{2} + \frac{\varepsilon}{2\sqrt{d}} \right) \le \frac{(1+\alpha)d}{2} \right) \\
&= -\sum_{0\le k\le (1+\alpha)d/2} \int_{-\frac{\varepsilon}{2\sqrt{d}}}^{\frac{\varepsilon}{2\sqrt{d}}} \frac{\mathrm{d}}{\mathrm{d}x} \bP\left(\sum_{j=1}^d \mathsf{Bern}\left( \frac{1+\alpha_j}{2} + x \right) = k\right) \mathrm{d}x \\
&= \sum_{i=1}^d \int_{-\frac{\varepsilon}{2\sqrt{d}}}^{\frac{\varepsilon}{2\sqrt{d}}} \bP\left(\sum_{j\neq i} \mathsf{Bern}\left( \frac{1+\alpha_j}{2} + x \right) = \left\lfloor \frac{(1+\alpha)d}{2}\right\rfloor \right) \mathrm{d}x. 
\end{align*}
We will show that the integrand is uniformly of the order $\Omega(1/\sqrt{d})$. To this end, note that for $|x|\le \varepsilon/(2\sqrt{d}) < 1/d$ as $d < 4/\varepsilon^2$, it holds that
\begin{align*}
\sum_{j\neq i} \left(\frac{1+\alpha_j}{2} + x\right) &< \sum_{j=1}^d \frac{1+\alpha_j}{2} + 1 \le \left\lfloor \frac{(1+\alpha)d}{2}\right\rfloor + 2, \\
\sum_{j\neq i} \left(\frac{1+\alpha_j}{2} + x\right) &> \sum_{j=1}^d \frac{1+\alpha_j}{2} - 1 - 1 \ge \left\lfloor \frac{(1+\alpha)d}{2}\right\rfloor - 2. 
\end{align*}
Consequently, by the last statement of Lemma \ref{lemma:bernoulli_convolution}, one has
\begin{align*}
&\bP\left(\left| \sum_{j\neq i} \mathsf{Bern}\left( \frac{1+\alpha_j}{2} + x \right) - \left\lfloor \frac{(1+\alpha)d}{2}\right\rfloor \right| \le 2 \right) \\ 
&\ge \bP\left(\left| \mathsf{B}\left(d-1, \frac{1}{d-1}\sum_{j\neq i} \left(\frac{1+\alpha_j}{2} + x\right) \right) - \left\lfloor \frac{(1+\alpha)d}{2}\right\rfloor \right| \le 2 \right) = \Omega_{\underline{\alpha},\overline{\alpha}}\left(\frac{1}{\sqrt{d}}\right), 
\end{align*}
where the last step is again due to \eqref{eq:stirling}. The above display is a lower bound for the probability of a size-$5$ set, and in view of the following lemma, the same $\Omega_{\underline{\alpha},\overline{\alpha}}(1/\sqrt{d})$ lower bound also holds for the probability of any singleton. Plugging this lower bound back into the integral then yields to $\varepsilon^\star(\calP,n,m) = \Omega_{\underline{\alpha},\overline{\alpha}}(\varepsilon)$, as desired. 

\begin{lemma}\label{lemma:bernoulli_conv_ratio}
Let $p_1,\cdots,p_n\in [a,b]$ with $0<a\le b<1$, and $cn\le k<k+1\le (1-c)n$ for some $c>0$. Define
\begin{align*}
    f(k) = \bP\left(\sum_{i=1}^n \mathsf{Bern}(p_i) = k\right). 
\end{align*}
Then there exists an absolute constant $C = C(a,b,c) < \infty$ such that
\begin{align*}
    C^{-1} \le \frac{f(k+1)}{f(k)} \le C. 
\end{align*}
\end{lemma}
\begin{proof}
Let $W_k = \{w\in \{0,1\}^n: \sum_{i=1}^n w_i = k\}$, and $g(w) = \prod_{i=1}^n p_i^{w_i}(1-p_i)^{1-w_i}$. Then it is clear that
\begin{align*}
    f(k) = \sum_{w\in W_k} g(w). 
\end{align*}
Call two binary vectors $w$ and $w'$ as \emph{neighbors} (denoted by $w\sim w'$) if $w$ and $w'$ only differ in one coordinate. It is clear that every $w\in W_k$ has $n-k$ neighbors in $W_{k+1}$, and every $w\in W_{k+1}$ has $k+1$ neighbors in $W_k$. Moreover, for $w\sim w'$, 
\begin{align*}
    \frac{g(w)}{g(w')} \le \max\left\{\frac{b}{a}, \frac{1-a}{1-b}\right\} =: \rho. 
\end{align*}
Consequently, by double counting, 
\begin{align*}
f(k+1) &= \sum_{w\in W_{k+1}} g(w) = \frac{1}{k+1}\sum_{w'\in W_k}\sum_{w\in W_{k+1}: w\sim w'} g(w) \\
&\le \frac{\rho}{k+1}\sum_{w'\in W_k}\sum_{w\in W_{k+1}: w\sim w'} g(w') = \frac{(n-k)\rho}{k+1} \sum_{w'\in W_k} g(w') \\
&= \frac{(n-k)\rho}{k+1} f(k) \le \frac{(1-c)\rho}{c} f(k). 
\end{align*}
The other inequality can be established analogously. 
\end{proof}

\subsection{Proof of Theorem \ref{thm:lower_product}}

For each $j\in [d]$, we pick two points $\theta_{j,+}, \theta_{j,-}\in \Theta_j$ in Assumption \ref{assump.hellinger}. We first aim to show that
\begin{align}
0.09 \triangleq \varepsilon_1' \le \| p_{\theta_{j,+}}^{\otimes n} - p_{\theta_{j,-}}^{\otimes n} \|_{\text{TV}} &\le \varepsilon_1 \triangleq 0.6, \label{eq:TV_n_sample} \\
0.99995 \triangleq \varepsilon_2' \ge \| p_{\theta_{j,+}}^{\otimes 20n} - p_{\theta_{j,-}}^{\otimes 20n} \|_{\text{TV}} &\ge \varepsilon_2 \triangleq 0.86. \label{eq:TV_20n_sample}
\end{align}

The proofs of \eqref{eq:TV_n_sample} and \eqref{eq:TV_20n_sample} rely on the tensorization property of the Hellinger distance
\begin{align*}
H^2\left(\prod_{i=1}^d P_i, \prod_{i=1}^d Q_i \right) = 1 - \prod_{i=1}^d\left(1 - H^2(P_i, Q_i) \right),
\end{align*}
and the relationship between TV and Hellinger distance in \eqref{eq:TV_Hellinger}. For example, for \eqref{eq:TV_n_sample}, we have
\begin{align*}
\| p_{\theta_{j,+}}^{\otimes n} - p_{\theta_{j,-}}^{\otimes n} \|_{\text{TV}} &\le \sqrt{1 - (1-H^2(p_{\theta_{j,+}}^{\otimes n},p_{\theta_{j,-}}^{\otimes n}))^2} \\
&= \sqrt{1 - \left(1-H^2(p_{\theta_{j,+}}, p_{\theta_{j,-}})\right)^{2n}} \\
&\le \sqrt{1 - \left(1-\frac{1}{5n}\right)^{2n}} \\
&\le \sqrt{\frac{9}{25}} = 0.6. 
\end{align*}
Applying the other inequality, for \eqref{eq:TV_20n_sample} we have
\begin{align*}
\| p_{\theta_{j,+}}^{\otimes 20n} - p_{\theta_{j,-}}^{\otimes 20n} \|_{\text{TV}} &\ge H^2(p_{\theta_{j,+}}^{\otimes 20n} , p_{\theta_{j,-}}^{\otimes 20n} ) \\
&\ge 1 - \left(1 - \frac{1}{10n}\right)^{20n} \\
&\ge 1 - \exp(-2) > 0.86. 
\end{align*}
The other inequalities involving $\varepsilon_1'$ and $\varepsilon_2'$ could be established analogously. 

Next, for the choice of $m = \lceil c\varepsilon n/\sqrt{d} \rceil$ in the statement of Theorem \ref{thm:lower_product}, we show that there exists $n_j\in [n,20n]$ such that 
\begin{align}\label{eq:TV_diff}
\| p_{\theta_{j,+}}^{\otimes (n_j+m)} - p_{\theta_{j,-}}^{\otimes (n_j+m)} \|_{\text{TV}} - \| p_{\theta_{j,+}}^{\otimes n_j} - p_{\theta_{j,-}}^{\otimes n_j} \|_{\text{TV}} \ge \frac{\varepsilon_2 - \varepsilon_1}{\lceil 19\sqrt{d}/(c\varepsilon) \rceil}. 
\end{align}
To prove \eqref{eq:TV_diff}, first note that $t\mapsto f_j(t)\triangleq \| p_{\theta_{j,+}}^{\otimes t} - p_{\theta_{j,-}}^{\otimes t} \|_{\text{TV}}$ is non-decreasing by the data-processing property of the TV distance. Moreover, by \eqref{eq:TV_n_sample} and \eqref{eq:TV_20n_sample}, we have $f_j(n)\le \varepsilon_1$ and $f_j(20n)\ge \varepsilon_2$. Consequently, 
\begin{align*}
\varepsilon_2 - \varepsilon_1 &\le f_j(20n) - f_j(n) \\
&\le \sum_{k=1}^{\lceil 19n/m\rceil - 1} [f_j(n+km) - f_j(n+(k-1)m)] + [f_j(20n) - f_j(20n-m)] \\
&\le \left\lceil \frac{19n}{m} \right\rceil\cdot \max_{n\le n_j\le 20n-m} [f_j(n_j+m) - f_j(n_j)], 
\end{align*}
which gives \eqref{eq:TV_diff}. In addition, we also have $\varepsilon_1'\le f_j(n_j)\le f_j(n_j+m) \le \varepsilon_2'$. 

Next we are about to apply Theorem \ref{thm:lower_bound}. By \eqref{eq:TV_diff}, there exist $\alpha_j\in (\varepsilon_1', \varepsilon_2')$ such that
\begin{align*}
 \| p_{\theta_{j,+}}^{\otimes n_j} - p_{\theta_{j,-}}^{\otimes n_j} \|_{\text{TV}} &\le \alpha_j - \Omega_c\left(\frac{\varepsilon}{\sqrt{d}}\right), \\
  \| p_{\theta_{j,+}}^{\otimes (n_j+m)} - p_{\theta_{j,-}}^{\otimes (n_j+m)} \|_{\text{TV}} &\ge \alpha_j + \Omega_c\left(\frac{\varepsilon}{\sqrt{d}}\right). 
\end{align*}
Therefore, Theorem \ref{thm:lower_bound} shows that there is an absolute constant $c'>0$ depending only on $(c,\varepsilon_1',\varepsilon_2')$ such that 
\begin{align*}
\varepsilon^\star(\calP, (n_1,\cdots,n_d), m) \ge c'\varepsilon, 
\end{align*}
where the above quantity denotes the minimax error in a new sample amplification problem: suppose we draw $n_j$ independent samples from $\calP_j$, also independently for each $j\in [d]$, and we aim to amplify into $n_j+m$ independent samples from $\calP_j$. In other words, the sample sizes for different dimensions may not be equal in the new problem, but the target is still to generate $m$ more samples. We claim that $\varepsilon^\star(\calP, (n_1,\cdots,n_d), m) \le \varepsilon^\star(\calP,n,m)$, and thereby complete the proof. To show the claim, note that $n_j\ge n$ for all $j\in [d]$, hence in the new problem we could keep $n_j-n$ samples unused for each $j$, use the remaining samples to amplify into $n+m$ vectors, and add the above unused samples back to form the final amplification.

\subsection{Proof of Theorem \ref{thm:distribution_top_element}}
We first prove the upper bound. Consider the distribution estimator $\widehat{P}_n = (1,0,\cdots,0)$, which has a $\chi^2$-divergence
\begin{align*}
\chi^2(\widehat{P}_n \| P) = \frac{1}{t} - 1, \qquad \forall P\in \calP_{d,t}. 
\end{align*}
Consequently, the $\chi^2$-estimation error $r_{\chi^2}(\calP_{d,t}, n)$ is at most $1/t$, and Theorem \ref{thm:shuffling_general} states that the random shuffling approach achieves an $(n,n+1,0.1)$ sample amplification if $n=\Theta(1/t)$. 

Next we prove the lower bound. Let $n = 1/(100t)$ and $d$ be a multiple of $100$. Consider the following prior $\mu$ over $\calP_{d,t}$: $\mu$ is the uniform prior over $(p_0,\cdots,p_d)\in \calP_{d,t}$ with $p_0=t$ and the remaining $1-t$ mass evenly distributed on a uniformly random subset of $[d]$ of size $d/100$. The action space $\calA$ is chosen to be $\calX^{n+1}$, and the loss function is
\begin{align*}
    L(P, x^{n+1}) &= 1 - \1(x^{n+1}\text{ belongs to the support of }P, \\
    &\qquad \text{ does not contain symbol }0\text{ or repeated symbols}). 
\end{align*}
We first show that $r_{\text{B}}(\calP_{d,t}, n+1, \mu) \le 0.1$. In fact, after observing $n+1$ samples $X^{n+1}$, we simply use $X^{n+1}$ as the estimator under the above loss. Clearly $X^{n+1}$ belongs to the support of $P$. For the remaining conditions, 
\begin{align*}
\bP(X^{n+1}\text{ contains symbol }0) &= 1- (1-t)^{n+1} = 1 - (1-t)^{1/(100t)+1} \le 0.05, \\
\bP(X^{n+1}\text{ contains repeated symbols}) &\le \binom{n+1}{2}\left(t^2 + \sum_{j=1}^{d/100} \frac{1}{(d/100)^2}\right) \\
&\le n^2\left(t^2 + \frac{100}{d}\right) \le \frac{1}{10000} + \frac{1}{100} < 0.05.
\end{align*}
By the union bound, this estimator achieves a Bayes risk at most $0.1$, and thus $r_{\text{B}}(\calP_{d,t}, n+1, \mu) \le 0.1$. 

Next we show that $r_{\text{B}}(\calP_{d,t},n,\mu) \ge 0.9$, which combined with Lemma \ref{lemma:general_lower_bound} gives the desired lower bound. To show this, consider the new symbol in the estimator $x^{n+1}$ not in $X^n$ when the learner observes $X^n$. Since $x^{n+1}$ has length $n+1>n$, there is at least one such symbol. In order to have $L(P,x^{n+1})=0$, this new symbol cannot be $0$ or appear in $X^n$. Moreover, the posterior distribution of the support of $P\sim \mu$ given $X^n$ is uniformly distributed over
\begin{align*}
    \{0, X_1, \cdots, X_n\} \cup \left\{S \subseteq [d] \backslash \{X_1, \cdots, X_n\}: |S| = \frac{d}{100} - n \right\}. 
\end{align*}
Therefore, the posterior probability of the new symbol being outside the support of $P$ (recall that it could be neither one of $\{X_1,\cdots,X_n\}$ nor $0$) is at least
\begin{align*}
1 - \frac{d/100 - n}{d - n} = \frac{0.99d}{d-n} \ge 0.99 > 0.9,
\end{align*}
giving the desired inequality $r_{\text{B}}(\calP_{d,t},n,\mu) \ge 0.9$. 

\subsection{Proof of Theorem \ref{thm:lowrank_covariance}}
The upper bound of sample amplification directly follows from that of learning, and it remains to show the lower bound $n\ge d$. If $n\le d-1$, with probability one the observations $X^n$ spans an $n$-dimensional subspace of the row space of $\Sigma$. An $(n,n+1,0.1)$ sample amplification calls for at least one additional observation not in $X^n$, which with probability $1$ should not belong to the $n$-dimensional subspace spanned by $X^n$ for $n\le d-1$. However, since $p\ge d+1$, under a uniformly chosen $d$-dimensional row space of $\Sigma$, the posterior probability of the additional observation belonging to the row space of $\Sigma$ is zero. Consequently, an $(n,n+1,0.1)$ sample amplification is impossible if $n<d$. 

\subsection{Proof of Theorem \ref{thm:lipschitz}}
We prove the three claims separately. 

\textbf{The proof of $m^\star(\calP_c, n)\gtrsim_c n^{5/6}$.} Classical theory of nonparametric estimation (see, e.g. \cite{tsybakov2009introduction}) tells that there exists a density estimator $\widehat{f}$ such that $\bE_f\|\widehat{f}-f\|_2^2 \lesssim n^{-2/3}$. Since $f$ is lower bounded by $c$, this implies that
\begin{align*}
    \chi^2(\widehat{f}, f) \le \frac{\|\widehat{f}-f\|_2^2}{c} \lesssim_c n^{-2/3}. 
\end{align*}
Consequently, $r_{\chi^2}(\calP_c,n) \lesssim n^{-2/3}$, and Theorem \ref{thm:shuffling_general} implies that $m^\star(\calP_c, n) \gtrsim_c n^{5/6}$. 

\textbf{The proof of $m^\star(\calP_c, n)\lesssim_c n^{5/6}$.} We construct a parametric subfamily of $\calP_c$ and invoke Theorem \ref{thm:lower_product}. Let $g$ be a $1$-Lipschitz function supported on $[0,1]$ with $\int_0^1 g(x)dx = 0$ and $\|g\|_2>0$; in particular $\|g\|_\infty \le 1$. Let $h=n^{-1/3}$ and assume that $M:=h^{-1}$ is an integer. For $u=(u_1,\cdots,u_M)\in \{\pm 1\}^M$, define
\begin{align*}
    f_u(x) = 1 + c_0\sum_{i=1}^M u_ihg\left(\frac{x-(i-1)h}{h}\right), 
\end{align*}
where $c_0\in (0,1)$ is a small constant satisfying $c_0 h\le 1-c$. Consequently, $f_u\in \calP_c$ for every $u\in \{\pm 1\}^M$. However, the density estimation model $X_1,\cdots,X_n\sim f_u$ is not a product model, so Theorem \ref{thm:lower_product} cannot be directly applied. 

To overcome the above issue, we consider a \emph{Poissonized} model as follows: first we draw $N\sim \text{Poi}(n)$, and then draw $N$ i.i.d. samples $X_1,\cdots,X_N\sim f_u$. The Poissonized model satisfies the following two properties: 
\begin{enumerate}
    \item For any measurable set $A\subseteq [0,1]$, we have
    \begin{align*}
        M(A) := |\{i \in [N]: X_i \in A \}| \sim \text{Poi}\left(n\int_A f(x) \mathrm{d}x\right).
    \end{align*}
    \item For any collection of disjoint subsets $\{A_i, i\ge 1\}$, the random variables $\{M(A_i), i\ge 1\}$ are mutually independent. 
\end{enumerate}
For $i\in [M]$, let $A_i=[(i-1)/M, i/M)$, so that $\int_{A_i} f_u(x) \mathrm{d}u = 1/M$ for every $u\in \{\pm 1\}^M$. Clearly, there is a one-to-one correspondence between $(X_1,\cdots,X_N)$ and $(Y_1,\cdots,Y_M)$, where $Y_i$ is the collection of observations in $X_1,\cdots,X_N$ that falls into the set $A_i$. By the above two properties, $(Y_1,\cdots,Y_M)$ are mutually independent, and $Y_i\sim f_{i,u_i}^{\otimes n}$ under $f_u$. Here $f_{i,u_i}$ is the probability distribution of the following process: sample $N_i\sim \text{Poi}(1/M)$, and draw $N_i$ i.i.d. samples from the density
\begin{align*}
    M\left(1 + c_0u_i hg\left(\frac{x-(i-1)h}{h}\right) \right)
\end{align*}
supported on $A_i$. In addition, 
\begin{align*}
    H^2(f_{i,+1}, f_{i,-1}) \asymp \bE[N_i]\cdot Mc_0^2 h^3\|g\|_2^2 \asymp \frac{1}{n}, 
\end{align*}
so the Poissonized model is a product model which satisfies the prerequisite of Theorem \ref{thm:lower_product}. 

Next we denote the i.i.d. sampling model by $\calP_c^{\otimes n}$, and the Poissonized model by $\calP_{c}^{\text{Poi}(n)}$. Let $n_1 = n+C\sqrt{n}, n_2 = n+m-C\sqrt{n+m}$, where $C>0$ is a large universal constant to be chosen later. By Theorem \ref{thm:lower_product} applied to $d=M\asymp n^{1/3}$, Le Cam's distance in Definition \ref{def:le_cam_distance} between Poissonized models satisfies
\begin{align*}
\Delta(\calP_{c}^{\text{Poi}(n_1)}, \calP_{c}^{\text{Poi}(n_2)}) = \Omega(1), 
\end{align*}
as long as $m=\Omega(n^{5/6})$. In addition, the model $\calP_{c}^{\text{Poi}(n_1)}$ is more informative than $\calP_{c}^{\otimes n}$ with probability at least $1-\varepsilon_1$, where $\varepsilon_1 = \bP(\text{Poi}(n_1) < n)$. Similarly, $\calP_{c}^{\otimes (n+m)}$ is more informative than $\calP_{c}^{\text{Poi}(n_2)}$ with probability at least $1-\varepsilon_2$, where $\varepsilon_2 = \bP(\text{Poi}(n_2) >n+m)$. Therefore, 
\begin{align*}
\Delta(\calP_c^{\otimes n}, \calP_c^{\otimes (m+n)}) \ge \Delta(\calP_{c}^{\text{Poi}(n_1)}, \calP_{c}^{\text{Poi}(n_2)}) - \varepsilon_1 - \varepsilon_2 = \Omega(1),
\end{align*}
where in the last step we have $\varepsilon_1, \varepsilon_2\to 0$ by choosing $C>0$ large enough. 

\textbf{The proof of $m^\star(\calP, n)\lesssim n^{3/4}$.} Suppose $n\ge m\ge Cn^{3/4}$ for a large constant $C>0$. Consider the following prior $\mu$ on the density $f$: let $M = \sqrt{n}$, and $u=(u_1,\cdots,u_M)$ be uniformly distributed over all vectors in $\{0,1\}^M$ such that the number of $1$'s is $M/100$. Given $u$, we construct
\begin{align*}
    f_u(x) = \begin{cases}  u_i(2/M - 4|x-(2i-1)/2M|) & \text{if } x \in A_i \triangleq [(i-1)/M, i/M), i\in [M], \\
     (8-2/(25M))(x - 1/2) & \text{if } 1/2\le x\le 1, 
    \end{cases}
\end{align*}
and let $f=f_u$. It is not hard to verify that each $f_u$ is a density and $8$-Lipschitz, so $f_u\in \calP$. 

Next under the context of Lemma \ref{lemma:general_lower_bound}, we choose the action space to be $[0,1]^{n+m}$, as well as the following loss: 
\begin{align*}
    L(f,x^{n+m}) &= 1 - \1(x^{n+m}\text{ belongs to the support of }f, \\
    &\qquad \text{ and fall into at least } [1-(1-n^{-1})^n]\sqrt{n}/100 + C_0n^{1/4} \text{ sets in } \{A_i\}_{i=1}^M ), 
\end{align*}
where $C_0>0$ is a large enough constant. We aim to show that under the loss $L$ and prior $\mu$, the Bayes risk satisfies $r_{\text{B}}(\calP^{\otimes (n+m)}, \mu) \le 0.1$ and $r_{\text{B}}(\calP^{\otimes n}, \mu) \ge 0.9$. Then an application of Lemma \ref{lemma:general_lower_bound} completes the proof of $m^\star(\calP,n)\lesssim n^{3/4}$. 

We first show that $r_{\text{B}}(\calP^{\otimes (n+m)}, \mu) \le 0.1$. We simply use the observed sample $X^{n+m}$ as the estimator $x^{n+m}$ under the above loss, then clearly $x^{n+m}$ belongs to the support of $f$. For the second event in the loss function, let $U_{n+m}$ be the number of sets in $\{A_i\}_{i=1}^M$ which receive any observations in $X^{n+m}$. By the linearity of expectation, as well as the negative dependence across different set counts, we compute for every $f_u$ that
\begin{align*}
\bE[U_{n+m}] &= \sum_{i=1}^{M/100} \left[1-\left(1-\frac{1}{M^2}\right)^{n+m}\right] = \frac{\sqrt{n}}{100}\left[1 - \left(1-\frac{1}{n}\right)^n\right] + \Omega\left(\frac{m}{\sqrt{n}}\right), \\
\var(U_{n+m}) &\le \bE[U_{n+m}] = O(\sqrt{n}). 
\end{align*}
By Chebyshev's inequality, for $m\ge Cn^{3/4}$ with a large enough $C>0$, we have $r_{\text{B}}(\calP^{\otimes (n+m)}, \mu) \le 0.1$. 

Next we show that $r_{\text{B}}(\calP^{\otimes n}, \mu) \ge 0.9$. Let $(i_1,\cdots,i_k)\in [M]^k$ be the indices such that the set $A_{i_j}$ is hit by the observations $X^n$. Clearly the posterior distribution of the support of the vector $u$ is uniformly distributed over
\begin{align*}
    \{i_1, \cdots, i_k\} \cup \left\{ S\subseteq [M]\backslash \{i_1,\cdots,i_k\}: |S| = \frac{M}{100} - k \right\}. 
\end{align*}
On one hand, if the estimator $x^{m+n}$ hits a set $A_j$ with $j\notin \{i_1, \cdots, i_k\}$, the probability that $x^{n+m}$ does not belong to the support of $f$ is at least
\begin{align*}
    1 - \frac{M/100-k}{M-k} = \frac{0.99M}{M-k} \ge 0.99. 
\end{align*}
On the other hand, if the estimator $x^{m+n}$ never hits a set $A_j$ with $j\notin \{i_1, \cdots, i_k\}$, then $x^{m+n}$ fall into at most $U_n$ sets in $\{A_i\}_{i=1}^M$, where $U_n$ is defined in a similar way as $U_{n+m}$. Similar to the previous computation, we have
\begin{align*}
    \bE[U_n] = \frac{\sqrt{n}}{100}\left[1 - \left(1-\frac{1}{n}\right)^n\right], \qquad \var(U_n) = O(\sqrt{n}). 
\end{align*}
By Chebyshev's inequality, for $C_0>0$ large enough, the probability that $x^{n+m}$ violates the second event in the loss function is at least $0.99$. Now a combination of the above two cases implies that $r_{\text{B}}(\calP^{\otimes n},\mu)\ge 0.9$, as desired.

\section{Proof of main lemmas}\label{appendix:lemma}
\subsection{Proof of Lemma \ref{lemma:self-concordance}}
Recall that the log-partition function $A(\theta)$ is defined as
\begin{align*}
A(\theta) = \log\int_{\calX} \exp(\theta^\top T(x))\mathrm{d}\mu(x), 
\end{align*}
for any vector $\lambda\in \bR^d$ with $\theta + (\nabla^2 A(\theta))^{-1/2}\lambda\in \Theta$, we have
\begin{align}\label{eq:MGF}
&\bE_\theta \left[\exp\left(\lambda^\top (\nabla^2 A(\theta))^{-1/2} (T(X) - \nabla A(\theta)) \right)\right] \nonumber\\
&= \int_{\calX} \exp\left( (\theta + (\nabla^2 A(\theta))^{-1/2}\lambda)^\top T(x) - A(\theta) - [(\nabla^2 A(\theta))^{-1/2}\lambda]^\top \nabla A(\theta) \right) \mathrm{d}\mu(x) \nonumber\\
&= \exp\left( A(\theta + (\nabla^2 A(\theta))^{-1/2}\lambda) - A(\theta) -  [(\nabla^2 A(\theta))^{-1/2}\lambda]^\top \nabla A(\theta)\right). 
\end{align}
It remains to show that when $\|\lambda\|_2$ is sufficiently small, we always have $\theta + (\nabla^2 A(\theta))^{-1/2}\lambda\in \Theta$, and the exponent of \eqref{eq:MGF} is uniformly bounded from above over $\theta\in \Theta$ and $\lambda\in \bR^d$. Then the existence of uniformly bounded MGF around zero implies a uniformly bounded moment of any order.

The result of \cite[Theorem 4.1.6]{nesterov2003introductory} shows that for a self-concordant and convex function $f$, we have $\nabla^2 f(y) \preceq 4\nabla^2 f(x)$ whenever $(y-x)^\top \nabla^2 f(x)(y-x) \le M^2/16$. Consequently, for $\|\lambda\|_2\le M/4$, a Taylor expansion with a Lagrange remainder gives 
\begin{align}\label{eq:taylor}
&A(\theta + (\nabla^2 A(\theta))^{-1/2}\lambda) - A(\theta) -  [(\nabla^2 A(\theta))^{-1/2}\lambda]^\top \nabla A(\theta) \nonumber\\
&= \frac{1}{2} \lambda^\top (\nabla^2 A(\theta))^{-1/2} \cdot \nabla^2A(\xi) \cdot (\nabla^2 A(\theta))^{-1/2}\lambda, 
\end{align}
where $\xi$ lies on the line segment between $\theta$ and $\theta + (\nabla^2 A(\theta))^{-1/2}\lambda$. Consequently, we have
\begin{align*}
(\xi - \theta)^\top \cdot \nabla^2 A(\theta)\cdot (\xi - \theta) \le \lambda^\top (\nabla^2 A(\theta))^{-1/2} \cdot \nabla^2A(\theta) \cdot (\nabla^2 A(\theta))^{-1/2}\lambda = \|\lambda\|_2^2 \le \frac{M^2}{16}, 
\end{align*}
and therefore $\nabla^2 A(\xi) \preceq 4\nabla^2 A(\theta)$. Plugging it back into \eqref{eq:taylor} establishes the boundedness of \eqref{eq:MGF}, as well as the finiteness of $A(\theta + (\nabla^2 A(\theta))^{-1/2}\lambda) $, or equivalently, $\theta + (\nabla^2 A(\theta))^{-1/2}\lambda\in \Theta$. 

\subsection{Proof of Lemma \ref{lemma:gaussian_covariance}}
For $X_1,\cdots,X_n\sim \calN(0,\Sigma)$ and $n\ge d$, it is well known that the empirical covariance $\widehat{\Sigma}_n$ follows a Wishart distribution $\calW_d(\Sigma/n, n)$, where the density of $\calW_d(V,n)$ is given by
\begin{align*}
f_{V,n}(X) = \frac{\det(X)^{(n-d-1)/2}\exp(-\mathsf{Tr}(V^{-1}X)/2)}{2^{nd/2}\det(V)^{n/2}\Gamma_d(n/2)}, 
\end{align*}
where $\Gamma_d(x) = \pi^{d(d-1)/4}\prod_{i=1}^d \Gamma(x - (i-1)/2)$ is the multivariate Gamma function, and $\Gamma(x) = \int_0^\infty t^{x-1}e^{-t}dt$ is the usual Gamma function. After some algebra, we have
\begin{align}\label{eq:Wishart_KL}
& D_{\text{KL}}(\calW_d(\Sigma/n,n) \| \calW_d(\Sigma/(n+m),n+m) ) \nonumber \\
&= \frac{d}{2}\left(m - (n+m)\log\left(1+\frac{m}{n}\right) \right) + \log\frac{\Gamma_d((n+m)/2)}{\Gamma_d(n/2)} - \frac{m}{2}\psi_d\left(\frac{n}{2}\right), 
\end{align}
where $\psi_d(x) = \frac{d}{dx}[\log \Gamma_d(x)]$ is the multivariate digamma function. Note that the above KL divergence \eqref{eq:Wishart_KL} does not depend on $\Sigma$; we denote it by $f(n,m,d)$. 

By \eqref{eq:TV_KL_chi}, it suffices to establish an upper bound of $f(n,m,d)$. Applying infinite Taylor series to $\log \Gamma_d(x)$ at $x=n/2$ yields 
\begin{align*}
\log\Gamma_{d}\left(\frac{n+m}{2}\right) &= \log\Gamma_d\left(\frac{n}{2}\right) + \frac{m}{2}\psi_d\left(\frac{n}{2}\right) + \sum_{t=2}^\infty \frac{1}{t!}\left(\frac{m}{2}\right)^t\psi_d^{(t-1)}\left(\frac{n}{2}\right) \\
&=\log\Gamma_d\left(\frac{n}{2}\right) + \frac{m}{2}\psi_d\left(\frac{n}{2}\right) + \sum_{t=2}^\infty \frac{1}{t!}\left(\frac{m}{2}\right)^t \sum_{k=1}^d \psi^{(t-1)}\left(\frac{n-k+1}{2}\right),
\end{align*}
where $\psi^{(t-1)}(x)=\frac{d^t}{dx^t}[\log \Gamma(x)]$ is the polygamma function. For any $t\ge 2$ and $x\ge 1$, the following inequality holds for the polygamma function \cite[Equation 6.4.10]{abramowitz1964handbook}: 
\begin{align*}
\left| \psi^{(t-1)}(x) - (-1)^t\frac{(t-2)!}{x^{t-1}} \right| \le \frac{(t-1)!}{x^t}. 
\end{align*}

As a result, for the following modification
\begin{align}\label{eq:g_nmd}
g(n,m,d) \triangleq \frac{d}{2}\left(m - (n+m)\log\left(1+\frac{m}{n}\right) \right) + m\sum_{k=1}^d\sum_{t=2}^\infty \frac{(-1)^t}{2t(t-1)}\left(\frac{m}{n-k+1}\right)^{t-1}, 
\end{align}
it holds that
\begin{align}\label{eq:difference_fg}
|f(n,m,d) - g(n,m,d)| &\le \sum_{k=1}^d\sum_{t=2}^\infty \frac{1}{t!}\left(\frac{m}{2}\right)^t\cdot \frac{(t-1)!}{[(n-k+1)/2]^{t}} \nonumber \\
&=\sum_{k=1}^d \sum_{t=2}^\infty \frac{1}{t}\left(\frac{m}{n-k+1}\right)^t \nonumber \\
&\le d\cdot \sum_{t=2}^\infty \frac{1}{2}\left(\frac{2m}{n}\right)^t \le \frac{4dm^2}{n^2}, 
\end{align}
where we have used the assumption $n\ge 4\max\{m,d\}$. 

Next we establish an upper bound of $g(n,m,d)$. Using the identity
\begin{align*}
h(x) \triangleq \frac{(1+x)\log(1+x) - x}{x} = \sum_{t=2}^\infty \frac{(-1)^t x^{t-1}}{t(t-1)}
\end{align*}
and some algebra, we have
\begin{align*}
g(n,m,d) = -\frac{dm}{2}h\left(\frac{m}{n} \right) + \frac{m}{2}\sum_{k=1}^d h\left(\frac{m}{n-k+1}\right) = \frac{m}{2}\sum_{k=1}^d \left[ h\left(\frac{m}{n-k+1}\right) - h\left(\frac{m}{n}\right) \right]. 
\end{align*}
Note that for $x\in [0,1]$, we have
\begin{align*}
h'(x) = \frac{x-\log(1+x)}{x^2} \in \left[0, \frac{1}{2}\right], 
\end{align*}
we conclude that
\begin{align}\label{eq:upper_g}
g(n,m,d) \le \frac{m}{2}\sum_{k=1}^d \frac{1}{2}\left[\frac{m}{n-k+1} - \frac{m}{n}\right] \le \frac{md}{2}\cdot \frac{2md}{n^2} = \frac{m^2d^2}{n^2}. 
\end{align}

Finally, by \eqref{eq:difference_fg}, \eqref{eq:upper_g} and \eqref{eq:TV_KL_chi}, we conclude that
\begin{align*}
\|\calL(\widehat{\Sigma}_n) - \calL(\widehat{\Sigma}_{n+m}) \|_{\text{TV}} \le \sqrt{\frac{f(n,m,d)}{2}} \le \frac{2md}{n}. 
\end{align*}

For the claim that $S_{n+m}$ follows the uniform distribution on the set $A = \{U\in \bR^{d\times (n+m)}: UU^\top = I_d \}$, we need the following auxiliary definitions and results. A topological group is a group $(G,+)$ with a topology such that the operation $+: G\times G\to G$ is continuous. A (right) group action of $G$ on $X$ is a function $\phi: X\times G\to X$ such that $\phi( \phi(x,g), g') = \phi(x,gg')$ and $\phi(x,e)=x$, where $e$ is the identity element of $G$. A group action is called \emph{transitive} if for every $x,x'\in A$, there exists some $g\in G$ such that $\phi(x,g)=x'$. A group action is called \emph{proper} if for any compact $K\subseteq X$ and $x\in X$, the map $\phi_x: G\to X$ with $g\mapsto \phi(x,g)$ satisfies that $\phi_x^{-1}(K)\subseteq G$ is compact. The following lemma is useful. 
\begin{lemma}[Chapter 14, Theorem 25 of \cite{royden1988real}]\label{lemma:haar_measure}
Let $G$ be a locally compact group acting transitively and properly on a locally compact Hausdorff space $X$. Then there is a unique (up to multiplicative factors) Baire measure on $X$ which is invariant under the action of $G$. 
\end{lemma}

Now for the claimed result, it is easy to verify that (a proper version of) $S_{n+m}$ always takes value in $A$, assuming $n+m\ge d$. To show that $S_{n+m}$ is uniform on $A$, note that for any orthogonal matrix $V\in \bR^{(n+m)\times (n+m)}$, it is easy to verify
\begin{align*}
[X_1,\cdots,X_{n+m}] \overset{d}{=} [X_1,\cdots,X_{n+m}]V. 
\end{align*}
Denote the RHS by $[Y_1,\cdots,Y_{n+m}]$, we also have $\widehat{\Sigma}_{n+m}(Y^{n+m}) = \widehat{\Sigma}_{n+m}(X^{n+m})$. Consequently, 
\begin{align*}
S_{n+m}(X^{m+n})\cdot V &= [(n+m)\widehat{\Sigma}_{n+m}(X^{n+m})]^{-1/2}[Y_1, Y_2, \cdots, Y_{n+m}] \\
&= [(n+m)\widehat{\Sigma}_{n+m}(Y^{n+m})]^{-1/2}[Y_1, Y_2, \cdots, Y_{n+m}] \\
&\overset{d}{=} [(n+m)\widehat{\Sigma}_{n+m}(X^{n+m})]^{-1/2}[X_1, X_2, \cdots, X_{n+m}] = S_{n+m}(X^{m+n}), 
\end{align*}
meaning that the distribution of $S_{n+m}$ is invariant with right multiplication of an orthogonal matrix. Let $G$ be the orthogonal group $\mathsf{O}(n+m)$, the map $\phi: A\times G\to A$ with $\phi(U,V) = UV$ is a group action of $G$ on $A$. This action is transitive as for any $U, U'\in A$, we could add more rows to $U, U'$ to obtain $\widetilde{U}, \widetilde{U}' \in G$, and then $V = \widetilde{U}^\top \widetilde{U}'\in G$ maps $U$ to $U'$. This action is also proper as $G$ itself is compact. Hence, Lemma \ref{lemma:haar_measure} below shows that $S_{n+m}$ is uniformly distributed on $A$. 

\subsection{Proof of Lemma \ref{lemma:gaussian_mean_covariance}}
The proof of Lemma \ref{lemma:gaussian_mean_covariance} is a simple consequence of several known results. First, Basu's theorem claims that $\overline{X}_n$ and $\widehat{\Sigma}_n$ are independent. Second, the computation in Example \ref{example:mean} shows
\begin{align*}
\| \calL(\overline{X}_n) - \calL(\overline{X}_{n+m}) \|_{\text{TV}} \le \frac{m\sqrt{d}}{n}. 
\end{align*}
Finally, as $\widehat{\Sigma}_n \sim \calW_d(\Sigma/(n-1), n-1)$ again follows a Wishart distribution, the proof of Lemma \ref{lemma:gaussian_covariance} shows that
\begin{align*}
\|\calL(\widehat{\Sigma}_n) - \calL(\widehat{\Sigma}_{n+m}) \|_{\text{TV}} \le \frac{2md}{n-1}. 
\end{align*}
In conclusion, we have
\begin{align*}
	\| \calL(\overline{X}_n, \widehat{\Sigma}_n) - \calL(\overline{X}_{n+m}, \widehat{\Sigma}_{n+m}) \|_{\text{\rm TV}} &\le \| \calL(\overline{X}_n) - \calL(\overline{X}_{n+m}) \|_{\text{TV}}  + \|\calL(\widehat{\Sigma}_n) - \calL(\widehat{\Sigma}_{n+m}) \|_{\text{TV}} \\
	&\le \frac{3md}{n-1}. 
\end{align*}

For the distribution of $S_{n+m}$, it is clear that (a proper version of) $S_{n+m}$ always takes value in $A$, and we show that the distribution of $S_{n+m}$ is invariant with proper group actions in $A$. Consider the set
\begin{align*}
G = \{V \in \bR^{(n+m)\times (n+m)}: VV^\top = I_{n+m}, V{\bf 1} = {\bf 1} \}
\end{align*}
with the usual matrix multiplication. We show that $G$ is a group: clearly $I_{n+m}\in G$; for $V,V'\in G$, it is clear that $VV'{\bf 1} = V{\bf 1}={\bf 1}$ and therefore $VV'\in G$; for $V\in G$, it holds that $V^{-1}{\bf 1} = V^{-1}V{\bf 1}={\bf 1}$ and thus $V^{-1}\in G$. Next we show that the action $\phi: A\times G\to A$ with $\phi(U,V) = UV$ is a group action on $A$, and it suffices to show that $UV\in A$. This is true as $(UV)(UV)^\top = UVV^\top U^\top = UU^\top = I_d$, and $UV{\bf 1} = U{\bf 1} = {\bf 0}$. This group action is also transitive: for any $U,U'\in A$, we may properly add rows to them and obtain $\widetilde{U}, \widetilde{U}'\in \mathsf{O}(n+m)$, where one of the added rows is a scalar multiple of ${\bf 1}^\top$, which is feasible as $U{\bf 1}={\bf 0}$. Consequently, the matrix $V = \widetilde{U}^{-1}\widetilde{U}' \in \mathsf{O}(n+m)$ maps $U$ to $U'$, and also $\bf 1^\top$ to $\bf 1^\top$; hence $V\in G$. Finally, we show that any group action of $G$ on $S_{m+n}$ does not change the distribution of $S_{m+n}$. To see this, for any $V\in G$ we have
\begin{align*}
[X_1,\cdots,X_{n+m}] &\overset{d}{=} [X_1,\cdots,X_{n+m}]V, \\
\overline{X}_{n+m}([X_1,\cdots,X_{n+m}]V) &= \overline{X}_{n+m}([X_1,\cdots,X_{n+m}]), \\
\widehat{\Sigma}_{n+m}([X_1,\cdots,X_{n+m}]V) &= \widehat{\Sigma}_{n+m}([X_1,\cdots,X_{n+m}]). 
\end{align*}
Therefore, following the same arguments as in Example \ref{example:gaussian_covariance} we arrive at the desired invariance, and the uniform distribution of $S_{m+n}$ is a direct consequence of Lemma \ref{lemma:haar_measure}. 

\subsection{Proof of Lemma \ref{lemma:shuffle_chi}}
For any subset $S\subseteq [n+m]$ with $|S| = m$, let $P_S$ be the distribution of $(Z_1,\cdots,Z_{n+m})$ when the samples $(Y_1,\cdots,Y_m)$ are placed in the index set $S$ of the pool $(Z_1,\cdots,Z_{n+m})$. Then it is clear that $P_{\text{mix}} = \bE[P_S]$, with $S$ uniformly distributed on all size-$m$ subsets of $[n+m]$. To compute the $\chi^2$-divergence where the first distribution is a mixture, the following identity holds \cite{ingster2012nonparametric}: 
$$
\chi^2\left(P_{\text{\rm mix}}, P^{\otimes (n+m)} \right)  = \bE_{S,S'}\left[ \int\frac{\mathrm{d}P_S\mathrm{d}P_{S'}}{\mathrm{d}P^{\otimes (n+m)}} \right] - 1, 
$$
where $S'$ is an independent copy of $S$. By the independence assumption, we have
$$
\frac{\mathrm{d}P_S}{\mathrm{d}P^{\otimes (n+m)}}(z_1,\cdots,z_{n+m}) = \prod_{i\in S} \frac{\mathrm{d}Q}{\mathrm{d}P}(z_i), 
$$
and consequently
\begin{align*}
\int\frac{\mathrm{d}P_S\mathrm{d}P_{S'}}{\mathrm{d}P^{\otimes (n+m)}} &= \bE_{P}\left[\prod_{i\in S} \frac{\mathrm{d}Q}{\mathrm{d}P}(z_i) \prod_{i\in S'} \frac{\mathrm{d}Q}{\mathrm{d}P}(z_i)\right] \\
&= \prod_{i\in S\cap S'} \bE_P\left[\left(\frac{\mathrm{d}Q}{\mathrm{d}P}(z_i)\right)^2\right] \cdot \prod_{i\in S\Delta S'} \bE_P\left[\frac{\mathrm{d}Q}{\mathrm{d}P}(z_i) \right] \\
&= (1 + \chi^2(Q,P))^{|S\cap S'|}. 
\end{align*}

It remains to upper bound the expectation with respect to the random variable $|S\cap S'|$. Note that $|S\cap S'|$ follows the hypergeometric distribution with parameter $(n+m,m,m)$, which corresponds to sampling \emph{without replacement}. The counterpart for sampling \emph{with replacement} corresponds to a Binomial distribution $\mathsf{B}(m, \frac{m}{n+m})$, and the following lemma shows that the latter dominates the former in terms of the convex order: 
\begin{lemma}\cite[Theorem 4]{hoeffding1963probability}\label{lem:convex1}
	Let the population be $\calC= \{c_1,\cdots,c_N\}$. Let $X_1,\cdots,X_n$ denote random samples without replacement from $\calC$ and $Y_1,\cdots,Y_n$ denote random samples with replacement. If $f(\cdot)$ is convex, then
	$$
	\bE \left[ f\left(\sum_{i=1}^n X_i\right) \right] \le \bE \left[ f\left(\sum_{i=1}^n Y_i\right) \right]. 
	$$
\end{lemma}

Applying Lemma \ref{lem:convex1} to the convex function $x\mapsto (1+\chi^2(Q,P))^x$ yields
\begin{align*}
\bE_{S,S'}[ (1 + \chi^2(Q,P))^{|S\cap S'|} ] &\le \bE[ (1 + \chi^2(Q,P))^{\mathsf{B}(m,m/(n+m))} ]  \\
&=\left( \bE[(1 + \chi^2(Q,P))^{\text{Bern}(m/(n+m))}] \right)^m \\
&= \left(1 + \frac{m}{n+m} \chi^2(Q,P)\right)^m, 
\end{align*}
as desired. 

\subsection{Proof of Lemma \ref{lemma:modified_chi2}}
Note that in the proof of Theorem \ref{thm:shuffling_general}, we have
\begin{align*}
\| P_{\text{mix}}(X^{n/2}) - P^{\otimes (n/2+m)} \|_{\text{TV}} \le \sqrt{\frac{m^2}{n}\chi^2(\widehat{P}_n(X^{n/2}),P)}. 
\end{align*}
Since the TV distance is always upper bounded by one, the following upper bound is also true: 
\begin{align*}
\| P_{\text{mix}}(X^{n/2}) - P^{\otimes (n/2+m)} \|_{\text{TV}} \le \sqrt{\frac{m^2}{n}\cdot\left( \chi^2(\widehat{P}_n(X^{n/2}),P) \wedge n\right)}. 
\end{align*}
Consequently, taking the expectation over $X^{n/2}$ leads to the claim. 

\subsection{Proof of Lemma \ref{lmm:sparse_gaussian_shuffle}}
First we note that
\begin{align*}
\chi^2\left( \calN(\widehat{\theta}_{n,j}, 1), \calN(\theta_j,1) \right)= \exp\left( (\widehat{\theta}_{n,j} - \theta_j)^2 \right) - 1,
\end{align*}
By the triangle inequality, 
\begin{align}\label{eq:error_mean}
\bE|\widehat{\theta}_{n,j} - \theta_j| &\le \bE\left|\widehat{\theta}_{n,j} - \frac{1}{n}\sum_{i=1}^n X_{i,j} \right| + \bE\left|\frac{1}{n}\sum_{i=1}^n X_{i,j} - \theta_j \right| \nonumber \\
&\le \sqrt{\frac{C\log n}{n}} + \sqrt{\frac{2}{n\pi}} = O\left(\sqrt{\frac{\log n}{n}}\right). 
\end{align}
Moreover, it is straightforward to verify that $|\widehat{\theta}_{n,j} - \theta_j|$ is $(1/n)$-Lipschitz with respect to $(X_{1,j},\cdots,X_{n,j})$. Therefore, the Gaussian Lipschitz concentration (see, e.g. \cite[Theorem 10.17]{Boucheron--Lugosi--Massart2013}) gives that
\begin{align}\label{eq:error_concentration}
\bP\left(|\widehat{\theta}_{n,j} - \theta_j| \ge \bE|\widehat{\theta}_{n,j} - \theta_j| + t \right) \le \exp\left(-\frac{nt^2}{2}\right)
\end{align}
for any $t\ge 0$. Hence, combining \eqref{eq:error_mean} and \eqref{eq:error_concentration}, we conclude that $|\widehat{\theta}_{n,j} - \theta_j| = O(\sqrt{\log(nd)/n})$ holds with probability at least $1-(nd)^{-2}$, and therefore
\begin{align*}
\bE\left[ \chi^2\left( \calN(\widehat{\theta}_{n,j}, 1), \calN(\theta_j,1) \right) \wedge n \right] \le 2\cdot \bE[(\widehat{\theta}_{n,j} - \theta_j)^2] + \frac{1}{nd}, 
\end{align*}
where we have used that $e^x\le 1+2x$ whenever $x\in [0,1]$. Summing over $j\in [d]$ and using the property of the soft-thresholding estimator $\sup_{\theta\in\Theta}\bE[\| \widehat{\theta}_n - \theta \|_2^2]=O(s\log d/n)$ gives the claimed result.

\subsection{Proof of Lemma \ref{lmm:shuffle_poisson}}
For the first claim, note that for Poisson models, we have
\begin{align*}
\chi^2\left( \mathsf{Poi}(\widehat{\lambda}_{n,j}) , \mathsf{Poi}(\lambda_j) \right) = \exp\left(\frac{(\widehat{\lambda}_{n,j} - \lambda_j)^2}{\lambda_j} \right) - 1.
\end{align*}
Consequently, for $\lambda_j = 1/(n^2\log n)$, we have
\begin{align*}
\bE\left[ \chi^2\left( \mathsf{Poi}(\widehat{\lambda}_{n,j}) , \mathsf{Poi}(\lambda_j) \right) \wedge n \right] &\ge \left[\left( \exp\left(\frac{(1/n - \lambda_j)^2}{\lambda_j} \right) - 1 \right) \wedge n\right]\cdot \bP(\widehat{\lambda}_{n,j} = 1/n) \\
&= \Omega(n)\cdot e^{-n\lambda_j}n\lambda_j = \Omega\left(\frac{1}{\log n} \right) \gg \frac{1}{n}, 
\end{align*}
establishing the first claim. 

For the second claim, note that conditioning on $X^{n/2}$, 
\begin{align*}
D_{\text{\rm KL}}\left(\calL(\widehat{T}_{n/2+m}) \| \calL(T_{n/2+m}) \right) &= \sum_{j=1}^d D_{\text{\rm KL}}\left(\Poi(n\lambda_j/2 + m\widehat{\lambda}_{n,j}) \| \Poi((n/2+m)\lambda_j)\right) \\
&\le \sum_{j=1}^d \frac{m^2(\widehat{\lambda}_{n,j} - \lambda_j)^2}{(n/2+m)\lambda_j}. 
\end{align*}
where we have used $D_{\text{\rm KL}}(\Poi(\lambda_1) \| \Poi(\lambda_2)) = \lambda_2 - \lambda_1 + \lambda_1\log(\lambda_1/\lambda_2)\le (\lambda_1-\lambda_2)^2/\lambda_2$. Hence, 
\begin{align*}
\bE_{X^{n/2}}\left[D_{\text{\rm KL}}\left(\calL(\widehat{T}_{n/2+m}) \| \calL(T_{n/2+m}) \right)  \right] \le \frac{dm^2}{(n/2+m)n/2} \le \frac{4dm^2}{n^2}, 
\end{align*}
and the TV distance satisfies 
\begin{align*}
\| \calL(\widehat{X}^{n+m}) - \calL(X^{n+m}) \|_{\text{\rm TV}} &= \bE_{X^{n/2}}\left[ \| \calL(\widehat{T}_{n/2+m}) - \calL(T_{n/2+m})\|_{\text{\rm TV}}\right]  \\
&\le \bE_{X^{n/2}} \left[\sqrt{\frac{1}{2} D_{\text{\rm KL}}\left(\calL(\widehat{T}_{n/2+m}) \| \calL(T_{n/2+m}) \right)} \right] \\
&\le \sqrt{\frac{1}{2}\bE_{X^{n/2}}\left[D_{\text{\rm KL}}\left(\calL(\widehat{T}_{n/2+m}) \| \calL(T_{n/2+m}) \right)  \right]} \\
&\le \frac{m\sqrt{2d}}{n}. 
\end{align*}

\subsection{Proof of Lemma \ref{lmm:lower_poissonization}}
The upper bound result is easy. The $m=O(n\varepsilon/\sqrt{k})$ upper bound is a consequence of the general product Poisson model considered in Example \ref{example:poisson}. For the $m=O(\sqrt{n}\varepsilon)$ upper bound, we consider the sufficient statistic $T_n=\sum_{i=1}^n X_i\sim \prod_{j=1}^k \mathsf{Poi}(np_j)$, and simply apply the sufficiency-based algorithm to $\widehat{T}_{n+m} = T_n$. Since
\begin{align*}
D_{\text{\rm KL}}(\calL(T_{n+m}) \| \calL(T_n)) = \sum_{j=1}^k D_{\text{\rm KL}}(\mathsf{Poi}((n+m)p_j) \| \mathsf{Poi}(np_j)) \le \sum_{j=1}^k \frac{(mp_j)^2}{np_j} = \frac{m^2}{n}, 
\end{align*}
where the last identity crucially makes use of the identity $\sum_{j=1}^k p_j = 1$, this procedure works. 

Next we show that sample amplification is impossible when $m=\omega(n\varepsilon/\sqrt{k})$ and $k=O(n)$. Note that this implies that $m=\omega(\max\{n\varepsilon/\sqrt{k}, \sqrt{n}\varepsilon\} )$ is impossible in general, for the $k>n$ case is always not easier than the $k=n$ case. To prove the above claim, w.l.o.g. we assume that $k=2k_0$ is even, and consider the following parametric submodel: 
\begin{align*}
P_{\theta} = \prod_{j=1}^{k_0} \left[\mathsf{Poi}\left(\frac{1}{k} + \theta_j \right) \times \mathsf{Poi}\left(\frac{1}{k} - \theta_j \right)\right], \qquad \theta \in \Theta \triangleq \left[-\frac{1}{k}, \frac{1}{k}\right]^{k_0}.
\end{align*}
Clearly $P_\theta$ is a parametric submodel, by setting $p_{2j-1} = 1/k + \theta_j$ and $p_{2j} = 1/k-\theta_j$ in the original model. This submodel is a product model, thus we could apply the result of Theorem \ref{thm:lower_product} after we have verified Assumption \ref{assump.hellinger}. Note that when $\theta, \theta'$ arbitrarily vary in $[-1/k, 1/k]$, the range of
\begin{align*}
H^2\left(\mathsf{Poi}\left(\frac{1}{k} + \theta \right) \times \mathsf{Poi}\left(\frac{1}{k} - \theta \right), \mathsf{Poi}\left(\frac{1}{k} + \theta' \right) \times \mathsf{Poi}\left(\frac{1}{k} - \theta' \right)\right)
\end{align*}
is $[0, \Theta(1/k)]$, so Assumption \ref{assump.hellinger} is fulfilled when $k\le cn$ for a small constant $c>0$. Consequently, Theorem \ref{thm:lower_product} establishes the desired bound $m=O(n\varepsilon/\sqrt{k})$.

\subsection{Proof of Lemma \ref{lemma:gaussian_prob}}
For the first inequality, note that for a large $C_0>0$, both $Z_1\le C_0$ and $\max_{2\le j\le d} Z_j \ge \sqrt{2\log d} - C_0$ hold with probability at least $0.99$. When both events hold, we have
\begin{align*}
\exp(t(t+Z_1)) &\le \exp(t(t+C_0)), \\
\sum_{j=2}^d \exp(tZ_j) &\ge \exp(t(\sqrt{2\log d} - C_0)). 
\end{align*}
Consequently, for $C=3C_0$ with $C_0>0$ large enough, we have
\begin{align*}
p_d(\sqrt{2\log d} - C) &\le 0.01 + \frac{\exp((\sqrt{2\log d} - 3C_0)(\sqrt{2\log d}-2C_0))}{\exp((\sqrt{2\log d}-3C_0)(\sqrt{2\log d} - C_0) )} \\
&= 0.01 + \exp\left(-C_0 (\sqrt{2\log d} - 3C_0) \right) < 0.1. 
\end{align*}

For the second inequality, note that Jensen's inequality yields that
\begin{align*}
p_d(t) &\ge \bE_{Z_1}\left[\frac{\exp(t(t+Z_1))}{ \exp(t(t+Z_1)) + \sum_{j=2}^d\bE[\exp(tZ_j)] }\right] \\
&\ge \bE_{Z_1}\left[\frac{\exp(t(t+Z_1))}{ \exp(t(t+Z_1)) + d\exp(t^2/2)}\right]. 
\end{align*}
Again, with probability at least $0.99$ we have $Z_1\ge -C_0$, and therefore for $t=\sqrt{2\log d} + 2C_0$, 
\begin{align*}
p_d(t) &\ge 0.99\times \frac{\exp(t(t-C_0))}{ \exp(t(t-C_0)) + d\exp(t^2/2)} \\
&= 0.99\times \frac{1}{1+\exp(t^2/2 + (t-2C_0)^2/2 - t(t-C_0) ) }\\
&=0.99\times \frac{1}{1+\exp(-C_0\cdot \sqrt{2\log d})} > 0.9, 
\end{align*}
for $C_0>0$ large enough. 

For the last claim, consider any $s<d/2$.  Let $d_0 \triangleq \lfloor d/s\rfloor\ge 2$, and consider the product prior $\mu^{\otimes s}$ to $s$ blocks each of dimension $d_0$. In other words, we set exactly one of the first $d_0$ coordinates of the mean vector to $t$ uniformly at random, and do the same for the next $d_0$ coordinates, and so on. Clearly the resulting mean vector is always $s$-sparse. Writing $\theta=(\theta_1,\cdots,\theta_{s})$ with each $\theta_i\in \bR^{d_0}$, let the loss function be
\begin{align*}
L(\theta, \widehat{\theta}) = \1\left(\sum_{j=1}^s \1(\theta_j \neq \widehat{\theta}_j) \ge N \right),
\end{align*}
with an integer $N$ to be specified later. Then in each block, we reduce to the case $s=1$, and the error probability for this block is $1-p_{d_0}(\sqrt{n}\cdot t)$ for sample size $n$. Moreover, the errors in different blocks are independent. Consequently, 
\begin{align}\label{eq:binomial_diff}
 &r_{\text{\rm B}}(\calP,n,\mu,L) - r_{\text{\rm B}}(\calP,n+m,\mu,L) \nonumber \\
 &= \bP\left(\mathsf{B}(s,1-p_{d_0}(\sqrt{n}\cdot t))\ge N \right) - \bP\left(\mathsf{B}(s,1-p_{d_0}(\sqrt{n+m}\cdot t))\ge N \right). 
\end{align}
Finally, again by the properties of $p_d(\cdot)$ summarized in Lemma \ref{lemma:gaussian_prob} and the pigeonhole principle, for $m=\lceil cn\varepsilon/\sqrt{s\log d_0}\rceil$, we could always find some $t>0$ such that $p_{d_0}(\sqrt{n+m}\cdot t) - p_{d_0}(\sqrt{n}\cdot t) = \Omega_c(\varepsilon/\sqrt{s})$, with both quantities in $[0.1,0.9]$. Now choosing 
\begin{align*}
    N = \left\lceil\frac{s(1-p_{d_0}(\sqrt{n}\cdot t))}{2} + \frac{s(1-p_{d_0}(\sqrt{n+m}\cdot t))}{2}\right\rceil
\end{align*}
with the above choice of $t$, the Bayes risk difference will be lower bounded by $\Omega_c(\varepsilon)$ by a similar argument to the proof of Theorem \ref{thm:lower_bound}. Now by \eqref{eq:binomial_diff} and Lemma \ref{lemma:general_lower_bound}, we see that $n=\Omega(\sqrt{s\log(d/s)}/\varepsilon)$ and $m=O(n\varepsilon/\sqrt{s\log(d/s)})$ are necessary for sample amplification.

\subsection{Proof of Lemma \ref{lmm:invariant_estimator}}
The following results are the key to the proof of Lemma \ref{lmm:invariant_estimator}, which we will assume for now and prove in the subsequent sections. 
\begin{lemma}\label{lemma:stein_loss}
For $D_n = \text{\rm diag}(\lambda_1,\cdots,\lambda_d)$, the following identities hold for the estimator $\widehat{\Sigma}_n$: 
\begin{align}
\bE[\ell(\Sigma, \widehat{\Sigma}_n)] &= \sum_{j=1}^d \left[(n+d+1-2j)\lambda_j - \log \lambda_j - \bE\left[ \log \chi_{n-j+1}^2 \right] \right] - d,  \label{eq:stein_mean} \\
\var(\ell(\Sigma, \widehat{\Sigma}_n)) &= \sum_{j=1}^d\left[ 2(n+d+1-2j)\lambda_j^2 -4\lambda_j + \psi'\left(\frac{n+1-j}{2}\right)\right], \label{eq:stein_variance}
\end{align}
where $\chi_m^2$ denotes the chi-squared distribution with $m$ degrees of freedom, and $\psi'(x)$ is the polygamma function of order $1$. In particular, when $n\ge 2d$, the following inequalities hold: 
\begin{align}
\left| \bE[\ell(\Sigma, \widehat{\Sigma}_n)] - \left(g(n+1-d,d) + \sum_{j=1}^d h((n+d+1-2j)\lambda_j) \right) \right| \le \frac{5d}{n}, \label{eq:stein_mean_bound} \\
\var( \ell(\Sigma, \widehat{\Sigma}_n) ) \le \frac{16d^2}{n^2} + \sum_{j=1}^d \frac{4((n+d+1-2j)\lambda_j - 1)^2}{n}, \label{eq:stein_variance_bound} 
\end{align}
where the functions $g$ and $h$ are given by \eqref{eq:g_function} and 
\begin{align}
h(u) &\triangleq u - \log u - 1. \label{eq:h_function}
\end{align}
\end{lemma}

\begin{lemma}\label{lemma:functions_g_h}
For the function $h$ in \eqref{eq:h_function} and $u_1,\cdots,u_d\in \bR_+$, it holds that
\begin{align}\label{eq:h_function_ineq}
\sum_{j=1}^d h(u_j) \ge \frac{1}{8} \min\left\{ \sum_{j=1}^d (u_j-1)^2, \sqrt{\sum_{j=1}^d (u_j-1)^2} \right\}. 
\end{align}
\end{lemma}

Returning to the proof of Lemma \ref{lmm:invariant_estimator}, the first claim is a direct application of \eqref{eq:stein_mean_bound} and \eqref{eq:stein_variance_bound}. For the second claim, let $D_n = \text{\rm diag}(\lambda_1,\cdots,\lambda_d)$, and
\begin{align*}
V\triangleq \sum_{j=1}^d ((n+d+1-2j)\lambda_j - 1)^2. 
\end{align*}
Then by \eqref{eq:stein_mean_bound}, \eqref{eq:stein_variance_bound} and Lemma \ref{lemma:functions_g_h}, we have
\begin{align*}
\bE[\ell(\Sigma, \widehat{\Sigma}_n)] &\ge g(n+1-d,d) + \frac{V\wedge \sqrt{V}}{8} - \frac{5d}{n}. 
\end{align*}
Meanwhile, the variance satisfies
\begin{align*}
    \sqrt{\var(\ell(\Sigma, \widehat{\Sigma}_n))} \le \frac{4d}{n} + 2\sqrt{\frac{V}{n}}. 
\end{align*}

Note that for $n,d$ larger than a constant depending only on $C_1$, we always have
\begin{align*}
\frac{d}{n} + \frac{V\wedge \sqrt{V}}{8} \ge 2C_1 \sqrt{\frac{V}{n}}. 
\end{align*}
Therefore, for $n,d=\Omega(1)$ large enough, the above inequalities implies
\begin{align*}
\ell(\Sigma, \widehat{\Sigma})  \ge g(n+1-d,d)  + C_1 \left(\sqrt{\var(\ell(\Sigma, \widehat{\Sigma}_n))} - \frac{4d}{n}\right) - \frac{6d}{n}, 
\end{align*}
establishing the second claim. 

For the last statement, note that
\begin{align*}
\frac{\partial g(u,v)}{\partial u} = \frac{\log(u+2v) + \log(v)}{2} - \log(u+v) = \frac{1}{2}\log\left( \frac{u(u+2v)}{(u+v)^2} \right) \le -\frac{v^2}{2(u+v)^2}, 
\end{align*}
the intermediate value theorem then implies that
\begin{align*}
g(n+1-d,d) - g(n+m+1-d,d) \ge m\cdot \min_{u\in [n+1-d,n+m+1-d]} \frac{d^2}{2(u+d)^2} \ge \frac{md^2}{13n^2}. 
\end{align*}

\subsection{Proof of Lemma \ref{lemma:stein_loss}}
We first recall the well-known Bartlett decomposition: for the lower triangular matrix $L_n = (L_{ij})_{i\ge j}$, the random variables $\{L_{ij}\}_{i\ge j}$ are mutually independent, with
\begin{align*}
L_{ij} \sim \calN(0, 1), \quad i>j; \qquad L_{jj}^2\sim \chi_{n+1-j}^2, \quad j\in [d]. 
\end{align*}
For $\Sigma = I_d$ and $\widehat{\Sigma}_n = L_nD_nL_n^\top$, simple algebra gives
\begin{align*}
\ell(\Sigma, \widehat{\Sigma}_n) = \sum_{j=1}^d \left( \sum_{i\ge j}  \lambda_{j} L_{ij}^2 - \log \lambda_j - \log L_{jj}^2 - 1\right). 
\end{align*}
Consequently, the identity \eqref{eq:stein_mean} follows from the above Bartlett decomposition. This identity was also obtained in \cite{james1961estimation}. 

For the identity \eqref{eq:stein_variance} on the variance, by the mutual independence we have
\begin{align}\label{eq:stein_variance_1}
\var(\ell(\Sigma, \widehat{\Sigma}_n)) = \sum_{j=1}^d \sum_{i\ge j} \lambda_j^2 \cdot \var(L_{ij}^2) + \sum_{j=1}^d \var(\log L_{jj}^2) - 2\sum_{j=1}^d \lambda_j\cdot \mathsf{Cov}(L_{jj}^2, \log L_{jj}^2). 
\end{align}
Next we evaluate each term of \eqref{eq:stein_variance_1}. Clearly $\var(L_{ij}^2) = \bE[Z^4] - \bE[Z^2]^2 = 2$ for $i>j$ and $Z\sim \calN(0,1)$. For the other random variables, we need to recall the following identity for $\chi_m^2$: 
\begin{align}\label{eq:chi_squared_MGF}
\Lambda(t) \triangleq \log \bE[(\chi_m^2)^t] = t\log 2 + \log \Gamma\left(\frac{m}{2} + t\right) - \log \Gamma\left(\frac{m}{2}\right). 
\end{align}
Based on \eqref{eq:chi_squared_MGF}, we have
\begin{align*}
\var(L_{jj}^2) = (n+1-j)(n+3-j) - (n+1-j)^2 = 2(n+1-j). 
\end{align*}
Moreover, differentiating $\Lambda(t)$ at $t=0$ gives
\begin{align}
\bE[\log \chi_m^2] &= \Lambda'(0) = \log 2 + \psi\left(\frac{m}{2}\right), \label{eq:first_derivative}\\
\bE[(\log \chi_m^2)^2] - \left( \bE[\log \chi_m^2] \right)^2 &= \Lambda''(0) = \psi'\left(\frac{m}{2}\right) \label{eq:second_derivative}.
\end{align}

Note that \eqref{eq:second_derivative} leads to
\begin{align*}
\var(\log L_{jj}^2)  = \psi'\left(\frac{n+1-j}{2} \right).
\end{align*}
Finally, differentiating $\Lambda(t)$ at $t=1$ gives that
\begin{align*}
\bE[L_{jj}^2\log L_{jj}^2] &= \bE[L_{jj}^2]\cdot \left(\log 2 + \psi\left(\frac{n+1-j}{2}+1\right)\right) \\
&= (n+1-j)\left(\log 2 + \psi\left(\frac{n+1-j}{2}+1\right)\right),
\end{align*}
and hence the identity $\psi(x+1) = \psi(x) + x^{-1}$ for the digamma function together with \eqref{eq:first_derivative} leads to
\begin{align*}
\mathsf{Cov}(L_{jj}^2, \log L_{jj}^2) = 2. 
\end{align*}
Therefore, plugging the above quantities in \eqref{eq:stein_variance_1} gives the identity \eqref{eq:stein_variance}. 

Next we prove the remaining inequalities when $n\ge 2d$. For \eqref{eq:stein_mean_bound}, note that \eqref{eq:stein_mean} gives an identity (together with \eqref{eq:first_derivative})
\begin{align*}
&\bE[\ell(\Sigma, \widehat{\Sigma}_n)] - \left( g(n+1-d,d) + \sum_{j=1}^d h((n+d+1-2j)\lambda_j) \right) \\
&= \sum_{j=1}^d \left(\log(n+d+1-2j) - \log 2 - \psi\left(\frac{n+1-j}{2}\right) \right) - g(n+1-d,d). 
\end{align*}
Since $|\psi(x)-\log(x)|\le 1/x$ for all $x\ge 1$, replacing $\psi(\cdot)$ by $\log(\cdot)$ in the above expression only incurs an absolute difference at most
\begin{align*}
\sum_{j=1}^d \frac{2}{n+1-j} \le 2\int_0^d \frac{\mathrm{d}x}{n-x} = 2\log\left(1+\frac{d}{n-d}\right) \le \frac{4d}{n}.
\end{align*}
For the remaining terms, it is not hard to verify that
\begin{align*}
g(n+1-d,d) = \int_0^d \log\frac{n+1-d+2x}{n+1-d+x}dx. 
\end{align*}
As $x\mapsto \log(n+1-d+2x) - \log(n+1-d+x)$ is increasing on $[0,\infty)$, we have
\begin{align*}
0\le g(n+1-d,d) - \sum_{x=0}^{d-1}  \log\frac{n+1-d+2x}{n+1-d+x} \le \log\frac{n+1-d+2d}{n+1-d+d} \le \frac{d}{n}. 
\end{align*}
Now \eqref{eq:stein_mean_bound} follows from a combination of the above inequalities. 

For the inequality \eqref{eq:stein_variance_bound}, we complete the square of \eqref{eq:stein_variance} to obtain
\begin{align*}
\var(\ell(\Sigma,\widehat{\Sigma}_{n})) &= \sum_{j=1}^d \frac{2((n+d+1-2j)\lambda_j - 1)^2}{n+d+1-2j} + \sum_{j=1}^d \left(\psi'\left(\frac{n+1-j}{2}\right) - \frac{2}{n+d+1-2j} \right) \\
&\le \frac{4}{n}\sum_{j=1}^d ((n+d+1-2j)\lambda_j - 1)^2 + \sum_{j=1}^d \left(\psi'\left(\frac{n+1-j}{2}\right) - \frac{2}{n+d+1-2j} \right). 
\end{align*}
To handle the second sum, note that \cite[Equation 6.4.10]{abramowitz1964handbook} gives $|\psi'(x) - x^{-1}| \le x^{-2}$ for $x\ge 1$. Therefore, the second term has an absolute value at most
\begin{align*}
\sum_{j=1}^d \left[ \frac{2(d-j)}{(n+1-j)(n+d+1-2j)} + \left( \frac{2}{n+1-j} \right)^2 \right] \le d\cdot \left[ \frac{2d}{(n/2)^2} + \left(\frac{2}{n/2}\right)^2 \right] \le \frac{16d^2}{n^2}. 
\end{align*}

\subsection{Proof of Lemma \ref{lemma:functions_g_h}}
Note that when $0<u\le 2$, Taylor expansion of $h(\cdot)$ at $u=1$ gives
\begin{align*}
h(u) \ge \frac{(u-1)^2}{2}\min_{\theta\in [0,2]} h''(\theta) = \frac{(u-1)^2}{8}. 
\end{align*}
For $u>2$, we have $u-1\ge \log_2 u > 10/7\cdot \log u$, and therefore
\begin{align*}
h(u) = u - \log u - 1 \ge \frac{3(u-1)}{10}. 
\end{align*}
Therefore, in both cases we have
\begin{align*}
h(u) \ge \frac{1}{8}\min\left\{(u-1)^2, |u-1| \right\}. 
\end{align*}

To prove \eqref{eq:h_function_ineq}, let $J \triangleq \{j\in [d]: |u_j-1|\le 1\}$, $S \triangleq \sum_{j=1}^d (u_j-1)^2$. Using the above inequality, we have
\begin{align*}
\sum_{j=1}^d h(u_j) &\ge \frac{1}{8} \sum_{j\in J} (u_j-1)^2 + \frac{1}{8}\sum_{j\notin J} |u_j - 1| \\
&\ge \frac{1}{8}\sum_{j\in J} (u_j-1)^2 + \frac{1}{8}\sqrt{ \sum_{j\notin J} (u_j - 1)^2} \\
&\ge \frac{1}{8} \min_{x\in [0,S]}\left(x + \sqrt{S-x} \right) \\
&= \frac{1}{8}\min\{S, \sqrt{S} \}, 
\end{align*}
which is precisely \eqref{eq:h_function_ineq}.

\bibliographystyle{alpha}
\bibliography{refs}

\end{document}